\documentclass[english]{smfart}

\usepackage{amsfonts, amssymb, amsmath, latexsym}
\usepackage{mathptmx}
\usepackage[english]{babel}
\usepackage[latin1]{inputenc}
\usepackage[T1]{fontenc}
\usepackage{textcomp}

\usepackage{graphicx}
\usepackage[dvips]{color}
\input xy
\xyoption{all}

\newcommand{\mr}{\mathrm}

\theoremstyle{plain}

\newtheorem{Thm}{Theorem}[section]
\newtheorem{Prop}[Thm]{Proposition}
\newtheorem{Lemma}[Thm]{Lemma}
\newtheorem{Def}[Thm]{Definition}

\newtheorem{Example}[Thm]{Example}
\newtheorem{Question}[Thm]{Question}
\newtheorem{Remark}[Thm]{Remark}
\newtheorem{PropDef}[Thm]{Proposition/Definition}

\theoremstyle{remark}

\newtheorem*{ack}{Acknowledgements}
\newtheorem*{convention}{Conventions}
\newtheorem*{remark}{Remark}

\newtheorem{Ex}[Thm]{Example}

\newcommand{\G}{\mathrm{G}}

\newcommand{\V}{\mathrm{V}}
\newcommand{\W}{\mathrm{W}}

\newcommand{\an}{\mathrm{an}}
\newcommand{\rP}{\mathrm{P}}

\newcommand{\Sr}{\mathrm{S}}

\newcommand{\T}{\mathrm{T}}
\newcommand{\N}{\mathrm{N}}
\newcommand{\I}{\mathrm{I}}

\newcommand{\Hr}{\mathrm{H}}

\DeclareMathOperator{\Stab}{Stab}

\newcommand{\be}{{\bf e}}
\newcommand{\ap}{\mathbb{A}_{\be}}
\newcommand{\compap}{\overline{\mathbb{A}}_\be}
\newcommand{\Bor}{{\rm Bor}}

\author{Bertrand R\'emy}
\address{Universit\'e de Lyon 1\\
Institut Camille Jordan \\
43 boulevard du 11 novembre 1918 \\
F-69622 Villeurbanne cedex}

\email{remy@math.univ-lyon1.fr}

\author{Amaury Thuillier}

\address{Universit\'e de Lyon 1\\
Institut Camille Jordan \\
43 boulevard du 11 novembre 1918 \\
F-69622 Villeurbanne cedex}

\email{thuillier@math.univ-lyon1.fr}

\author{Annette Werner}

\address{Institut f\"ur Mathematik \\
Goethe-Universit\"at Frankfurt \\
Robert-Mayer-Str. 6-8 \\
D-60325 Frankfurt a. M.}

\email{werner@math.uni-frankfurt.de}

 \title{Bruhat-Tits buildings and analytic geometry}

 \alttitle{Immeubles de Bruhat-Tits et g\'eom\'etrie analytique}

\begin{document}

\frontmatter

\begin{abstract}
This paper provides an overview of the theory of Bruhat-Tits buildings. Besides, we explain how Bruhat-Tits buildings can be realized inside Berkovich spaces. In this way, Berkovich analytic geometry can be used to compactify buildings. We discuss in detail the example of the special linear group. 
\end{abstract}
 
\begin{altabstract}
Ce texte introduit les immeubles de
Bruhat-Tits associ\'es aux groupes r\'eductifs sur les corps valu\'es et
explique comment les r\'ealiser et les compactifier au moyen de la
g\'eom\'erie analytique de Berkovich. Il contient une discussion d\'etaill\'ee
du cas du groupe sp\'ecial lin\'eaire. 
\end{altabstract} 

\subjclass{20E42,
51E24,
14L15,
14G22.
}

\keywords{algebraic group, valued field, Berkovich analytic geometry, Bruhat-Tits building, compactification} 
\altkeywords{groupe alg\'ebrique, corps valu\'e, g\'eom\'etrie analytique au sens de Berkovich, immeuble de Bruhat-Tits, compactification}

 \maketitle
\setcounter{tocdepth}{3}
\tableofcontents


\mainmatter

\section*{Introduction}
\label{s - intro}

This paper is mainly meant to be a survey on two papers written by the same authors, namely \cite{RTW1} and \cite{RTW2}.
The general theme is to explain what the theory of analytic spaces in the sense of Berkovich brings to the problem of compactifying Bruhat-Tits buildings.

\vspace{0.1cm}

{\bf 1.}~{\it Bruhat-Tits buildings}.---~
The general notion of a building was introduced by J. Tits in the 60ies \cite{TitsICM}, \cite[Exercises for IV.2]{Lie456}.
These spaces are cell complexes, required to have some nice symmetry properties so that important classes of groups may act on them.
More precisely, it turned out in practice that for various classes of algebraic groups and generalizations, a class of buildings is adapted in the sense that any group from such a class admits a very transitive action on a suitable building.
The algebraic counterpart to the transitivity properties of the action is the possibility to derive some important structure properties for the group.

This approach is particularly fruitful when the class of groups is that of simple Lie groups over non-Archimedean fields, or more generally reductive groups over non-Archimedean valued fields -- see Sect. \ref{s - Bruhat-Tits general}.
In this case the relevant class of buildings is that of Euclidean buildings (\ref{ss - Euclidean buildings}).
{\it This is essentially the only situation in building theory we consider in this paper}.
Its particularly nice features are, among others, the facts that in this case the buildings are (contractible, hence simply connected) gluings of Euclidean tilings and that deep (non-positive curvature) metric  arguments are therefore available; moreover, on the group side, structures are shown to be even richer than expected.
For instance, topologically the action on the buildings enables one to classify and understand maximal compact subgroups (which is useful to representation theory and harmonic analysis) and, algebraically, it enables one to define important integral models for the group (which is again useful to representation theory, and which is also a crucial step towards analytic geometry).

One delicate point in this theory is merely to prove that for a suitable non-Archimedean reductive group, there does exist a nice action on a suitable Euclidean building: this is the main achievement of the work by F.~Bruhat and J.~Tits in the 70ies \cite{BT1a}, \cite{BT1b}.
Eventually, Bruhat-Tits theory suggests to see the Euclidean buildings attached to reductive groups over valued fields (henceforth called {\it Bruhat-Tits buildings}) as non-Archimedean analogues of the symmetric spaces arising from real reductive Lie groups, from many viewpoints at least.

\vspace{0.1cm}

{\bf 2.}~{\it Some compactification procedures}.---~
Compactifications of symmetric spaces were defined and used in the 60ies; they are related to the more difficult problem of compactifying locally symmetric spaces \cite{Satake2}, to probability theory \cite{Furst}, to harmonic analysis...
One group-theoretic outcome is the geometric parametrization of classes of remarkable closed subgroups \cite{Moore}.
For all the above reasons and according to the analogy between Bruhat-Tits buildings and symmetric spaces, it makes therefore sense to try to construct compactifications of Euclidean buildings.

When the building is a tree, its compactification is quite easy to describe \cite{SerreArbres}. In general, and for the kind of compactifications we consider here, the first construction is due to E. Landvogt \cite{La}: he uses there the fact that the construction of the Bruhat-Tits buildings themselves, at least at the beginning of Bruhat-Tits theory for the simplest cases, consists in defining a suitable gluing equivalence relation for infinitely many copies of a well-chosen Euclidean tiling.
In Landvogt's approach, the equivalence relation is extended so that it glues together infinitely many compactified copies of the Euclidean tiling used to construct the building.
Another approach is more group-theoretic  and relies on the analogy with symmetric spaces: since the symmetric space of a simple real Lie group can be seen as the space  of maximal compact subgroups of the group, one can compatify this space by taking its closure in the (compact) Chabauty space of all closed subgroups.
This approach is carried out by Y. Guivarc'h and the first author \cite{GuiRem}; it leads to statements in group  theory which are analogues of \cite{Moore} (e.g., the virtual geometric classification of maximal amenable subgroups) but the method contains an intrinsic limitation due to which one cannot compactify more than the set of vertices of the Bruhat-Tits buildings.

The last author of the present paper also constructed compactifications of Bruhat-Tits buildings, in at least two different ways.
The first way is specific to the case of the general linear group:
going back to Bruhat-Tits' interpretation of Goldman-Iwahori's work
\cite{GoldmanIwahori}, it starts by seeing the Bruhat-Tits building of
${\rm GL}({\rm V})$ -- where ${\rm V}$ is a vector space over a
discretely valued non-Archimedean field -- as the space of (homothety classes of) non-Archimedean norms on ${\rm V}$.
The compactification consists then in adding at infinity the (homothety classes of) non-zero non-Archimedean seminorms on ${\rm V}$.
Note that the symmetric space of ${\rm SL}_n({\bf R})$ is the set of
normalized scalar products on ${\bf R}^n$ and a natural
compactification consists in projectivizing the cone of positive nonzero semidefinite bilinear forms: what is done in \cite{Wer04} is the non-Archimedean analogue of this; it has some connection with Drinfeld spaces and is useful to our subsequent compactification in the vein of Satake's work for symmetric spaces.
The second way is related to representation theory \cite{Wer07}: it provides, for a given group, a finite family of compactifications of the Bruhat-Tits building.
The compactifications, as in E. Landvogt's monograph, are defined by gluing compactified Euclidean tilings but the variety of possibilities comes from exploiting various possibilities of compactifying equivariantly these tilings in connection with highest weight theory.

\vspace{0.1cm}

{\bf 3.}~{\it Use of Berkovich analytic geometry}.---~
The compactifications we would like to introduce here make a crucial use of Berkovich analytic geometry.
There are actually two different ways to use the latter theory for compactifications.

The first way is already investigated by V.~Berkovich himself when the algebraic group under consideration is split \cite[Chap. 5]{Ber1}.
One intermediate step for it consists in defining a map from  the
building to the analytic space attached to the algebraic group: this
map attaches to each point $x$ of the building an affinoid subgroup
${\rm G}_x$, which is characterized by a unique maximal point $\vartheta(x)$  in the ambient analytic space of the group.
The map $\vartheta$ is a closed embedding when the ground field is local; a compactification is obtained when $\vartheta$ is composed with the (analytic map) associated to a fibration from the group to one of its flag varieties.
One obtains in this way the finite family of compactifications described in \cite{Wer07}.
One nice feature is the possibility to obtain easily maps between compactifications of a given group but attached to distinct flag varieties.
This enables one to understand in combinatorial Lie-theoretic terms which boundary components are shrunk when going from a ``big" compactification to a smaller one.

The second way mimics I.~Satake's work in the real case.
More precisely, it uses a highest weight representation of the group in order to obtain a map from the building of the group to the building of the general linear group of the representation space which, as we said before, is nothing else than a space of non-Archimedean norms.
Then it remains to use the seminorm compactification mentioned above by taking the closure of the image of the composed map from the building to the compact space of (homothety classes of) seminorms on the non-Archimedean representation space.

For a given group, these two methods lead to the same family of compactifications, indexed by the conjugacy classes of parabolic subgroups.
One interesting point in these two approaches is the fact that the compactifications are obtained by taking the closure of images of equivariant maps.
The construction of the latter maps is also one of the main difficulties; it is overcome thanks to the fact that Berkovich geometry has a rich formalism which combines techniques from algebraic and analytic geometry (the possibility to use field extensions, or the concept of Shilov boundary, are for instance crucial to define the desired equivariant maps).

\vspace{0.1cm}

{\bf Structure of the paper.}~
In Sect.~1, we define (simplicial and non-simplicial) Euclidean buildings and illustrate the notions in the case of the groups ${\rm SL}_n$; we also show in these cases how the natural group actions on the building encode information on the group structure of rational points.
In Sect.~2, we illustrate general notions thanks to the examples of spaces naturally associated to special linear groups (such as projective spaces); this time the notions are relevant to Berkovich analytic geometry and to Drinfeld upper half-spaces.
We also provide specific examples of compactifications which we generalize later.
In Sect.~3, we sum up quickly what we need from Bruhat-Tits theory, including the existence of integral models for suitable bounded open subgroups; following the classical strategy, we first show how to construct a Euclidean building in the split case by gluing together Euclidean tilings, and then how to rely on Galois descent arguments for non-necessarily split groups.
In Sect.~4, we finally introduce the maps that enable us to obtain compactifications of Bruhat-Tits buildings (these maps from buildings to analytifications of flag varieties have been previously defined by V. Berkovich in the split case); a variant of this embedding approach, close to Satake's ideas using representation theory to compactify symmetric spaces, is also quickly presented.
In the last section, we correct a mistake in the proof of an auxiliary lemma in \cite{RTW1} which requires us to introduce an additional hypothesis for two results of \cite{RTW2}.

\vspace{0.1cm}

\begin{ack}
We warmly thank the organizers of the summer school ``Berkovich spaces" held in Paris in July 2010. We are grateful to the referee for many comments, corrections and some relevant questions, one of which led to Proposition 5.11. Finally, we thank Tobias Schmidt for pointing out that Lemma A.10 of \cite{RTW1} needed to be corrected. 
\end{ack}

\begin{convention}
In this paper, as in \cite{Ber1}, valued fields are assumed to be
non-Archimedean and complete, the valuation ring of such a field $k$
is denoted by $k^\circ$, its maximal ideal is by
$k^{\circ\circ}$ and its residue field by $\widetilde{k} =
k^\circ/k^{\circ \circ}$.
Moreover a \textit{local field} is a non-trivially valued
non-Archimedean field which is locally compact for the topology given
by the valuation (i.e., it is complete, the valuation is discrete and the residue field is finite).
\end{convention}


\section{Buildings and special linear groups}
\label{s - SL(n) Bruhat-Tits}

We first provide a (very quick) general treatment of Euclidean buildings; general references for this notion are \cite{RousseauGrenoble} and \cite{WeissAffine}.
It is important for us to deal with the simplicial as well as the
non-simplicial version of the notion of a Euclidean building because
compactifying Bruhat-Tits buildings via Berkovich techniques uses huge valued fields.
The second part illustrates these definitions for special linear groups; in particular, we show how to interpret suitable spaces of norms to obtain concrete examples of buildings in the case when the algebraic group under consideration is the special linear group of a vector space.
These spaces of norms will naturally be extended to spaces of (homothety classes of) seminorms when buildings are considered in the context of analytic projective spaces.

\subsection{Euclidean buildings}
\label{ss - Euclidean buildings}

Euclidean buildings are non-Archimedean analogues of Riemannian symmetric spaces of the non-compact type, at least in the following sense: if ${\rm G}$ is a simple algebraic group over a valued field $k$, Bruhat-Tits theory (often) associates to ${\rm G}$ and $k$ a metric space, called a Euclidean building, on which ${\rm G}(k)$ acts by isometries in a ``very transitive" way.
This is a situation which is very close to the one where a (non-compact) simple real Lie group acts on its associated (non-positively curved) Riemannian symmetric space.
In this more classical case, the transitivity of the action, the explicit description of fundamental domains for specific (e.g., maximal compact) subgroups and some non-positive curvature arguments lead to deep conjugacy and structure results -- see \cite{MaubonGrenoble} and \cite{ParadanGrenoble} for a modern account.
Euclidean buildings are singular spaces but, by and large, play a similar role for non-Archimedean Lie groups ${\rm G}(k)$ as above.

\subsubsection{Simplicial definition}
\label{sss - simplicial}
The general reference for building theory from the various ``discrete" viewpoints is \cite{AbramenkoBrown}.
Let us start with an affine reflection group, more precisely a {\it Coxeter group of affine type} \cite{Lie456}.
The starting point to introduce this notion is a locally finite family of hyperplanes -- called {\it walls} -- in a Euclidean space [{\bf loc. cit.}, V \S 1 introduction].
An affine Coxeter group can be seen as a group generated by the reflections in the walls, acting properly on the space and stabilizing the collection of walls [{\bf loc. cit.}, V \S 3 introduction]; it is further required that the action on each irreducible factor of the ambient space be via an infinite {\it essential}~group (no non-zero vector is fixed by the group).

\begin{Example}
\label{ex - simplicial apartments}
\begin{itemize}
\item[1.] The simplest (one-dimensional) example of a Euclidean tiling is provided by the real line tesselated by the integers.
The corresponding affine Coxeter group, generated by the reflections in two consecutive vertices (i.e., integers), is the infinite dihedral group ${\rm D}_\infty$.
\item[2.] The next simplest (irreducible) example is provided by the tesselation of the Euclidean plane by regular triangles.
The corresponding tiling group is the Coxeter group of affine type $\widetilde{{\rm A}_2}$; it is generated by the reflections in the three lines supporting the edges of any fundamental triangle.
\end{itemize}
\end{Example}

Note that Poincar\'e's theorem is a concrete source of Euclidean tilings: start with a Euclidean polyhedron in which each dihedral angle between codimension 1 faces is of the form $\frac{\pi}{ m}$ for  some integer $m \geqslant 1$ (depending on the pair of faces), then the group generated by the reflections in these faces is an affine Coxeter group \cite[IV.H.11]{Maskit}.

In what follows, $\Sigma$ is a Euclidean tiling giving rise to a Euclidean reflection group by Poincar\'e's theorem (in Bourbaki's terminology, it can also be seen as the natural geometric realization of the Coxeter complex of an affine Coxeter group, that is the affinization of the Tits' cone of the latter group \cite{Lie456}).

\begin{Def}
\label{defi - simplicial building}
Let $(\Sigma, W)$ be a Euclidean tiling and its associated Euclidean reflection group.
A {\rm (discrete) Euclidean builiding} of type $(\Sigma, W)$ is a polysimplicial complex, say $\mathcal{B}$, which is covered by subcomplexes all isomorphic to $\Sigma$ -- called the {\rm apartments} -- such that the following incidence properties hold.
\begin{enumerate}
\item[SEB 1] Any two cells of $\mathcal{B}$ lie in some apartment.
\item[SEB 2] Given any two apartments, there is an isomorphism between them fixing their intersection in $\mathcal{B}$.
\end{enumerate}
The cells in this context are called {\rm facets} and the group $W$ is called the {\rm Weyl group} of the building $\mathcal{B}$.
The facets of maximal dimension are called {\rm alcoves}.
\end{Def}

The axioms of a Euclidean building can be motivated by metric reasons.
Indeed, once the choice of a $W$-invariant Euclidean metric on $\Sigma$ has been made, there is a natural way the define a distance on the whole building: given any two points $x$ and $x'$ in $\mathcal{B}$, by {\rm (SEB 1)} pick an apartment $\mathbb{A}$ containing them and consider the distance between $x$ and $x'$ taken in $\mathbb{A}$; then {\rm (SEB 2)} implies that the so--obtained non-negative number doesn't depend on the choice of $\mathbb{A}$.
It requires further work to check that one defines in this way a distance on the building (i.e., to check that the triangle inequality holds \cite[Prop. II.1.3]{Parreau}).

\begin{Remark}
\label{rk - metric motivation}
The terminology ``polysimplicial" refers to the fact that a building can be a direct product of simplicial complexes rather than merely a simplicial complex; this is why we provisionally used the terminology ``cells" instead of ``polysimplices" to state the axioms (as already mentioned, cells will henceforth be called facets -- alcoves when they are top-dimensional).
\end{Remark}

Let us provide now some examples of discrete buildings corresponding to the already mentioned examples of Euclidean tilings.

\begin{Example}
\label{ex - simplicial buildings}
\begin{itemize}
\item[1.] The class of buildings of type $({\bf R}, {\rm D}_\infty)$ coincides with the class of trees without terminal vertex
(recall that a tree is a $1$-dimensional simplicial complex -- i.e., the geometric realization of a graph -- without non-trivial loop \cite{SerreArbres}).
\item[2.] A $2$-dimensional $\widetilde{{\rm A}_2}$-building is already impossible to draw, but roughly speaking it can be constructed by gluing half-tilings to an initial one along {\rm walls} (i.e., fixed point sets of reflections) and by iterating these gluings infinitely many times provided a prescribed ``shape" of neighborhoods of vertices is respected -- see Example \ref{ex - SL3bis} for further details on the local description of a building in this case.
\end{itemize}
\end{Example}

It is important to note that axiom (ii) does {\it not} require that the isomorphism between apartments extends to a global automorphism of the ambient building.
In fact, it may very well happen that for a given Euclidean building $\mathcal{B}$ we have ${\rm Aut}(\mathcal{B}) = \{1 \}$ (take for example a tree in which any two distinct vertices have distinct valencies).
However, J. Tits' classification of Euclidean buildings \cite{TitsCome} implies that in dimension $\geqslant 3$ any irreducible building comes -- via Bruhat-Tits theory, see next remark --  from a simple algebraic group over a local field, and therefore admits a large automorphism group.
At last, note that there do exist $2$-dimensional exotic Euclidean buildings, with interesting but unexpectedly small automorphism groups \cite{Barre}.

\begin{Remark}
\label{rk - BrT}
In Sect. \ref{s - Bruhat-Tits general}, we will briefly introduce Bruhat-Tits theory.
The main outcome of this important part of algebraic group theory is that, given a semisimple algebraic group $ {\rm G}$ over a local field $k$, there exists a discrete Euclidean building $\mathcal{B} = \mathcal{B}({\rm G},k)$ on which the group of rational points ${\rm G}(k)$ acts by isometries and {\rm strongly transitively} (i.e., transitively on the inclusions of an alcove in an apartment).
\end{Remark}

\begin{Example}
Let ${\rm G}$ as above be the group ${\rm SL}_3$.
Then the Euclidean building associated to  ${\rm SL}_3$ is a Euclidean building in which every apartment is a Coxeter complex of type $\widetilde{{\rm A}_2}$, that is the previously described $2$-dimensional tiling of the Euclidean space ${\bf R}^2$ by regular triangles.
Strong transitivity of the ${\rm SL}_3(k)$-action means here that given any alcoves (triangles) $c, c'$ and any apartments $\mathbb{A}, \mathbb{A}'$ such that $c \subset \mathbb{A}$ and $c' \subset \mathbb{A}'$ there exists $g \in {\rm SL}_3(k)$ such that $c'=g.c$ and $\mathbb{A}'=g.\mathbb{A}$.
\end{Example}

The description of the apartments doesn't depend on the local field $k$ (only on the Dynkin diagram of the semisimple group in general), but the field $k$ plays a role when one describes the combinatorial neighborhoods of facets, or small metric balls around vertices.
Such subsets, which intersect finitely many facets when $k$ is a local field, are known to be realizations of some (spherical) buildings: these buildings are naturally associated to semisimple groups (characterized by some subdiagram of the Dynkin diagram of ${\rm G}$) over the residue field $\widetilde{k}$ of $k$.

\begin{Example}
\label{ex - SL3bis}
For ${\rm G}={\rm SL}_3$ and $k={\bf Q}_p$, each sufficiently small ball around a vertex is the flag complex of a $2$-dimensional vector space over ${\bf Z}/p{\bf Z}$ and any edge in the associated Bruhat-Tits building is contained in the closure of exactly $p+1$ triangles.
A suitably small metric ball around any point in the relative interior of an edge can be seen as a projective line over ${\bf Z}/p{\bf Z}$, that is the flag variety of ${\rm SL}_2$ over ${\bf Z}/p{\bf Z}$.
\end{Example}

\subsubsection{Non-simplicial generalization}
\label{sss - non-simplicial}
We will see, e.g. in \ref{ss - closed embedding}, that it is often necessary to understand and use reductive algebraic groups over valued fields for {\it non-discrete} valuations even if in the initial situation the ground field is discretely valued.
The geometric counterpart to this is the necessary use of non-discrete Euclidean buildings.
The investigation of such a situation is already covered by the fundamental work by F.~Bruhat and J.~Tits as written in \cite{BT1a} and \cite{BT1b}, but the intrinsic definition of a non-discrete Euclidean building is not given there -- see \cite{TitsCome} though, for a reference roughly appearing at the same time as Bruhat-Tits' latest papers.

The definition of a building in this generalized context is quite similar to the discrete one (\ref{sss - simplicial}) in the sense that it replaces an atlas by a collection of ``slices" which are still called {\it apartments} and turn out to be maximal flat (i.e., Euclidean) subspaces once the building is endowed with a natural distance.
What follows can be found for instance in A. Parreau's thesis \cite{Parreau}.

Let us go back to the initial question.

\begin{Question}
\label{q - non-simplicial}
Which geometry can be associated to a group ${\rm G}(k)$ when ${\rm G}$ is a reductive group over $k$, a (not necessarily discretely) valued field?
\end{Question}

The answer to this question is a long definition to swallow, so we will provide some explanations immediately after stating it.

The starting point is again a $d$-dimensional Euclidean space, say $\Sigma_{\rm vect}$, together with a finite group $\overline W$ in the group of isometries ${\rm Isom}(\Sigma_{\rm vect}) \simeq {\rm O}_d({\bf R})$.
By definition, a {\it vectorial wall} in $\Sigma_{\rm vect}$ is the
fixed-point set in $\Sigma_{\rm vect}$ of a reflection in
$\overline{W}$ and a {\it vectorial Weyl chamber} is a connected component of the complement of the union of the walls in $\Sigma_{\rm vect}$, so that Weyl chambers are simplicial cones.

Now assume that we are given an affine Euclidean space $\Sigma$ with underlying Euclidean vector space $\Sigma_{\rm vect}$.
We have thus ${\rm Isom}(\Sigma) \simeq {\rm Isom}(\Sigma_{\rm vect}) \ltimes \Sigma_{\rm vect}\simeq {\rm O}_d({\bf R}) \ltimes {\bf R}^d$.
We also assume that we are given a group $W$ of (affine) isometries in $\Sigma$ such that the vectorial part of $W$ is $\overline W$ and such that there exists a point $x \in \Sigma$ and a subgroup ${\rm T} \subset {\rm Isom}(\Sigma)$ of translations satisfying $W = W_x \cdot {\rm T}$; we use here the notation $W_x = {\rm Stab}_W(x)$. A point $x$ satisfying this condition is called \emph{special}.

\begin{Def}
\label{defi - non-simplicial building}
Let $\mathcal{B}$ be a set and let $\mathcal{A} = \{f : \Sigma \to \mathcal{B} \}$ be a collection of injective maps, whose images are called {\rm apartments}.
We say that $\mathcal{B}$ is a {\rm Euclidean building} of type $(\Sigma,W)$ if the apartments satisfy the following axioms.
\begin{itemize}
\item[EB 1] The family $\mathcal{A}$ is stable by precomposition with any element of $W$ (i.e., for any $f \in \mathcal{A}$ and any $w \in W$, we have $f \circ w \in \mathcal{A}$).
\item[EB 2] For any $f,f' \in \mathcal{A}$ the subset $\mathcal{C}_{f,f'}= f'^{-1} \bigl( f(\Sigma) \bigr)$ is convex in $\Sigma$ and there exists $w \in W$ such that we have the equality of restrictions $(f^{-1} \circ f') \mid_{\mathcal{C}_{f,f'}} = w \mid_{\mathcal{C}_{f,f'}}$.
\item[EB 3] Any two points of $\mathcal{B}$ are contained in a suitable apartment.
\end{itemize}
At this stage, there is a well-defined map $d : \mathcal{B}\times \mathcal{B} \to {\bf R}_{\geqslant 0}$ and we further require:
\begin{itemize}
\item[EB 4] Given any (images of) Weyl chambers, there is an apartment of $X$ containing sub-Weyl chambers of each.
\item[EB 5] Given any apartment $\mathbb{A}$ and any point $x \in \mathbb{A}$, there is a $1$-lipschitz retraction map $r = r_{x,\mathbb{A}} : \mathcal{B} \to \mathbb{A}$ such that
$r \mid_\mathbb{A} = {\rm id}_\mathbb{A}$ and $r^{-1}(x) = \{x \}$.
\end{itemize}
\end{Def}

The above definition is taken from \cite[II.1.2]{Parreau}; in these axioms a {\it Weyl chamber} is the affine counterpart to the previously defined notion of a {\it Weyl chamber} and a ``sub-Weyl chamber" is a translate of the initial Weyl chamber which is completely contained in the latter.

\begin{Remark}
A different set of axioms is given in G. Rousseau's paper \cite[\S 6]{RousseauGrenoble}.
It is interesting because it provides a unified approach to simplicial and non-simplicial buildings via incidence requirements on apartments.
The possibility to obtain a non-discrete building with Rousseau's axioms is contained in the model for an apartment and the definition of a facet as a filter.
The latter axioms are adapted to some algebraic situations which cover the case of Bruhat-Tits theory over non-complete valued fields -- see \cite[Remark 9.4]{RousseauGrenoble} for more details and comparisons.
\end{Remark}

\begin{Remark}
In this paper we do not use the plain word ``chamber" though it is standard terminology in abstract building theory.
This choice is made to avoid confusion: alcoves here are chambers (in the abstract sense) in Euclidean buildings and parallelism classes of Weyl chambers here are chambers (in the abstract sense) in spherical buildings at infinity of Euclidean buildings \cite[Chap.~8]{WeissAffine}, \cite[11.8]{AbramenkoBrown}.
\end{Remark}

It is easy to see that, in order to prove that the map $d$ defined thanks to axioms {\rm (EB 1)-(EB 3)} is a distance, it remains to check that the triangle inequality holds; this is mainly done by using the retraction given by axiom {\rm (EB 5)}.
The previously quoted metric motivation (Remark \ref{rk - metric motivation}) so to speak became a definition.
Note that the existence of suitable retractions is useful to other purposes.

The following examples of possibly non-simplicial Euclidean buildings correspond to the examples of simplicial ones given in Example \ref{ex - simplicial buildings}.

\begin{Example}
\label{ex - non-simplicial buildings}
\begin{itemize}
\item[1.] Consider the real line $\Sigma = {\bf R}$ and its isometry group ${\bf Z}/2{\bf Z} \ltimes {\bf R}$.
Then a Euclidean building of type $({\bf R}, {\bf Z}/2{\bf Z} \ltimes {\bf R})$ is a real tree -- see below.
\item[2.] For a $2$-dimensional case extending simplicial $\widetilde{{\rm A}_2}$-buildings, a model for an apartment can be taken to be a maximal flat in the symmetric space of ${\rm SL}_3({\bf R})/{\rm SO}(3)$ acted upon by its stabilizer in ${\rm SL}_3({\bf R})$ (using the notion of singular geodesics to distinguish the walls).
There is a geometric way to define the Weyl group and Weyl chambers (six directions of simplicial cones) in this differential geometric context -- see \cite{MaubonGrenoble} for the general case of arbitrary symmetric spaces.
\end{itemize}
\end{Example}

Here is a (purely metric) definition of real trees.
It is a metric space $({\rm X},d)$ with the following two properties:
\begin{itemize}
\item[(i)] it is {\it geodesic}: given any two points $x, x' \in {\rm X}$ there is a (continuous) map $\gamma : [0;d] \to {\rm X}$, where $d = d(x,x')$, such that $\gamma(0)=x$, $\gamma(d)=x'$ and
$d\bigl(\gamma(s),\gamma(t) \bigr) = \, \mid\! s-t \!\mid$ for any $s,t \in [0;d]$;
\item[(ii)] any geodesic triangle is a tripod (i.e., the union of three geodesic segments with a common end-point).
\end{itemize}

\begin{Remark}
\label{rk - asymptotic cones}
Non-simplicial Euclidean buildings became more popular since recent work of geometric (rather than algebraic) nature, where non-discrete buildings appear as asymptotic cones of symmetric spaces and Bruhat-Tits buildings \cite{KleLee}.
\end{Remark}

The remark implies in particular that there exist non-discrete Euclidean buildings in any dimension, which will also be seen more concretely by studying spaces of non-Archimedean norms on a given vector space -- see \ref{ss - SL(n) Bruhat-Tits}.

\begin{Remark}
\label{rk - BrT bis}
Note that given a reductive group ${\rm G}$ over a valued field $k$, Bruhat-Tits theory ``often" provides a Euclidean building on which the group ${\rm G}(k)$ acts strongly transitively in a suitable sense (see Sect. \ref{s - Bruhat-Tits general} for an introduction to this subject).
\end{Remark}

\subsubsection{More geometric properties}
\label{sss - geometry of buildings}
We motivated the definitions of buildings by metric considerations, therefore we must mention the metric features of Euclidean buildings once these spaces have been defined.
First, a Euclidean building always admits a metric whose restriction to any apartment is a (suitably normalized) Euclidean distance \cite[Prop. 6.2]{RousseauGrenoble}.
Endowed with such a distance, a Euclidean building is always a geodesic metric space as introduced in the above metric definition of real trees (\ref{sss - non-simplicial}).

\vspace{0,1cm}

{\it Recall that we use the axioms {\rm (EB)} from Definition~\ref{defi - non-simplicial building} to define a building; moreover we assume that the above metric is complete.}
This is sufficient for our purposes since we will eventually deal with Bruhat-Tits buildings associated to algebraic groups over complete non-Archimedean fields.

\vspace{0,1cm}

Let $(\mathcal{B}, d)$ be a Euclidean building endowed with such a metric.
Then $(\mathcal{B}, d)$ satisfies moreover a remarkable non-positive curvature property, called the {\it ${\rm CAT}(0)$-property} (where ``CAT" seems to stand for Cartan-Alexandrov-Toponogov).
Roughly speaking, this property says that geodesic triangles are at least as thin as in Euclidean planes.
More precisely, the point is to compare a geodesic triangle drawn in $\mathcal{B}$ with ``the" Euclidean triangle having the same edge lengths.
A geodesic space is said to have the {\it ${\rm CAT}(0)$-property}, or to {\it be} CAT(0), if a median segment in each geodesic triangle is at most as long as the corresponding median segment in the comparison triangle drawn in the Euclidean plane ${\bf R}^2$ (this inequality has to be satisfied for all geodesic triangles).
Though this property is stated in elementary terms, it has very deep consequences \cite[\S 7]{RousseauGrenoble}.

\vspace{0,1cm}

One first consequence is the uniqueness of a geodesic segment between any two points \cite[Chap. II.1, Prop. 1.4]{BriHae}.

\vspace{0,1cm}

The main consequence is a famous and very useful fixed-point property.
The latter statement is itself the consequence of a purely geometric one: any bounded subset in a complete, CAT(0)-space has a unique, metrically characterized, circumcenter \cite[11.3]{AbramenkoBrown}.
This implies that if a group acting by isometries on such a space (e.g., a Euclidean building) has a bounded orbit, then it has a fixed point.
This is the {\it Bruhat-Tits fixed point lemma}; it applies for instance to any compact group of isometries.

\vspace{0,1cm}

Let us simply mention two very important applications of the Bruhat-Tits fixed point lemma (for simplicity, we assume that the building under consideration is discrete and locally finite -- which covers the case of Bruhat-Tits buildings for reductive groups over local fields).

\begin{itemize}
\item[1.] The Bruhat-Tits fixed point lemma is used to classify maximal bounded subgroups in the isometry group of a building.
Indeed, it follows from the definition of the compact open topology on the isometry group ${\rm Aut}(\mathcal{B})$ of a building $\mathcal{B}$, that a facet stabilizer is a compact subgroup in ${\rm Aut}(\mathcal{B})$.
Conversely, a compact subgroup has to fix a point and this point can be sent to a point in a given fundamental domain for the action of ${\rm Aut}(\mathcal{B})$ on $\mathcal{B}$ (the isometry used for this conjugates the initial compact subgroup into the stabilizer of a point in the fundamental domain).
\item[2.] Another consequence is that any Galois action on a Bruhat-Tits building has ``sufficiently many" fixed points, since a Galois group is profinite hence compact.
These Galois actions are of fundamental use in Bruhat-Tits theory, following the general idea -- widely used in algebraic group theory --
that an algebraic group ${\rm G}$ over $k$ is nothing else than a split algebraic group over the separable closure $k^s$, namely ${\rm G} \otimes_k k^s$, together with a semilinear action of ${\rm Gal}(k^s/k)$ on ${\rm G} \otimes_k k^s$ \cite[AG \S~\S11-14]{Borel}.
\end{itemize}

Arguments similar to the ones mentioned in 1. imply that, when $k$ is a local field, there are exactly $d+1$ conjugacy classes of maximal compact subgroups in ${\rm SL}_{d+1}(k)$.
They are parametrized by the vertices contained in the closure of a given alcove (in fact, they are all isomorphic to ${\rm SL}_{d+1}(k^\circ)$ and are all conjugate under the action of ${\rm GL}_{d+1}(k)$ by conjugation).

\begin{Remark} One can make 2. a bit more precise.
The starting point of Bruhat-Tits theory is indeed that a reductive group ${\rm G}$ over any field, say $k$, splits -- hence in particular is very well understood -- after extension to the separable closure $k^s$ of the ground field.
Then, in principle, one can go down to the group ${\rm G}$ over $k$ by means of suitable Galois action -- this is one leitmotiv in  \cite{BoTi}.
In particular, Borel-Tits theory provides a lot of information about the group ${\rm G}(k)$ by seeing it as the fixed-point set ${\rm G}(k^s)^{{\rm Gal}(k^s/k)}$.
When the ground field $k$ is a valued field, then one can associate a
Bruhat-Tits building $\mathcal{B} = \mathcal{B}({\rm G},k^s)$ to ${\rm
  G} \otimes_k k^s$ together with an action by isometries of ${\rm
  Gal}(k^s/k)$. The Bruhat-Tits building of ${\rm G}$ over $k$ is
contained in the Galois fixed-point set $\mathcal{B}^{{\rm
    Gal}(k^s/k)}$, but this is inclusion is strict in general: the
Galois fixed-point set is bigger than the desired building
\cite[III]{RousseauOrsay}. Still, this may be a good first approximation of Bruhat-Tits theory to have in mind.
We refer to \ref{sss - descent and functoriality} for further details.
\end{Remark}

\subsection{The ${\rm SL}_n$ case}
\label{ss - SL(n) Bruhat-Tits}

We now illustrate many of the previous notions in a very explicit situation, of arbitrary dimension.
Our examples are spaces of norms on a non-Archimedean vector space.
They provide the easiest examples of Bruhat-Tits buildings, and are also very close to spaces occurring in Berkovich analytic geometry.
In this section, we denote by ${\rm V}$ a $k$-vector space and by $d+1$ its (finite) dimension over $k$.

\emph{Note that until Remark \ref{rk - extended building} we assume that $k$ is a local field.}

\subsubsection{Goldman-Iwahori spaces}
\label{sss - GI spaces}
The materiel of this subsection is classical and could be find, for example, in \cite{Weil2}.

We are interested in the following space.

\begin{Def}
\label{defi - GI}
The {\rm Goldman-Iwahori} space of the $k$-vector space ${\rm V}$ is the space of non-Archimedean norms on ${\rm V}$; we denote it by $\mathcal{N}({\rm V},k)$.
We denote by $\mathcal{X}({\rm V},k)$ the quotient space $\displaystyle {\mathcal{N}({\rm V},k) \big/ \sim}$, where $\sim$ is the equivalence relation which identifies two homothetic norms.
\end{Def}

To be more precise, let $\parallel \cdot \parallel$ and $\parallel \cdot \parallel'$ be norms in  $\mathcal{N}({\rm V},k)$.
We have $\parallel \cdot \parallel \sim \parallel \cdot \parallel'$ if and only if there exists $c > 0$ such that $\parallel\! x \!\parallel \, = \, c \parallel\! x \!\parallel'$ for all $x \in {\rm V}$.
In the sequel, we use the notation $[\cdot]_\sim$ to denote the class with respect to the homothety equivalence relation.

\begin{Example}
\label{ex - norm}
Here is a simple way to construct non-Archimedean norms on ${\rm V}$.
Pick a basis $\mathbf{e} = (e_0, e_1, \dots, e_d)$ in ${\rm V}$.
Then for each choice of parameters $\underline c = (c_0, c_1, \dots, c_d) \in {\bf R}^{d+1}$, we can define the non-Archimedean norm which sends each vector $x = \sum_i \lambda_i e_i$ to $\max_i \{\exp(c_i) \mid\! \lambda_i \!\mid \}$, where $\mid \cdot \mid$ denotes the absolute value of $k$.
We denote this norm by $\parallel \cdot \parallel_{\mathbf{e},{\underline c}}$.

\end{Example}

We also introduce the following notation and terminology.

\begin{Def}
\label{defi - adapted}
\begin{itemize}
\item[(i)] Let $\parallel \cdot \parallel$ be a norm and let $\mathbf{e}$ be a basis in ${\rm V}$.
We say that $\parallel \cdot \parallel$ is {\rm diagonalized} by $\mathbf{e}$ if there exists $\underline c \in {\bf R}^{d+1}$ such that $\parallel \cdot \parallel = \parallel \cdot \parallel_{\mathbf{e},{\underline c}}$; in this case, we also say that the basis $\mathbf{e}$ is {\rm adapted} to the norm $\parallel \cdot \parallel$.
\item[(ii)] Given a basis $\mathbf{e}$, we denote by $\widetilde{\mathbb{A}_{\mathbf{e}}}$ the set of norms diagonalized by $\mathbf{e}$:
$$\widetilde{\mathbb{A}_{\mathbf{e}}} = \{ \parallel \cdot \parallel_{\mathbf{e},{\underline c}} \, : \, \underline c \in {\bf R}^{d+1} \}.$$
\item[(iii)] We denote by ${\mathbb{A}_{\mathbf{e}}}$ the quotient of $\widetilde{\mathbb{A}_{\mathbf{e}}}$ by the homothety equivalence relation:
$\displaystyle {\mathbb{A}_{\mathbf{e}}} = {\widetilde{\mathbb{A}_{\mathbf{e}}} / \sim}$.
\end{itemize}
\end{Def}

Note that the space $\widetilde{\mathbb{A}_{\mathbf{e}}}$ is naturally an affine space with underlying vector space ${\bf R}^{d+1}$: the free transitive ${\bf R}^{d+1}$-action is by shifting the coefficients $c_i$ which are the logarithms of the ``weights"  $\exp(c_i)$ for the norms $\parallel \cdot \parallel_{\mathbf{e},{\underline c}} : \sum_i \lambda_i e_i \mapsto \max_{0 \leqslant i \leqslant d}\{\exp(c_i) \mid\! \lambda_i \!\mid \}$.
Under this identification of affine spaces, we have: $\displaystyle {\mathbb{A}_{\mathbf{e}}} \simeq {{\bf R}^{d+1} / {\bf R}(1,1,\dots, 1)} \simeq {\bf R}^d$.

\begin{Remark}
\label{rk - GI apartment}
The space $\mathcal{X}({\rm V},k)$ will be endowed with a Euclidean building structure (Th.\ref{th - GI building general}) in which the spaces ${\mathbb{A}_{\mathbf{e}}}$ -- with $\mathbf{e}$ varying over the bases of ${\rm V}$ -- will be the apartments.
\end{Remark}

The following fact can be generalized to more general valued fields than local fields but is {\it not}~true in general (Remark \ref{rk - Berko and split norms}).

\begin{Prop}
\label{prop - adapted basis} Every norm of $\mathcal{N}({\rm V},k)$ admits an adapted basis in ${\rm V}$.
\end{Prop}

\begin{proof}
Let $\parallel \cdot \parallel$ be a norm of $\mathcal{N}({\rm V},k)$.
We prove the result by induction on the dimension of the ambient $k$-vector space.
Let $\mu$ be any non-zero linear form on ${\rm V}$.
The map ${\rm V}\setminus\{0 \} \to {\bf R}_+$ sending $y$ to $\frac{\mid\! \mu(y) \!\mid}{\parallel\! y \!\parallel}$ naturally provides, by homogeneity, a continuous map $\phi : \mathbf{P}({\rm V})(k) \to {\bf R}_+$.
Since $k$ is locally compact, the projective space $\mathbf{P}({\rm V})(k)$ is compact, therefore there exists an element $x \in {\rm V}\setminus\{ 0 \}$ at which $\phi$ achieves its supremum, so that
\begin{equation*}\label{eq*}
\frac{\mid\! \mu(z) \!\mid}{ \mid\! \mu(x) \!\mid}
\parallel\! x \!\parallel \leqslant \parallel\! z \!\parallel
\tag{*}
\end{equation*}
 for any $z \in {\rm V}$.

Let $z$ be an arbitrary vector of ${\rm V}$.
We write $\displaystyle z = y + \frac{\mu(z)}{\mu(x)} x$ according to the direct sum decomposition ${\rm V} = {\rm Ker}(\mu) \oplus k x$.
By the ultrametric inequality satisfied by $\parallel \cdot \parallel$, we have

\begin{equation*}\label{eq**}
\parallel\! z \!\parallel \leqslant \max \{\parallel\! y \!\parallel; \frac{\mid\! \mu(z) \!\mid}{\mid\! \mu(x) \!\mid}\parallel\! x \!\parallel \}
\tag{**}
\end{equation*}
and
\begin{equation*}\label{eq***}
\parallel\! y \!\parallel \leqslant \max \{\parallel\! z \!\parallel; \frac{\mid\! \mu(z) \!\mid}{\mid\! \mu(x) \!\mid}\parallel\! x \!\parallel \}~.
\tag{***}
\end{equation*}

Inequality~\eqref{eq*} says that $\displaystyle \max \{\parallel\! z \!\parallel; \frac{\mid\! \mu(z) \!\mid}{\mid\! \mu(x) \!\mid}\parallel\! x \!\parallel \} = \parallel\! z \!\parallel$, so~\eqref{eq***} implies $\parallel\! z \!\parallel \geqslant \parallel\! y \!\parallel$.
The latter inequality together with~\eqref{eq*} implies that $\displaystyle \parallel\! z \!\parallel \geqslant \max \{\parallel\! y \!\parallel; \frac{\mid\! \mu(z) \!\mid}{\mid\! \mu(x) \!\mid}\parallel\! x \!\parallel \}$.
Combining this with~\eqref{eq**} we obtain the equality $\displaystyle \parallel\! z \!\parallel = \max \{\parallel\! y \!\parallel; \frac{\mid\! \mu(z) \!\mid}{\mid\! \mu(x) \!\mid}\parallel\! x \!\parallel \}$.
Applying the induction hypothesis to ${\rm Ker}(\mu)$, we obtain a basis adapted to the restriction of $\parallel \cdot \parallel$ to ${\rm Ker}(\mu)$.
Adding $x$ we obtain a basis adapted to $\parallel \cdot \parallel$, as required (note that $ \frac{\mu(z)}{\mu(x)}$ is the coordinate corresponding to the vector $x$ in any such basis).
\end{proof}

\smallskip

Actually, we can push a bit further this existence result about adapted norms.

\begin{Prop}
\label{prop - adapted basis - 2}
For any two norms of $\mathcal{N}({\rm V},k)$ there is a basis of ${\rm V}$ simultaneously adapted to them.
\end{Prop}

\begin{proof}
We are now given two norms, say $\parallel \cdot \parallel$ and $\parallel \cdot \parallel'$, in $\mathcal{N}({\rm V},k)$.
In the proof of Proposition~\ref{prop - adapted basis}, the choice of a non-zero linear form $\mu$ had no importance.
In the present situation, we will take advantage of this freedom of choice.
We again argue by induction on the dimension of the ambient $k$-vector space.

By homogeneity, the map ${\rm V}\setminus \{0 \} \to {\bf R}_+$ sending $y$ to $\displaystyle \frac{\parallel\! y \!\parallel}{\parallel\! y \!\parallel'}$ naturally provides a continuous map $\psi : \mathbf{P}({\rm V})(k) \to {\bf R}_+$.
Again because the projective space $\mathbf{P}({\rm V})(k)$ is compact, there exists $x \in {\rm V}\setminus \{ 0 \}$ at which $\psi$ achieves its supremum, so that

\smallskip

\hskip 6mm $\displaystyle \frac{\parallel\! y \!\parallel}{\parallel\! x \!\parallel} \leqslant  \frac{\parallel\! y \!\parallel'}{\parallel\! x \!\parallel'}$ for any $y \in {\rm V}$.

\smallskip

\noindent Now we endow the dual space ${\rm V}^*$ with the operator norm $\parallel \cdot \parallel^*$ associated to $\parallel \cdot \parallel$ on ${\rm V}$.
Since ${\rm V}$ is finite-dimensional, by biduality (i.e. the normed vector space version of ${\rm V}^{**} \simeq {\rm V}$), we have the equality
$\displaystyle \parallel\! x \!\parallel = \sup_{\mu \in {\rm V}^*\setminus \{ 0\}} \frac{\mid\! \mu(x) \!\mid}{\parallel\! \mu \!\parallel^*}$.
By homogeneity and compactness, there exists $\lambda \in {\rm V}^* \setminus\{0 \}$ such that $\displaystyle \parallel\! x \!\parallel = \frac{\mid\! \lambda(x) \!\mid}{\parallel\! \lambda \!\parallel^*}$.
For arbitrary $y \in {\rm V}$ we have $\mid\! \lambda(y) \!\mid\, \leqslant \, \parallel\! y \!\parallel \cdot \parallel\! \lambda \!\parallel^*$, so the definition of $x$ implies that

\smallskip

\hskip 6mm $\displaystyle \frac{\mid\! \lambda(y) \!\mid}{\mid\! \lambda(x) \!\mid} \leqslant \frac{\parallel\! y \!\parallel}{\parallel\! x \!\parallel}$ for any $y \in {\rm V}$.

\smallskip

\noindent In other words, we have found $x \in {\rm V}$ and $\lambda \in {\rm V}^*$ such that

\smallskip

\hskip 6mm $\displaystyle \frac{\mid\! \lambda(y) \!\mid}{\mid\! \lambda(x) \!\mid} \leqslant \frac{\parallel\! y \!\parallel}{\parallel\! x \!\parallel}\leqslant  \frac{\parallel\! y \!\parallel'}{\parallel\! x \!\parallel'}$ for any $y \in {\rm V}$.

\smallskip

Now we are in position to apply the arguments of the proof of Proposition~\ref{prop - adapted basis} to both $\parallel \cdot \parallel$ and $\parallel \cdot \parallel'$ to obtain that
$\displaystyle \parallel\! z \!\parallel = \max \{\parallel\! y \!\parallel; \frac{\mid\! \lambda(z) \!\mid}{\mid\! \lambda(x) \!\mid}\parallel\! x \!\parallel \}$ and
$\displaystyle \parallel\! z \!\parallel' = \max \{\parallel\! y \!\parallel'; \frac{\mid\! \lambda(z) \!\mid}{\mid\! \lambda(x) \!\mid}\parallel\! x \!\parallel' \}$ for any $z \in {\rm V}$ decomposed as $z = x + y$ with $y \in {\rm Ker}(\lambda)$.
It remains then to apply the induction hypothesis (i.e., that the desired statement holds in the ambient dimension minus 1).
\end{proof}

\subsubsection{Connection with building theory}
\label{sss - GI building}
It is now time to describe the connection between Goldman-Iwahori spaces and Euclidean buildings.
As already mentioned, the subspaces ${\mathbb{A}_{\mathbf{e}}}$ will be the apartments in $\mathcal{X}({\rm V},k)$ (Remark \ref{rk - GI apartment}).

Let us fix a basis $\mathbf{e}$ in ${\rm V}$ and consider first the bigger affine space $\widetilde{\mathbf{A}_{\mathbf{e}}} = \{ \parallel \cdot \parallel_{\mathbf{e},{\underline c}} \, : \, \underline c \in {\bf R}^{d+1} \} \simeq {\bf R}^{d+1}$.
The symmetric group $\mathcal{S}_{d+1}$ acts on this affine space by permuting the coefficients $c_i$.
This is obviously a faithful action and we have another one given by the affine structure.
We obtain in this way an action of the group $\mathcal{S}_{d+1} \ltimes {\bf R}^{d+1}$ on $\widetilde{\mathbf{A}_{\mathbf{e}}}$ and, after passing to the quotient space, we can see
${\mathbb{A}_{\mathbf{e}}}$ as the ambient space of the Euclidean tiling attached to the affine Coxeter group of type $\widetilde{{\rm A}_d}$ (the latter group is isomorphic to $\mathcal{S}_{d+1} \ltimes {\bf Z}^d$).
The following result is due to Bruhat-Tits, elaborating on Goldman-Iwahori's investigation of the space of norms $\mathcal{N}({\rm V},k)$ \cite{GoldmanIwahori}.

\begin{Thm}
\label{th - GI building simplicial}
The space $\displaystyle \mathcal{X}({\rm V},k) = \displaystyle {\mathcal{N}({\rm V},k) / \sim}$ is a simplicial Euclidean building of type $\widetilde{{\rm A}_d}$, where $d+1 = {\rm dim}({\rm V})$; in particular, the apartments are isometric to ${\bf R}^d$ and the Weyl group is isomorphic to $\mathcal{S}_{d+1} \ltimes {\bf Z}^d$.
\end{Thm}

\begin{proof}[Reference for the proof]
In \cite[10.2]{BT1a} this is stated in group-theoretic terms, so one has to combine the quoted statement with [{\bf loc. cit.}, 7.4] in order to obtain the above theorem.
This will be explained in Sect. \ref{s - Bruhat-Tits general}.
\end{proof}

The $0$-skeleton (i.e., the vertices) for the simplicial structure corresponds to the {\it $k^\circ$-lattices}~in the $k$-vector space ${\rm V}$, that is the free $k^\circ$-submodules in ${\rm V}$ of rank $d+1$.
To a lattice $\mathcal{L}$ is attached a norm $\parallel \cdot \parallel_\mathcal{L}$ by setting
$\parallel\! x \!\parallel_\mathcal{L} = \inf \{\mid\! \lambda \!\mid \, :\, \lambda \in k^\times$ and $\lambda^{-1} x \in \mathcal{L} \}$.
One recovers the $k^\circ$-lattice $\mathcal{L}$ as the unit ball of the norm $\parallel \cdot \parallel_\mathcal{L}$.

\begin{Remark}
\label{rk - extended building}
Note that the space $\mathcal{N}({\rm V},k)$ is an extended building in the sense of \cite{TitsCorvallis}; this is, roughly speaking, a building to which is added a Euclidean factor in order to account geometrically for the presence of a center of positive dimension.
\end{Remark}

Instead of trying to prove this result, let us mention that Proposition~\ref{prop - adapted basis - 2} says, in our building-theoretic context, that any two points are contained in an apartment.
In other words, this proposition implies axiom (SEB 1) of Definition~\ref{defi - simplicial building}: it is the non-Archimedean analogue of the fact that any two real scalar products are diagonalized in a suitable common basis (Gram-Schmidt).

\smallskip

Now let us skip the hypothesis that $k$ is a local field.  If $k$ is a not  discretely valued, then it is not true in general that every norm in  $\mathcal{N}({\rm V},k)$ can be diagonalized in some suitable basis.
Therefore we introduce the following subspace:

\smallskip

\centerline{$\mathcal{N}({\rm V},k)^{\rm diag} = \{$norms in $\mathcal{N}({\rm V},k)$ admitting an adapted basis$\}$.}

\begin{Remark}
\label{rk - Berko and split norms}
We will see (Remark \ref{rk-approximation}) that the connection between Berkovich projective spaces and Bruhat-Tits buildings helps to understand why
$\mathcal{N}({\rm V},k) \,\mathbf{-}\, \mathcal{N}({\rm V},k)^{\rm diag} \neq \varnothing$ if and only if the valued field $k$ is {\rm not} maximally complete (one also says spherically complete).
\end{Remark}

\noindent Thanks to the subspace $\mathcal{N}({\rm V},k)^{\rm diag}$, we can state the result in full generality.

\begin{Thm}
\label{th - GI building general}
The space $\displaystyle \mathcal{X}({\rm V},k) = \displaystyle {\mathcal{N}({\rm V},k)^{\rm diag} / \sim}$ is a Euclidean building of type $\widetilde{{\rm A}_d}$ in which the apartments are isometric to ${\bf R}^d$ and the Weyl group is isomorphic to $\mathcal{S}_{d+1} \ltimes \Lambda$ where $\Lambda$ is a translation group, which is discrete if and only if so is the valuation of $k$.
\end{Thm}

\begin{proof}[Reference for the proof]
This is proved for instance in \cite[III.1.2]{Parreau}; see also \cite{BT1b} for a very general treatment.
\end{proof}

\begin{Ex}
\label{ex - GI real tree}
For $d=1$, i.e. when ${\rm V} \simeq k^2$, the Bruhat-Tits building $$\mathcal{X}({\rm V},k) = \displaystyle {\mathcal{N}({\rm V},k)^{\rm diag} / \sim}$$ given by Theorem~\ref{th - GI building general} is a tree, which is a (non-simplicial) real tree whenever $k$ is not discretely valued.
\end{Ex}

\subsubsection{Group actions}
\label{sss - actions on GI}

After illustrating the notion of a building thanks to Goldman-Iwahori spaces, we now describe the natural action of a general linear group over the valued field $k$ on its Bruhat-Tits building.
We said that buildings are usually used to better understand groups which act sufficiently transitively on them.
We therefore have to describe the ${\rm GL}({\rm V},k)$-action on $\mathcal{X}({\rm V},k)$ given by precomposition on norms (that is, $g.\parallel \cdot \parallel \, = \, \parallel \cdot \parallel \circ \, g^{-1}$ for any $g \in {\rm GL}({\rm V},k)$ and any $\parallel \cdot \parallel \in  \mathcal{N}({\rm V},k)$).
Note that we have the formula

\smallskip

\centerline{$g.\parallel \cdot \parallel_{\mathbf{e},{\underline c}} = \parallel \cdot \parallel_{g.\mathbf{e},{\underline c}}$.}

\smallskip

We will also explain how this action can be used to find interesting decompositions of  ${\rm GL}({\rm V},k)$. Note that the ${\rm GL}(\V,k)$-action on $\mathcal{X}(\V,k)$ factors through an action by the group ${\rm PGL}(\V,k)$.

\smallskip

{\it For the sake of simplicity, we assume that $k$ is discretely valued until the rest of this section}.

\smallskip

We describe successively: the action of monomial matrices on the corresponding apartment, stabilizers, fundamental domains and the action of elementary unipotent matrices on the buildings (which can be thought of as ``foldings" of half-apartments fixing complementary apartments).

First, it is very useful to restrict our attention to apartments.
Pick a basis $\mathbf{e}$ of ${\rm V}$ and consider the associated apartment $\mathbb{A}_\mathbf{e}$.
The stabilizer of $\mathbb{A}_\mathbf{e}$ in ${\rm GL}({\rm V},k)$ consists of the subgroup of linear automorphisms $g$ which are {\it monomial}~with respect to $\mathbf{e}$, that is whose matrix expression with respect to $\mathbf{e}$ has only one non-zero entry in each row and in each column; we denote
${\rm N}_\mathbf{e} = {\rm Stab}_{{\rm GL}({\rm V},k)}(\mathbb{A}_\mathbf{e})$.
Any automorphism in ${\rm N}_\mathbf{e}$ lifts a permutation of the indices of the vectors $e_i$ ($0 \leqslant i \leqslant d)$ in $\mathbf{e}$.
This defines a surjective homomorphism ${\rm N}_\mathbf{e}\twoheadrightarrow \mathcal{S}_{d+1}$ whose kernel is the group, say ${\rm D}_\mathbf{e}$, of the linear automorphisms diagonalized by $\mathbf{e}$.
The group ${\rm D}_\mathbf{e} \cap {\rm SL}({\rm V},k)$ lifts the translation subgroup of the (affine) Weyl group $\mathcal{S}_{d+1} \ltimes {\bf Z}^d$ of $\mathcal{X}({\rm V},k)$.
Note that the latter translation group consists of the translations
contained in the group generated by the reflections in the codimension
1 faces of a given alcove, therefore this group is (of finite index
but) smaller than the ``obvious" group given by translations with
integral coefficients with respect to the basis $\mathbf{e}$. For any
$\underline{\lambda} \in (k^{\times})^n$, we have the following ``translation formula":

\smallskip

\centerline{$\underline{\lambda}.\parallel \cdot \parallel_{\mathbf{e},{\underline c}} = \parallel \cdot \parallel_{\mathbf{e},(c_i - \log \mid\! \lambda_i \!\mid)_i}$,}

\smallskip

\begin{Example}
\label{ex - respect type tree}
When $d=1$ and when $k$ is local, the translations of smallest displacement length in the (affine) Weyl group of the corresponding tree are translations whose displacement length along their axis is equal to twice the length of an edge.
\end{Example}

The fact stated in the example corresponds to the general fact that the ${\rm SL}({\rm V},k)$-action on $\mathcal{X}({\rm V},k)$ is {\it type} (or {\it color}){\it -preserving}: choosing $d+1$ colors, one can attach a color to each {\it panel} (= codimension 1 facet) so that each color appears exactly once in the closure of any alcove; a panel of a given color is sent by any element of ${\rm SL}({\rm V},k)$ to a panel of the same color. Note that the action of ${\rm GL}(\V,k)$, hence also of ${\rm PGL}(\V,k)$, on  $\mathcal{X}({\rm V},k)$ is not type-preserving since ${\rm PGL}(\V,k)$ acts transitively on the set of vertices.

\smallskip

It is natural to first describe the isotropy groups for the action we are interested in.

\begin{Prop}
\label{prop - stabilizers}
We have the following description of stabilizers:

\smallskip

\centerline{
${\rm Stab}_{{\rm GL}({\rm V},k)}(\parallel \cdot \parallel_{\mathbf{e},{\underline c}})
= \{g \in {\rm GL}({\rm V},k) : |{\rm det}(g)| = 1$ and $\log(\mid\! g_{ij} \!\mid) \leqslant c_j - c_i \}$,}

\smallskip

\noindent where $[ g_{ij} ]$ is the matrix expression of ${\rm GL}({\rm V},k)$ with respect to the basis $\mathbf{e}$.
\end{Prop}

\begin{proof}[Reference for the proof]
This is for instance \cite[Cor. III.1.4]{Parreau}.
\end{proof}

There is also a description of the stabilizer group in ${\rm SL}(\V,k)$ as the set of matrices stabilizing a point with respect to a tropical matrix operation \cite[Prop. 2.4]{Wer11}.

\smallskip

We now turn our attention to fundamental domains.
Let $x$ be a vertex in $\mathcal{X}({\rm V},k)$.
Fix a basis $\mathbf{e}$ such that $x = [\parallel \cdot \parallel_{\mathbf{e},{\underline 0}}]_\sim$.
Then we have an apartment $\mathbb{A}_\mathbf{e}$ containing $x$ and the inequalities

\smallskip

\centerline{$c_0 \leqslant c_1 \leqslant \dots \leqslant c_d$}

\smallskip

\noindent
define a Weyl chamber with tip $x$ (after passing to the homothety classes).
The other Weyl chambers with tip $x$ contained in $\mathbb{A}_\mathbf{e}$ are obtained by using the action of the spherical Weyl group $\mathcal{S}_{d+1}$, which amounts to permuting the indices of the $c_i$'s (this action is lifted by the action of monomial matrices with coefficients $\pm 1$ and determinant 1).

Accordingly, if we denote by $\varpi$ a uniformizer of $k$, then the inequalities

\smallskip

\centerline{$c_0 \leqslant c_1 \leqslant \dots \leqslant c_d$ \quad and \quad $c_d - c_0 \leqslant -\log \mid\! \varpi \!\mid$}

\smallskip

\noindent
define an alcove (whose boundary contains $x$) and any other alcove in $\mathbb{A}_\mathbf{e}$ is obtained by using the action of the affine Weyl group $\mathcal{S}_{d+1} \ltimes {\bf Z}^d$.

\begin{Prop}
\label{prop - fundamental domains}
Assume $k$ is local.
We have the following description of fundamental domains.
\begin{itemize}
\item[(i)] Given a vertex $x$, any Weyl chamber with tip $x$ is a fundamental domain for the action of the maximal compact subgroup ${\rm Stab}_{{\rm SL}({\rm V},k)}(x)$ on $\mathcal{X}({\rm V},k)$.

\item[(ii)] Any alcove is a fundamental domain for the natural action of ${\rm SL}({\rm V},k)$ on the building $\mathcal{X}({\rm V},k)$.
\end{itemize}
\end{Prop}

If we abandon the hypothesis that $k$ is a local field and assume the absolute value of $k$ is surjective (onto ${\bf R}_{\geqslant 0}$), then the ${\rm SL}({\rm V},k)$-action on $\mathcal{X}({\rm V},k)$ is transitive.

\smallskip

\begin{proof}[Sketch of proof].
Property~(ii) follows from (i) and from the previous description of the action of the monomial matrices of ${\rm N}_\mathbf{e}$ on $\mathbb{A}_\mathbf{e}$ (note that ${\rm SL}({\rm V},k)$ is type-preserving, so a fundamental domain cannot be strictly smaller than an alcove).

(i). A fundamental domain for the action of the symmetric group $\mathcal{S}_{d+1}$ as above on the apartment $\mathbb{A}_\mathbf{e}$ is given by a Weyl chamber with tip $x$, and the latter symmetric group is lifted by elements in ${\rm Stab}_{{\rm SL}({\rm V},k)}(x)$.
Therefore it is enough to show that any point of the building can be mapped into $\mathbb{A}_\mathbf{e}$ by an element of ${\rm Stab}_{{\rm SL}({\rm V},k)}(x)$.
Pick a point $z$ in the building and consider a basis $\mathbf{e'}$ such that  $\mathbb{A}_\mathbf{e'}$ contains both $x$ and $z$ (Proposition~\ref{prop - adapted basis - 2}).
We can write $x = \, \parallel \cdot \parallel_{\mathbf{e},0} \, = \, \parallel \cdot \parallel_{\mathbf{e'},{\underline c}}$, with weights $\underline c$ in $\log \mid\! k^\times \!\mid$ since $x$ is a vertex.
After dilation, if necessary, of each vector of the basis $\mathbf{e'}$, we may -- and shall -- assume that $\underline c = 0$.
Pick $g \in {\rm SL}({\rm V},k)$ such that $g.\mathbf{e}=\mathbf{e'}$.
Since $\mathbf{e}$ and $\mathbf{e'}$ span the same lattice $L$ over
$k^\circ$, which is the unit ball for $x$ (see comment after Th. \ref{th - GI building simplicial}), we have $g.L=L$ and therefore $g$ stabilizes $x$.
We have therefore found $g \in {\rm Stab}_{{\rm SL}({\rm V},k)}(x)$ with $g.\mathbb{A}_\mathbf{e} = \mathbb{A}_\mathbf{e'}$, in particular $g^{-1}.z$ belongs to $\mathbb{A}_\mathbf{e}$.
\end{proof}

\begin{Remark}
\label{rk - geometric Cartan}
Point (i) above is the geometric way to state the so-called Cartan decomposition:
${\rm SL}({\rm V},k) = {\rm Stab}_{{\rm SL}({\rm V},k)}(x) \cdot \overline{{\rm T}^+} \cdot {\rm Stab}_{{\rm SL}({\rm V},k)}(x)$, where $\overline{{\rm T}^+}$ is the semigroup of linear automorphisms $t$ diagonalized by $\mathbf{e}$ and such that $t.x$ belongs to a fixed Weyl chamber in $\mathbb{A}_\mathbf{e}$ with tip $x$.
The Weyl chamber can be chosen so that $\overline{{\rm T}^+}$ consists of the diagonal matrices whose diagonal coefficients are powers of some given uniformizer with the exponents increasing along the diagonal.
Let us recall how to prove this by means of elementary arguments \cite[\S 3.4 p. 152]{PlatonovRapinchuk}.
Let $g \in {\rm SL}({\rm V},k)$; we pick $\lambda \in k^\circ$ so that $\lambda g$ is a matrix of ${\rm GL}({\rm V},k)$ with coefficients in $k^\circ$.
By interpreting left and right multiplication by elementary unipotent matrices as matrix operations on rows and columns, and since $k^\circ$ is a principal ideal domain, we can find $p,p' \in {\rm SL}_{d+1}(k^\circ)$ such that $p^{-1}\lambda g p'^{-1}$ is a diagonal matrix (still with coefficients in $k^\circ$), which we denote by $d$.
Therefore, we can write $g = p \lambda^{-1}d p'$; and since $g$, $p$ and $p'$ have determinant 1, so does $t=\lambda^{-1}d$.
It remains to conjugate $\lambda^{-1}d$ by a suitable monomial matrix with coefficients $\pm 1$ and determinant 1 in order to obtain the desired decomposition.
\end{Remark}

\smallskip

At the beginning of this subsection, we described the action of linear automorphisms on an apartment when the automorphisms are diagonalized by a basis defining the apartment.
One last interesting point is the description of the action of elementary unipotent matrices (for a given basis).
The action looks like a ``folding" in the building, fixing a suitable closed half-apartment.

More precisely, let us introduce the elementary unipotent matrices  $u_{ij}(\nu) = {\rm id}+ \nu {\rm E}_{ij}$ where $\nu \in k$ and ${\rm E}_{ij}$ is the matrix whose only non-zero entry is the $(i,j)$-th one, equal to 1.

\begin{Prop}
\label{prop - folding}
The intersection $\widetilde{\mathbb{A}_{\mathbf{e}}} \cap u_{ij}(\lambda).\widetilde{\mathbb{A}_{\mathbf{e}}}$ is the half-space of $\widetilde{\mathbb{A}_{\mathbf{e}}}$ consisting of the norms $\parallel \cdot \parallel_{\mathbf{e},{\underline c}}$ satisfying $c_j - c_i \geqslant \log \mid\! \lambda \!\mid$.
The isometry given by the matrix $u_{ij}(\lambda)$ fixes pointwise this intersection and
the image of the open half-apartment $\widetilde{\mathbb{A}_{\mathbf{e}}} \mathbf{-} \{\parallel\! \cdot \!\parallel_{\mathbf{e},{\underline c}} : c_j - c_i \geqslant \log \mid\! \lambda \!\mid \}$ is (another half-apartment) disjoint from $\widetilde{\mathbb{A}_{\mathbf{e}}}$.
\end{Prop}

\begin{proof}
In the above notation, we have $u_{ij}(\nu)(\sum_i \lambda_i e_i) = \sum_{k \neq i} \lambda_k e_k + (\lambda_i + \nu\lambda_j) e_i$ for any $\nu \in k$.

\smallskip

First, we assume that we have $u_{ij}(\lambda).\parallel \cdot \parallel_{\mathbf{e},{\underline c}} = \parallel \cdot \parallel_{\mathbf{e},{\underline c}}$.
Then, applying this equality of norms to the vector $e_j$ provides  $e^{c_j} = \max \{ e^{c_j} ; e^{c_i} \mid\! \lambda \!\mid \}$,
hence the inequality $c_j - c_i \geqslant \log \mid\! \lambda \!\mid$.

\smallskip

Conversely, pick a norm $\parallel\! \cdot \!\parallel_{\mathbf{e},{\underline c}}$ such that $c_j - c_i \geqslant \log \mid\! \lambda \!\mid$ and let $x = \sum_i \lambda_i e_i$.
By the ultrametric inequality, we have $e^{c_i}\mid\! \lambda_i-\lambda\lambda_j\!\mid \leqslant \max \{ e^{c_i}\mid\! \lambda_i \!\mid ; e^{c_i} \mid\! \lambda \!\mid \mid\! \lambda_j \!\mid \}$, and the assumption $c_j - c_i \geqslant \log \mid\! \lambda \!\mid$ implies that $e^{c_i}\mid\! \lambda_i-\lambda\lambda_j\!\mid \leqslant \max \{ e^{c_i}\mid\! \lambda_i \!\mid ; e^{c_j}\mid\! \lambda_j \!\mid\}$, so that
$e^{c_i}\mid\! \lambda_i-\lambda\lambda_j\!\mid \leqslant \max_{1 \leqslant \ell \leqslant d} e^{c_\ell} \mid\! \lambda_\ell \!\mid$.
Therefore we obtain that $u_{ij}(\lambda).\parallel x \parallel_{\mathbf{e},{\underline c}} \leqslant \parallel x \parallel_{\mathbf{e},{\underline c}}$ for any vector $x$.
Replacing $\lambda$ by $-\lambda$ and $x$ by $u_{ij}(-\lambda).x$, we finally see that the norms $u_{ij}(\lambda).\parallel \cdot \parallel_{\mathbf{e},{\underline c}}$ and $\parallel \cdot \parallel_{\mathbf{e},{\underline c}}$ are the same when $c_j - c_i \geqslant \log \mid\! \lambda \!\mid$.
We have thus proved that the fixed-point set of $u_{ij}(\lambda)$ in $\widetilde{\mathbb{A}_{\mathbf{e}}}$ is the closed half-space
${\rm D}_\lambda = \{\parallel\! \cdot \!\parallel_{\mathbf{e},{\underline c}} : c_j - c_i \geqslant \log \mid\! \lambda \!\mid \}$.

\smallskip

It follows from this that $\widetilde{\mathbb{A}_{\mathbf{e}}} \cap u_{ij}(\lambda).\widetilde{\mathbb{A}_{\mathbf{e}}}$ contains ${\rm D}_\lambda$.
Assume that $\widetilde{\mathbb{A}_{\mathbf{e}}} \cap u_{ij}(\lambda).\widetilde{\mathbb{A}_{\mathbf{e}}} \supsetneq {\rm D}_\lambda$ in order to obtain a contradiction.
This would provide norms $\parallel \cdot \parallel$ and $\parallel \cdot \parallel'$ in $\widetilde{\mathbb{A}_{\mathbf{e}}} - {\rm D}_\lambda$ with the property that $u_{ij}(\lambda).\parallel \cdot \parallel=\parallel \cdot \parallel'$.
But we note that a norm in $\widetilde{\mathbb{A}_{\mathbf{e}}}-{\rm D}_\lambda$ is characterized by its orthogonal projection onto the boundary hyperplane $\partial {\rm D}_\lambda$ and by its distance to $\partial {\rm D}_\lambda$.
Since $u_{ij}(\lambda)$ is an isometry which fixes ${\rm D}_\lambda$ we conclude that $\parallel \cdot \parallel = \parallel \cdot \parallel'$, which is in contradiction with the fact that the fixed-point set of $u_{ij}(\lambda)$ in $\widetilde{\mathbb{A}_{\mathbf{e}}}$ is exactly  ${\rm D}_\lambda$.
\end{proof}


\section{Special linear groups, Berkovich and Drinfeld spaces}
\label{s - SL(n) Berkovich-Drinfeld}

We ended the previous section by an elementary construction of
the building of special linear groups over discretely valued
non-Archimedean field. The generalization to an arbitrary reductive
group over such a field is significantly harder and requires the full
development of Bruhat-Tits, which will be the topic of Section \ref{s
  - Bruhat-Tits general}. Before diving into the subtelties of
buildings construction, we keep for a moment the particular case of
special linear groups and describe a realization of their buildings
in the framework of Berkovich's analytic geometry, which leads very
naturally to a compactification of those buildings. The general
picture, namely Berkovich realizations and compactifications of
general Bruhat-Tits buildings will be dealt with in Sect. \ref{s -
  general compactifications}).

\smallskip
Roughly speaking understanding the realization (resp. compactification) described below
of the building of a special linear group amounts to understanding
(homothety classes of) norms on a non-Archimedean vector space
(resp. their degenerations), using the viewpoint of multiplicative seminorms on the corresponding symmetric algebra.

A useful reference for Berkovich theory is \cite{Temkin}.
{\it Unless otherwise indicated, we assume in this section that $k$ is a local field}.

\subsection{Drinfeld upper half spaces and Berkovich affine and projective spaces}
\label{ss - Drinfeld}

Let ${\rm V}$ be a finite-dimensional vector space over $k$, and let $\Sr^{\bullet}\V$ be the symmetric algebra of $\V$.
It is a graded $k$-algebra of finite type. Every choice of a basis $v_0, \ldots, v_d$ of ${\rm V}$ induces an isomorphism of $\Sr^{\bullet}\V$ with the polynomial ring over $k$ in $d+1$ indeterminates.
The affine space $\mathbf{A}(\V)$ is defined as the spectrum $\mathrm{Spec}(\Sr^{\bullet}\V)$, and the projective space $\mathbf{P}(\V)$ is defined as the projective spectrum $\mathrm{Proj}(\Sr^{\bullet}\V)$.
These algebraic varieties give rise to analytic spaces in the sense of Berkovich, which we briefly describe below.

\subsubsection{Drinfeld upper half-spaces in analytic projective spaces}
\label{sss - Drinfeld}
As a topological space, the Berkovich affine space $\mathbf{A}(\V)^{\an}$ is the set of all multiplicative seminorms on $\Sr^{\bullet}\V$ extending the absolute value on $k$ together with the topology of pointwise convergence.
The Berkovich projective space $\mathbf{P}(\V)^{\rm an}$ is the quotient of $\mathbf{A}(\V)^{\an} - \{0\}$ modulo the equivalence relation $\sim$ defined as follows: $\alpha \sim \beta$, if and only if there exists a constant  $c>0$ such that for all $f$ in $\Sr^n \V$ we have $\alpha(f) = c^n \beta(f)$.
There is a natural ${\rm PGL}({\rm V})$-action on $\mathbf{P}(\V)^{\rm an}$ given by $g \alpha = \alpha \circ g^{-1}$.
From the viewpoint of Berkovich geometry, Drinfeld upper half-spaces can be introduced as follows \cite{Ber3}.

\begin{Def}
We denote by $\Omega$  the complement of the union of  all $k$-rational hyperplanes in $\mathbf{P}(\V)^{\rm an}$.
The analytic space $\Omega$ is called Drinfeld upper half space.
\end{Def}

Our next goal is now to mention some connections between the above analytic spaces and the Euclidean buildings defined in the previous section.

\subsubsection{Retraction onto the Bruhat-Tits building}
\label{sss - retraction}
Let $\alpha$ be a point in $\mathbf{A}(\V)^{\rm an}$, i.e.  $\alpha$ is a multiplicative seminorm on $\Sr^{\bullet}\V$.
If $\alpha$ is not contained in any $k$-rational hyperplane of $\mathbf{A}(\V)$, then by definition $\alpha$ does not vanish on any element of $\Sr^1 \V = \V$.
Hence the restriction of the seminorm $\alpha$ to the degree one part $\Sr^1 \V = \V$ is a norm.
Recall that the Goldman-Iwahori space $\mathcal{N}({\rm V},k)$ is defined as the set of all non-Archimedean norms on ${\rm V}$, and that $\mathcal{X}(\V,k)$ denotes the quotient space after the homothety relation (\ref{sss - GI spaces}).
Passing to the quotients we see that restriction of seminorms induces a map
\[\tau: \Omega \longrightarrow \mathcal{X}({\rm V},k).\]
If we endow the Goldman-Iwahori space $\mathcal{N}({\rm V},k)$ with the coarsest topology, so that all evaluation maps on a fixed $v \in \V$ are continuous,  and $\mathcal{X}(\V,k)$ with the quotient topology, then $\tau$ is continuous.
Besides, it is equivariant with respect to the action of ${\rm PGL}({\rm V},k)$.
We refer to \cite[\S 3]{RTW2} for further details.

\subsubsection{Embedding of the building (case of the special linear group)}
\label{sss - embedding PGL}
Let now $\gamma$ be a non-trivial norm on ${\rm V}$.
By Proposition \ref{prop - adapted basis}, there exists a basis $e_0, \ldots, e_{d}$ of $\V$ which is adapted to $\gamma$, i.e. we have

\vspace{3mm}
\centerline{$\gamma\bigl(\sum_{i} \lambda_i e_i\bigr) = \max_i \{\exp(c_i) | \lambda_i | \}$}
\vspace{3mm}

\noindent for some real numbers $c_0, \ldots, c_d$.
We can associate to $\gamma$ a multiplicative seminorm $j(\gamma) $ on $\Sr^{\bullet}\V$ by mapping the polynomial $\sum_{I= (i_0, \ldots, i_d)} a_I e_0^{i_0} \ldots e_d^{i_d}$ to $\max_I \{|a_I| \exp(i_0 c_0 + \ldots + i_d c_d)\}$.
Passing to the quotients, we get a continuous map
\[j: \mathcal{X}(\V,k) \longrightarrow \Omega \] satisfying $\tau\bigl(j(\alpha)\bigr) = \alpha$.

Hence $j$ is injective and a homeomorphism onto its image.
Therefore the map $j$ can be used to realize  the Euclidean building  $\mathcal{X}({\rm V},k)$ as a subset of a Berkovich analytic space.
This observation is due to Berkovich, who used it to determine the automorphism group of $\Omega$ \cite{Ber3}.

\begin{Remark}
\label{rk-approximation}
In this remark, we remove the assumption that $k$ is local and we recall that the building $\mathcal{X}({\rm V},k)$ consists of homothety classes of \emph{diagonalizable} norms on ${\rm V}$ (Theorem \ref{th - GI building general}). Assuming ${\rm dim}({\rm V})=2$ for simplicity, we want to rely on analytic geometry to prove the existence of non-diagonalizable norms on ${\rm V}$ for some $k$.

\smallskip 

The map $j : \mathcal{X}({\rm V},k) \rightarrow \mathbf{P}^1({\rm V})^{\rm an}$ can be defined without any assumption on $k$. Given any point $x \in \mathcal{X}({\rm V},k)$, we pick a basis $\mathbf{e} = (e_0,e_1)$ diagonalizing $x$ and define $j(x)$ to be the multiplicative norm on ${\rm S}^\bullet({\rm V})$ mapping an homogeneous polynomial $f = \sum_\nu a_\nu e_0^{\nu_0} e_1^{\nu_1}$ to $\max_{\nu} \{|a_\nu| \cdot |e_0|(x)^{\nu_0} \cdot |e_1|(x)^{\nu_1}\}$. We do not distinguish between $\mathcal{X}({\rm V},k)$ and its image by $j$ in $\mathbf{P}({\rm V})^{\rm an}$, which consists only of points of types 2 and 3 (this follows from \cite[3.2.11]{Temkin}).

\smallskip Let us now consider the subset $\Omega'$ of $\Omega = \mathbf{P}({\rm V})^{\rm an} - \mathbf{P}({\rm V})(k)$ consisting of multiplicative norms on ${\rm S}^\bullet({\rm V})$ whose restriction to ${\rm V}$ is diagonalizable. The map $\tau$ introduced above is well-defined on $\Omega'$ by $\tau(z) = z_{|{\rm V}}$. This gives a continuous retraction of $\Omega'$ onto $\mathcal{X}({\rm V},k)$. The inclusion $\Omega' \subset \Omega$ is strict in general, i.e. if $k$ is not local. For example, assume that $k = \mathbf{C}_p$ is the completion of an algebraic closure of $\mathbf{Q}_p$; this non-Archimedean field is algebraically closed but not spherically complete. In this situation, $\Omega$ contains a point $z$ of type 4 \cite[2.3.13]{Temkin}, which we can approximate by a sequence $(x_n)$ of points in $\mathcal{X}({\rm V},k)$ (this is the translation of the fact that $z$ corresponds to a decreasing sequence of closed balls in $k$ with empty intersection \cite[2.3.11.(iii)]{Temkin}). Now, if $z \in \Omega'$, then $\tau(z)=\tau~(\lim x_n) = \lim~\tau(x_n) = \lim~x_n$ and therefore $z$ belongs to $\mathcal{X}({\rm V},k)$. Since the latter set contains only points of type 2 or 3, this cannot happen and $z \notin \Omega'$; in particular, the restriction of $z$ to ${\rm V}$ produces a norm which is not diagonalizable.
\end{Remark}

\subsection{A first compactification}
\label{ss - seminorm compactification}
Let us now turn to compactification of the building $\mathcal{X}({\rm
  V},k)$. We give an outline of the construction and refer to \cite[\S
  3]{RTW2} for additional details. The generalization to arbitrary
reductive groups is the subject of \ref{ss-maps_to_flags}.
Recall that we assume that $k$ is a local field.
\subsubsection{The space of seminorms}
\label{sss - embedding seminorms}
Let us consider the set  $\mathcal{S}(\V,k)$ of non-Archimedean seminorms on $\V$.
Every non-Archimedean seminorm $\gamma$ on ${\rm V}$ induces a norm on
the quotient space $\V / \mathrm{ker}(\gamma)$. Hence using
Proposition \ref{prop - adapted basis}, we find that there exists a
basis $e_0, \ldots, e_d$ of $\V$ such that $\alpha\bigl(\sum_{i}
\lambda_i e_i\bigr) =  \max_i \{r_i \mid\! \lambda_i \!\mid \}$ for
some non-negative real numbers $r_0, \ldots, r_d$. In this case we say that $\alpha$ is diagonalized by $\be$.
Note that in contrast to Definition \ref{defi - adapted} we do no longer  assume that the $r_i$ are non-zero and hence exponentials.


We can extend $\gamma$ to a seminorm $j(\gamma)$ on the symmetric algebra $\Sr^\bullet \V \simeq k[e_0,\ldots, e_d]$ as follows:

\vspace{3mm}
\centerline{$j(\gamma)\Bigl(\sum_{I= (i_0,\ldots, i_d)}  a_I e_0^{i_0} \ldots e_d^{i_d}\Bigr) = \max \{|a_I| r_0^{i_0} \ldots r_d^{i_d}\}$.}
\vspace{3mm}

We denote by $\overline{\mathcal{X}}(\V,k)$ the quotient of $\mathcal{S}(\V,k) \setminus \{0\}$ after the equivalence relation $\sim$ defined as follows: $\alpha \sim \beta$ if and only if there exists a real constant $c$ with $\alpha = c \beta$. We equip $\mathcal{S}(\V,k)$ with the topology of pointwise convergence and $\overline{\mathcal{X}}(V,k)$ with the quotient topology.
Then the association $\gamma \mapsto j(\gamma)$  induces a continuous and ${\rm PGL}(\V,k)$-equivariant map
\[j: \overline{\mathcal{X}}(\V,k) \rightarrow \mathbf{P}(\V)^{\rm an}\]
which extends the map $j: \mathcal{X}(\V,k) \rightarrow \Omega$ defined in the previous section.

\subsubsection{Extension of the retraction onto the building}
\label{sss - extension of retraction}
Moreover, by restriction to the degree
one part $\Sr^1 \V = \V$, a non-zero multiplicative seminorm on $\Sr^\bullet \V$ yields an element in $\mathcal{S}(\V,k) - \{0 \}$.
Passing to the quotients, this induces a map
\[\tau:  \mathbf{P}(\V)^{\rm an} \longrightarrow \overline{\mathcal{X}}(\V,k)\]
extending the map $\tau: \Omega \rightarrow \mathcal{X}(\V,k) $ defined in section \ref{ss - Drinfeld}.

As in section \ref{ss - Drinfeld}, we see that $\tau \circ j$ is the
identity on $\overline{\mathcal{X}}(\V,k)$, which implies that $j$ is
injective: it is a homeomorphism onto its (closed) image in
$\mathbf{P}(\V)^{\rm an}$. Since $\mathbf{P}(\V)^{\rm an}$ is compact,
we deduce that the image of $j$, and hence $\overline{\mathcal{X}}(\V,k)$, is compact. As $\mathcal{X}(\V,k)$ is an open subset of $\overline{\mathcal{X}}(\V,k)$, the latter space is a compactification of the Euclidean building $\mathcal{X}(\V,k)$; it was studied in \cite{Wer04}.

\subsubsection{The strata of the compactification}
\label{sss - strata}
For every proper subspace ${\rm W}$ of ${\rm V}$ we can extend norms on ${\rm V}/{\rm W}$ to non-trivial seminorms on ${\rm V}$ by composing the norm with the quotient map $\V \rightarrow {\rm V}/{\rm W}$.
This defines a continuous embedding
\[\mathcal{X}(\V/\W,k) \rightarrow \mathcal{\overline{X}}(\V,k).\]
Since every seminorm on ${\rm V}$ is induced in this way from a norm
on the quotient space after its kernel, we find that
$\overline{\mathcal{X}}({\rm V},k)$ is the disjoint union of all Euclidean buildings
$\mathcal{X}({\rm V}/{\rm W},k)$, where ${\rm W}$ runs over all proper subspaces of ${\rm V}$. Hence our compactification of the Euclidean building $\mathcal{X}(\V,k)$ is a union of Euclidean buildings of smaller rank.

\subsection{Topology and group action }

\label{ss - topology}

We will now investigate the convergence of sequences in $\overline{\mathcal{X}}({\rm V},k)$ and deduce that it is compact. We also analyze the action of the group ${\rm SL}({\rm V},k)$ on this space.

\subsubsection{Degeneracy of norms to seminorms and compactness}

Let us first investigate convergence to the boundary of $\mathcal{X}(\V, k)$ in $\overline{\mathcal{X}} (\V,k) = (\mathcal{S}(\V,k) \backslash \{0\})/ \sim$.
We fix a basis $\be= ( e_0, \ldots, e_d)$  of $\V$ and denote by $\mathbb{A}_{\be} $ the corresponding apartment associated to the norms diagonalized by $\be$ as in Definition \ref{defi - adapted}. We denote by
$\compap \subset \overline{\mathcal{X}}(V,k)$ all classes of {\it seminorms} which are diagonalized by $\be$.

We say that a sequence $(z_n)_n$ of points in $\compap$ is distinguished, if there
exists a non-empty subset $\I$ of $\{0,\ldots, d\}$ such that: \begin{itemize}
\item[(a)] For all $i \in I$ and all $n$ we have $z_n(e_i) \neq  0$.
 \item[(b)] for any $i,j \in \I$, the sequence $\left(\frac{z_n (e_j)}{z_n(e_i)}\right)_n$ converges to a positive real number; \item[(c)] for any $i \in \I$ and $j \in \{0,\ldots, d\} - \I$, the sequence $\left(\frac{z_n(e_j)}{z_n(e_i)}\right)_n$ converges to $0$.\end{itemize}
Here we define $\left(\frac{z_n(e_i)}{z_n(e_j)}\right)_n$ as $\left(\frac{x_n(e_i)}{x_n(e_j)}\right)_n$ for an arbitrary representative $x_n \in \mathcal{S}(\V,k)$ of the class $z_n$. Note that this expression does not depend on the choice of the representative $x_n$.

\begin{Lemma} \label{distinguished}
Let $(z_n)_n$ be a distinguished sequence of points in $\compap$. Choose some element $i \in I$. We define a point $z_\infty$ in $\mathcal{S}(\V,k)$ as the homothety class of the seminorm $x_\infty$ defined as follows:
\[x_\infty(e_j) = \left\{ \begin{array}{ll}
 \lim_n \left(\frac{z_n(e_j)}{z_n(e_i)}\right) & \text{ if } j \in \I \\
 0 & \textrm{ if } j \notin \I \end{array} \right. \]
and $x_\infty(\sum_j a_j e_j) = \max |a_j| x_\infty (e_j)$.
Then $z_\infty$ does not depend on the choice of $i$, and the sequence $(z_n)_n$ converges to $z_\infty$ in $\overline{\mathcal{X}}(\V,k)$.
\end{Lemma}

\begin{proof}
Let $x_n$  be a representative of $z_n$ in $\mathcal{S}(\V,k)$.
For $i,j$ and $\ell$ contained in $I$ we have

\[\lim_n \left(\frac{x_n(e_j)}{x_n(e_\ell)}\right) = \lim_n \left(\frac{x_n(e_j)}{x_n(e_i)}\right) \lim_n \left(\frac{x_n(e_i)}{x_n(e_\ell)}\right),\]
which implies that the definition of the seminorm class $z_\infty$ does not depend on the choice of $i \in I$.

The convergence statement is obvious, since the seminorm $x_n$ is equivalent to $ (x_n(e_i))^{-1} x_n$.
\end{proof}

Hence the distinguished sequence of norm classes $(z_n)_n$ considered in the Lemma converges to
a seminorm class whose kernel $W_I$ is spanned by all $e_j$ with $j \notin I$. Therefore the limit point $z_\infty$ lies in the Euclidean building $\mathcal{X}(\V/W_I)$ at the boundary.

\smallskip

Note that the preceding Lemma implies that $\compap$ is the closure of $\ap$ in $\overline{\mathcal{X}}(\V,k)$. Namely, consider $z \in \compap$, i.e. $z$ is the class of a seminorm $x$ on $\V$ which is diagonalizable by $\be$. For every $n$ we define a norm $x_n$ on $\V$ by
\[ x_n(e_i) = \left\{ \begin{array}{ll} x(e_i), & \text{ if } x(e_i) \neq 0\\
                       \frac{1}{n},  & \text{ if } x(e_i) = 0 \end{array} \right.
                     \]
 and
\[x_n(\sum_i a_i e_i) = \max_i |a_i| x_n(e_i).\]

Then the sequence of norm classes $x_n= [z_n]_\sim$ in $\ap$ is distinguished with respect to the set $I = \{i: x(e_i) \neq 0\}$ and it converges towards $z$.

We will now deduce from these convergence results that the space of seminorms is compact. We begin by showing that $\compap$ is compact.
\begin{Prop}
\label{subsequences}Let $(z_n)_n$ be a sequence of points in $\compap$.
Then $(z_n)_n$ has a converging subsequence.
\end{Prop}

\begin{proof}
Let $x_n$ be seminorms representing the points $z_n$.
By the box principle, there exists an index $i \in \{0, \ldots, d\}$ such that
after passing to a subsequence we have
\[ x_n(e_i) \geqslant x_n(e_j) \text{ for all } j = 0, \ldots, d,
n\geqslant 0.\]
In particular we have $x_n(e_i) > 0$. For each $j = 0, \ldots, d$ we look at
the sequence
\[\beta(j)_n = \frac{x_n(e_j)}{x_n(e_i)}\]
which lies between zero and one. In particular, $\beta(i)_n = 1$ is constant.

After passing to a subsequence of $(z_n)_n$ we may -- and shall -- assume that all sequences $\beta(j)_n$ converge to some $\beta(j)$ between zero and one.
Let $I$ be the set of all $j= 0, \ldots, n$ such that $\beta(j) > 0$.
Then a  subsequence of $(z_n)_n$ is distinguished with respect to $I$, hence it converges by Lemma \ref{distinguished}.
\end{proof}

Since $\compap$ is metrizable, the preceding proposition shows that $\compap$ is compact.

We can now describe the ${\rm SL}({\rm V},k)$-action on the seminorm
compactification of the Goldman-Iwahori space of ${\rm V}$. As before,
we fix a basis $\mathbf{e} = (e_0, \ldots, e_n)$.

Let $o$ be the homothety class of the norm on $\V$ defined
by $$\left|\sum_{i=0}^{d} a_i e_i \right|(o) = \max_{0\leqslant i
  \leqslant d} |a_i|$$ and let $${\rm P}_o = \{g \in {\rm SL}(\V,k) \ ; \ g \cdot o \sim o\}$$
be the stabilizer of $o$. It follows from Proposition \ref{prop -
  stabilizers} that ${\rm P}_o= {\rm SL}_{d+1}(k^0)$ with respect to the basis $\be$.

\begin{Lemma}
\label{closedchamber}
The map ${\rm P}_o \times \compap \rightarrow \overline{\mathcal{X}}(\V,k)$ given by the ${\rm SL}(\V,k)$-action is surjective.
\end{Lemma}

\begin{proof}
Let $[x]_\sim$ be an arbitrary point in $\overline{\mathcal{X}}(\V,k)$.
The seminorm $x$ is diagonalizable with respect to some basis $\be'$ of $\V$.
A similar argument as in the proof of Proposition~\ref{prop -
  fundamental domains} shows that there exists an element $h \in {\rm P}_o$ such that $hx$ lies in $\compap$ (actually $hx$ lies in the closure, taken in the seminorm compactification, of a Weyl chamber with tip $o$).
\end{proof}
\smallskip

The group ${\rm P}_o$ is closed and bounded in ${\rm SL}({\rm V},k)$,
hence compact. Since $\overline{\mathbb{A}}_{\mathbf{e}}$ is compact by Proposition
\ref{subsequences}, the previous Lemma proves that $\overline{\mathcal{X}}({\rm V},k)$
is compact.

\subsubsection{Isotropy groups}
\label{ss - stabilizers}
Let $z$ be a point in $\overline{\mathcal{X}}(\V,k)$ represented by a seminorm $x$ with kernel ${\rm W} \subset \V$. By $\overline{x}$ we denote the norm induced by $x$ on the quotient space $\V/{\rm W}$. By definition, an element $g \in {\rm PGL}(\V,k)$ stabilizes $z$ if and only if one (and hence any) representative $h$ of $g$ in ${\rm GL}(\V,k)$ satisfies $h x \sim x$, i.e. if and only if there exists some $\gamma > 0$ such that
\begin{equation*}
 x (h^{-1} (v)) = \gamma x(v) \text{ for all } v \in \V.
\tag{$\ast$}
\end{equation*}
This is equivalent to saying that $h$ preserves the subspace ${\rm W}$ and that the induced element $\overline{h}$ in ${\rm GL}(\V/{\rm W},k)$ stabilizes the equivalence class of the norm $\overline{x}$ on $\V/{\rm W}$. Hence we find
\[\Stab_{{\rm PGL}(\V,k)}(z) = \{ h \in {\rm GL}(\V,k): h \text{ fixes the subspace }{\rm W} \text{ and }
\overline{h}\overline{x} \sim \overline{x}\} / k^\times.\]

Let us now assume that $z$ is contained in the compactified apartment $\compap$ given by the basis $\be$ of $\V$. Then there are non-negative real numbers $r_0, r_1, \ldots, r_d $ such that
\[x ( \sum_i a_i e_i) = \max_i \{r_i |a_i|\}.\]
The space ${\rm W}$ is generated by all vectors $e_i$ such that $r_i = 0$. We assume that if $r_i$ and $r_j$ are both non-zero, the element $r_j/r_i$ is contained in the value group $|k^\ast|$ of $k$. In this case,
if $h$ stabilizes $z$, we find  that
$ \gamma  = x (h^{-1} e_i) / r_i $ is contained in the value group $|k^\ast|$ of $k$, i.e.
we have $\gamma = |\lambda|$ for some $\lambda \in k^\ast$.
Hence $(\lambda h) x = x$.
Therefore in this case the stabilizer of $z$ in ${\rm PGL}(\V,k)$ is equal to the image of
\[\{ h \in {\rm GL}(\V,k): h \text{ fixes the subspace } {\rm W} \text{ and }
\overline{h}\overline{x} = \overline{x}\}\]
under the natural map from ${\rm GL}(\V,k)$ to ${\rm PGL}(\V,k)$.

\begin{Lemma}
Assume that $z$ is  contained in the closed Weyl chamber $\overline{\mathcal{C}} = \{[x]_\sim \in \compap: x(e_0) \leqslant x(e_1) \leqslant \ldots \leqslant x(e_d) \}$, i.e. using the previous notation we have
$ r_0 \leqslant r_1 \leqslant  \ldots \leqslant r_d $.
Let $d-\mu$ be the index such that $r_{d-\mu} = 0$ and $r_{d-\mu+1} > 0$. (If $z$ is contained in $\ap$, then we put $\mu = d+1$. )
Then the space $W$ is generated by the vectors $e_i$ with $i \leqslant d-\mu$.
We assume as above that $r_j / r_i$ is contained in $|k^\ast|$ if $i>
d-\mu$ and $j > d-\mu$. Writing elements in ${\rm GL}(\V)$ as matrices
with respect to the basis $\be$, we find that ${\rm Stab}_{{\rm PGL}(\V,k)}(z)$ is the image of
\begin{multline*}
  \left\{ \left( \begin{array}{ll}
{\rm A} & {\rm B} \\ 0 &  {\rm D} \end{array} \right) \in {\rm GL}_{d+1}(k) : {\rm D} = (\delta_{ij}) \in {\rm GL}_\mu(k), \right. 
\\
\left. \text{ with } \mid\! {\rm det(D)} \!\mid \, = 1 \text{ and } |\delta_{ij}| \leqslant r_j/r_i \text{ for all } i,j \leqslant \mu. \right\} 
\end{multline*}
in ${\rm PGL}(\V,k)$.
\end{Lemma}

\begin{proof}
This follows directly from the previous considerations combined with Proposition \ref{prop - stabilizers}  which describes the stabilizer groups of norms.\end{proof}

\smallskip
The isotropy groups of the boundary points can also be described in terms of tropical linear algebra, see \cite[Proposition~3.8]{ Wer11}.


\section{Bruhat-Tits theory}
\label{s - Bruhat-Tits general}

We provide now a very short survey of Bruhat-Tits theory.
The main achievement of the latter theory is the existence, for many reductive groups over valued fields, of a combinatorial structure on the rational points; the geometric viewpoint on this is the existence of a strongly transitive action of the group of rational points on a Euclidean building.
Roughly speaking, one half of this theory (the one written in  \cite{BT1a}) is of geometric and combinatorial nature and involves group actions on Euclidean buildings: the existence of a strongly transitive action on such a building is abstractly shown to come from the fact that the involved group can be endowed with the structure of a valued root group datum.
The other half of the theory (the one written in \cite{BT1b}) shows that in many situations, in particular when the valued ground field is local, the group of rational points can be endowed with the structure of a valued root group datum.
This is proved by subtle arguments of descent of the ground field and the main tool for this is provided by group schemes over the ring of integers of the valued ground field.
Though it concentrates on the case when the ground field is local, the survey article \cite{TitsCorvallis} written some decades ago by J.~Tits himself is still very useful. For a very readable introduction covering also the case of a
non-discrete valuation, we recommend the recent text of Rousseau
{\cite{RousseauGrenoble}}.


\subsection{Reductive groups}
\label{ss - Reductive groups}

We introduce a well-known family of algebraic groups which contains most classical groups (i.e., groups which are automorphism groups of suitable bilinear or sesquilinear forms, possibly taking into account an involution, see \cite{Weil} and \cite{KMRT}).
The ground field here is not assumed to be endowed with any absolute value.
The structure theory for rational points is basically due to C.~Chevalley over algebraically closed fields \cite{Bible}, and to A.~Borel and J.~Tits over arbitrary fields \cite{BoTi} (assuming a natural isotropy hypothesis).

\subsubsection{Basic structure results}
\label{sss - structure reductive}
We first need to recall some facts about general linear algebraic groups, up to quoting classical conjugacy theorems and showing how to exhibit a root system in a reductive group.
Useful references are A. Borel's \cite{Borel}, Demazure-Gabriel's \cite{DemazureGabriel} and W.C. Waterhouse's \cite{Waterhouse} books.

\smallskip

{\it Linear algebraic groups}.---~By convention, unless otherwise stated, an ``algebraic group" in what follows means a ``linear algebraic group over some ground field"; being a linear algebraic group amounts to being a smooth affine algebraic group scheme (over a field).
Any algebraic group can be embedded as a closed subgroup of some group ${\rm GL}({\rm V})$ for a suitable vector space over the same ground field (see \cite[3.4]{Waterhouse} for a scheme-theoretic statement and \cite[Prop. 1.12 and Th. 5.1]{Borel} for stronger statements but in a more classical context).

\smallskip

Let ${\rm G}$ be such a group over a field $k$; we will often consider the group ${\rm G}_{k^a}= {\rm G}\otimes_k k^a$ obtained by extension of scalars from $k$ to an algebraic closure $k^a$.

\smallskip

{\it Unipotent and diagonalizable groups}.---~
We say that $g \in {\rm G}(k^a)$ is {\it unipotent} if it is sent to a unipotent matrix in some ({\it a posteriori} any) linear embedding $\varphi : {\rm G}\hookrightarrow {\rm GL}({\rm V})$: this means that $\varphi(g) - {\rm id}_{\rm V}$ is nilpotent.
The group ${\rm G}_{k^a}$ is called {\it unipotent} if so are all its elements; this is equivalent to requiring that the group fixes a vector in any finite-dimensional linear representation as above  \cite[8.3]{Waterhouse}.

The group ${\rm G}$ is said to be a {\it torus} if it is connected and if ${\rm G}_{k^a}$ is \emph{diagonalizable}, which is to say that the algebra of regular functions $\mathcal{O}({\rm G}_{k^a})$ is generated by the characters of ${\rm G}_{k^a}$, i.e., $\mathcal{O}({\rm G}_{k^a}) \simeq k^a[{\rm X}({\rm G}_{k^a})]$ \cite[\S 8]{Borel}.
Here, ${\rm X}({\rm G}_{k^a})$ denotes the finitely generated abelian group of characters ${\rm G}_{k^a}\to {\bf G}_{{\bf m}, k^a}$ and $k^a[{\rm X}({\rm G}_{k^a})]$ is the corresponding group algebra over $k^a$.
A torus ${\rm G}$ defined over $k$ (also called a $k$-{\it torus}) is said to be {\it split over $k$} if the above condition holds over $k$, i.e., if its coordinate ring $\mathcal{O}({\rm G})$ is the group algebra of the abelian group ${\rm X}^*({\rm G}) = {\rm Hom}_{k,{\bf Gr}}({\rm G}, \mathbf{G}_{m,k})$.
In other words, a torus is a connected group of simultaneously diagonalizable matrices in any linear embedding over $k^a$ as above, and it is $k$-split if it is diagonalized in any linear embedding defined over $k$ \cite[\S 7]{Waterhouse}.

\smallskip

{\it Lie algebra and adjoint representation}.---~
One basic tool in studying connected real Lie groups is the Lie algebra of such a group, that is its tangent space at the identity element \cite[3.5]{Borel}.
In the context of algebraic groups, the definition is the same but it is conveniently introduced in a functorial way \cite[\S 12]{Waterhouse}.

\begin{Def}
\label{def - Lie algebra}
Let ${\rm G}$ be a linear algebraic group over a field $k$.
The {\rm Lie algebra} of ${\rm G}$, denoted by $\mathcal{L}({\rm G})$, is the kernel of the natural map ${\rm G}(k[\varepsilon]) \to {\rm G}(k)$, where $k[\varepsilon]$ is the $k$-algebra $k[X]/(X)$ and $\varepsilon$ is the class of $X$; in particular, we have $\varepsilon^2=0$.
\end{Def}

We have $k[\varepsilon] = k \oplus k\varepsilon$ and the natural map above is obtained by applying the functor of points ${\rm G}$ to the map $k[\varepsilon]\to k$ sending $\varepsilon$ to $0$.
The bracket for $\mathcal{L}({\rm G})$ is given by the commutator (group-theoretic) operation \cite[12.2-12.3]{Waterhouse}.

\begin{Example}
For ${\rm G}={\rm GL}({\rm V})$, we have $\mathcal{L}({\rm G}) \simeq {\rm End}({\rm V})$ where ${\rm End}({\rm V})$ denotes the $k$-vector space of all linear endomorphisms of ${\rm V}$.
More precisely, any element of $\mathcal{L}\bigl({\rm GL}({\rm V})\bigr)$ is of the form ${\rm id}_{\rm V} + u \varepsilon$ where $u \in {\rm End}({\rm V})$ is arbitrary.
The previous isomorphism is simply given by $u \mapsto {\rm id}_{\rm V} + u \varepsilon$ and the usual Lie bracket for ${\rm End}({\rm V})$ is recovered thanks to the following computation in ${\rm GL}({\rm V},k[\varepsilon])$:
$[{\rm id}_{\rm V} + u \varepsilon,  {\rm id}_{\rm V} + u' \varepsilon]=  {\rm id}_{\rm V} + (uu'-u'u) \varepsilon$ -- note that the symbol $[.,.]$ on the left hand-side stands for a commutator and that
$({\rm id}_{\rm V} + u \varepsilon)^{-1} = {\rm id}_{\rm V} - u \varepsilon$ for any $u \in {\rm End}({\rm V})$.
\end{Example}

An important tool to classify algebraic groups is the adjoint representation \cite[3.13]{Borel}.

\begin{Def}
\label{def - adjoint representation}
Let ${\rm G}$ be a linear algebraic group over a field $k$.
The {\rm adjoint representation} of ${\rm G}$ is the linear representation ${\rm Ad} : {\rm G} \to {\rm GL}\bigl(\mathcal{L}({\rm G})\bigr)$ defined by
${\rm Ad}(g) = {\rm int}(g) \!\mid_{\mathcal{L}({\rm G})}$ for any $g \in {\rm G}$, where ${\rm int}(g)$ denotes the conjugacy $h \mapsto ghg^{-1}$ -- the restriction makes sense since, for any $k$-algebra ${\rm R}$, both ${\rm G}({\rm R})$ and $\mathcal{L}({\rm G}) \otimes_k {\rm R}$ can be seen as subgroups of ${\rm G}({\rm R}[\varepsilon])$ and the latter one is normal.
\end{Def}

In other words, the adjoint representation is the linear representation provided by differentiating conjugacies at the identity element.

\begin{Example}
For ${\rm G}={\rm SL}({\rm V})$, we have $\mathcal{L}({\rm G}) \simeq \{u \in {\rm End}({\rm V}) : {\rm tr}(u) = 0\}$ and ${\rm Ad}(g).u = gug^{-1}$ for any $g \in {\rm SL}({\rm V})$ and any $u \in \mathcal{L}({\rm G})$.
In this case, we write sometimes $\mathcal{L}({\rm G}) = \mathfrak{sl}({\rm V})$.
\end{Example}

\smallskip

{\it Reductive and semisimple groups}.---~
The starting point for the definition of reductive and semisimple groups consists of the following existence statement \cite[11.21]{Borel}.

\begin{PropDef}
\label{propdef - radicals}
Let ${\rm G}$ be a linear algebraic group over a field $k$.
\begin{itemize}
\item[(i)] There is a unique connected, unipotent, normal subgroup in ${\rm G}_{k^a}$, which is maximal for these properties.
It is called the {\rm unipotent radical} of ${\rm G}$ and is denoted by $\mathcal{R}_u({\rm G})$.
\item[(ii)]  There is a unique connected, solvable, normal subgroup in ${\rm G}_{k^a}$, which is maximal for these properties.
It is called the {\rm radical} of ${\rm G}$ and is denoted by $\mathcal{R}({\rm G})$.
\end{itemize}
\end{PropDef}

The statement for the radical is implied by a finite dimension argument and the fact that the Zariski closure of the product of two connected, normal, solvable subgroups is again connected, normal and solvable.
The unipotent radical is also the unipotent part of the radical: indeed, in a connected solvable group (such as $\mathcal{R}({\rm G})$), the unipotent elements form a closed, connected, normal subgroup \cite[10.3]{Waterhouse}.
Note that by their very definitions, the radical and the unipotent radical depend only on the $k^a$-group ${\rm G}_{k^a}$ and not on the $k$-group ${\rm G}$.

\begin{Def}
\label{def - reductive and ss groups}
Let ${\rm G}$ be a linear algebraic group over a field $k$.
\begin{itemize}
\item[(i)] We say that ${\rm G}$ is {\rm reductive} if we have $\mathcal{R}_u({\rm G}) = \{1 \}$.
\item[(ii)]  We say that ${\rm G}$ is {\rm semisimple} if we have $\mathcal{R}({\rm G}) = \{1 \}$.
\end{itemize}
\end{Def}

\begin{Example}
\label{ex - reductive groups}
For any finite-dimensional $k$-vector space ${\rm V}$, the group ${\rm GL}({\rm V})$ is reductive and ${\rm SL}({\rm V})$ is semisimple.
The groups ${\rm Sp}_{2n}$ and ${\rm SO}(q)$ (for most quadratic forms $q$) are semisimple.
\end{Example}

If, taking into account the ground field $k$, we had used a rational version of the unipotent radical, then we would have obtained a weaker notion of reductivity.
More precisely, it makes sense to introduce the {\it rational unipotent radical}, denoted by $\mathcal{R}_{u,k}({\rm G})$ and contained in $\mathcal{R}_u({\rm G})$, defined to be the unique maximal connected, unipotent  subgroup in ${\rm G}$ {\it defined over $k$}.
Then ${\rm G}$ is called {\it $k$-pseudo-reductive} if we have $\mathcal{R}_{u,k}({\rm G}) = \{1 \}$.
This class of groups is considered in the note \cite{BoTiCRAS}, it is first investigated in some of J.~Tits' lectures (\cite{Tits9192} and \cite{Tits9293}).
A thorough study of pseudo-reductive groups and their classification are written in B.~Conrad, O.~Gabber and G.~Prasad's book \cite{CGP} (an available survey is for instance \cite{RemyBBK}).

\smallskip

{\it In the present paper, we are henceforth interested in reductive groups}.

\smallskip

{\it Parabolic subgroups}.---~
The notion of a parabolic subgroup can be defined for any algebraic group \cite[11.2]{Borel} but it is mostly useful to understand the structure of rational points of reductive groups.

\begin{Def}
\label{def - parabolic subgroup}
Let ${\rm G}$ be a linear algebraic group over a field $k$ and let ${\rm H}$ be a Zariski closed subgroup of $G$.
The subgroup ${\rm H}$ is called {\rm parabolic} if the quotient space ${\rm G}/{\rm H}$ is a complete variety.
\end{Def}

It turns out {\it a posteriori} that for a parabolic subgroup ${\rm H}$, the variety ${\rm G}/{\rm H}$ is actually a projective one; in fact, it can be shown that ${\rm H}$ is a parabolic subgroup if and only if it contains a {\it Borel subgroup}, that is a maximal connected solvable subgroup \cite[11.2]{Borel}.

\begin{Example}
\label{ex - parabolics for GL}
For ${\rm G}={\rm GL}({\rm V})$, the parabolic subgroups are, up to conjugacy, the various groups of upper triangular block matrices (there is one conjugacy class for each ``shape" of such matrices, and these conjugacy classes exhaust all possibilities).
\end{Example}

The completeness of the quotient space ${\rm G}/{\rm H}$ is used to have fixed-points for some subgroup action, which eventually provides conjugacy results as stated below \cite[IV, \S 4, Th. 3.2]{DemazureGabriel}.

\smallskip

{\it Conjugacy theorems}.---~
We finally mention a few results which, among other things, allow one to formulate classification results independent from the choices made to construct the classification data (e.g., the root system -- see \ref{sss - RS and RD} below) \cite[Th. 20.9]{Borel}.

\begin{Thm}
\label{th - conjugacy}
Let ${\rm G}$ be a linear algebraic group over a field $k$.
We assume that ${\rm G}$ is reductive.
\begin{itemize}
\item[(i)] Minimal parabolic $k$-subgroups are conjugate over $k$, that is any two minimal parabolic $k$-subgroups are conjugate by an element of ${\rm G}(k)$.
\item[(ii)]  Accordingly, maximal $k$-split tori are conjugate over $k$.
\end{itemize}
\end{Thm}

For the rational conjugacy of tori, the reductivity assumption can be dropped and simply replaced by a connectedness assumption; this more general result is stated in \cite[C.2]{CGP}.
In the general context of connected groups (instead of reductive ones), one has to replace parabolic subgroups by {\it pseudo-parabolic} ones in order to obtain similar conjugacy results \cite[Th. C.2.5]{CGP}.

\subsubsection{Root system, root datum and root group datum}
\label{sss - RS and RD}
The notion of a root system is studied in detail in \cite[IV]{Lie456}.
It is a combinatorial notion which encodes part of the structure of rational points of semisimple groups.
It also provides a nice uniform way to classify semisimple groups over algebraically closed fields up to isogeny, a striking fact being that the outcome does not depend on the characteristic of the field \cite{Bible}.

In order to state the more precise classification of reductive groups up to isomorphism (over algebraically closed fields, or more generally of split reductive groups), it is necessary to introduce a more subtle notion, namely that of a {\it root datum}:

\begin{Def}
\label{def - root datum}
Let ${\rm X}$ be a finitely generated free abelian group; we denote by ${\rm X}\,\,\check{}$ its ${\bf Z}$-dual and by $\langle \cdot, \cdot \rangle$ the duality bracket.
Let $R$ and $R\,\,\check{}$ be two finite subsets in ${\rm X}$ and ${\rm X}\,\,\check{}$, respectively.
We assume we are given a bijection $\,\,\check{} : \alpha \mapsto \alpha\,\,\check{}$ from $R$ onto $R\,\,\check{}$.
We have thus, for each $\alpha \in R$, endomorphisms

\smallskip
\centerline{
$s_\alpha : x \mapsto x - \langle x, \alpha\,\,\check{} \,\, \rangle \alpha$
\qquad
and
\qquad
$s_\alpha\check{} : x\,\check{}\, \mapsto x\,\,\check{}\, - \langle \alpha, x\,\check{}\,\, \rangle \alpha\,\,\check{}$
}
\smallskip
\noindent of the groups ${\rm X}$ and ${\rm X}\,\,\check{}$, respectively.
The quadruple $\Psi = ({\rm X}, R, {\rm X}\,\,\check{}, R\,\,\check{}\,)$ is said to be a {\rm root datum} if it satisfies the following axioms:
\begin{itemize}
\item[RD 1] For each $\alpha \in R$, we have $\langle \alpha, \alpha\,\,\check{} \,\, \rangle=2$.
\item[RD 2] For each $\alpha \in R$, we have $s_\alpha(R)=R$ and $s_\alpha\check{}\,(R\,\,\check{}\,) = R\,\,\check{}$.
\end{itemize}
\end{Def}

This formulation is taken from \cite{Springer}. The elements of $R$ are called {\rm roots} and the reflections $s_\alpha$ generate a finite group ${\rm W}$ of automorphisms of ${\rm X}$, called the \emph{Weyl group} of $\Psi$.

Let ${\rm Q}$ denote the subgroup of ${\rm X}$ generated by $R$. Up to introducing ${\rm V} = {\rm Q} \otimes_{\bf Z} {\bf R}$ and choosing a suitable ${\rm W}$-invariant scalar product on ${\rm V}$, we can see that $R$ is a root system in the following classical sense:

\begin{Def}
\label{def - root system}
Let ${\rm V}$ be a finite-dimensional real vector space endowed with a scalar product which we denote by $\langle \cdot, \cdot \rangle$. We say that a finite subset $R$ of ${\rm V} - \{0\}$ is a \emph{root system} if it spans ${\rm V}$ and if it satisfies the following two conditions.
\begin{itemize}
\item[RS 1]  To each $\alpha \in R$ is associated a reflection $s_\alpha$ which stabilizes $R$ and switches $\alpha$ and $-\alpha$.
\item[RS 2] For any $\alpha, \beta \in R$, we have $s_\alpha(\beta) - \beta \in \mathbf{Z}\alpha$.
\end{itemize}
\end{Def}

The Weyl group of $\Psi$ is identified with the group of automorphisms of ${\rm V}$ generated by the euclidean reflections $s_\alpha$.

Let $R$ be a root system.
For any subset $\Delta$ in $R$, we denote by $R^+(\Delta)$ the set of roots which can be written as a linear combination of elements of $\Delta$ with non-negative integral coefficients.
We say that $\Delta$ is a {\it basis} for the root system $R$ if it is a basis of ${\rm V}$ and if we have $R = R^+(\Delta) \sqcup R^-(\Delta)$, where $R^-(\Delta) = - R^+(\Delta)$.
Any root system admits a basis and any two bases of a given root system are conjugate under the Weyl group action \cite[VI.1.5, Th. 2]{Lie456}.
When $\Delta$ is a basis of the root system $R$, we say that $R^+(\Delta)$ is a {\it system of positive roots} in $R$; the elements in $\Delta$ are then called {\it simple roots} (with respect to the choice of $\Delta$).
The {\it coroot} associated to $\alpha$ is the linear form $\alpha^\vee$ on ${\rm V}$ defined by $\beta - s_\alpha(\beta) = \alpha^\vee(\beta) \alpha$; in particular, we have $\alpha^\vee(\alpha)=2$.

\begin{Example}
\label{ex - RS of type A}
Here is a well-known concrete construction of the root system of type ${\rm A}_n$.
Let ${\bf R}^{n+1} = \bigoplus_{i=0}^n {\bf R}\varepsilon_i$ be equipped with the standard scalar product, making the basis $(\varepsilon_i)$ orthonormal.
Let us introduce the hyperplane ${\rm V} = \{ \sum_i \lambda_i \varepsilon_i : \sum_i \lambda_i = 0 \}$; we also set $\alpha_{i,j}= \varepsilon_i - \varepsilon_j$ for $i \neq j$.
Then $R = \{\alpha_{i,j} : i \neq j\}$ is a root system in ${\rm V}$ and $\Delta = \{\alpha_{i,i+1} : 0 \leqslant i \leqslant n-1\}$ is a basis of it for which $R^+(\Delta) =  \{\alpha_{i,j} : i < j\}$.
The Weyl group is isomorphic to the symmetric group $\mathcal{S}_{n+1}$; canonical generators leading to a Coxeter presentation are for instance given by transpositions $i \leftrightarrow i+1$.
\end{Example}

Root systems in reductive groups appear as follows.
The restriction of the adjoint representation (Definition~\ref{def - adjoint representation}) to a maximal $k$-split torus ${\rm T}$ is simultaneously diagonalizable over $k$, so that we can write:

\smallskip

$$\mathcal{L}({\rm G}) = \bigoplus_{\varphi \in {\rm X}^*({\rm T})} \mathcal{L}({\rm G})_\varphi$$
where
$$\mathcal{L}({\rm G})_\varphi  = \{v \in \mathcal{L}({\rm G}) : {\rm Ad}(t).v = \varphi(t)v {\rm \ for \ all \ } t \in {\rm T}(k)\}.$$

\smallskip

The normalizer ${\rm N}={\rm N}_{\rm G}({\rm T})$ acts on ${\rm X}^*({\rm T})$ via its action by (algebraic) conjugation on ${\rm T}$, hence it permutes algebraic characters.
The action of the centralizer ${\rm Z}={\rm Z}_{\rm G}({\rm T})$ is trivial, so the group actually acting is the finite quotient ${\rm N}(k)/{\rm Z}(k)$ (finiteness follows from rigidity of tori \cite[7.7]{Waterhouse}, which implies that the identity component ${\rm N}^\circ$ centralizes ${\rm T}$; in fact, we have ${\rm N}^\circ = {\rm Z}$ since centralizers of tori in connected groups are connected).

\smallskip

\centerline{$R = R({\rm T},{\rm G}) = \{\varphi \in {\rm X}^*({\rm T}) : \mathcal{L}({\rm G})_\varphi \neq \{0 \}\}$.}

\smallskip

\noindent It turns out that \cite[Th. 21.6]{Borel}:

\smallskip

\begin{enumerate}
\item the ${\bf R}$-linear span of $R$ is ${\rm V} = {\rm Q} \otimes_{\bf Z} {\bf R}$, where ${\rm Q} \subset {\rm X}^\ast({ \rm T})$ is generated by ${R}$;
\item there exists an ${\rm N}(k)/{\rm Z}(k)$-invariant scalar product ${\rm V}$;
\item the set $R$ is a root system in ${\rm V}$ for this scalar product;
\item the Weyl group ${\rm W}$ of this root system is isomorphic to ${\rm N}(k)/{\rm Z}(k)$.
\end{enumerate}

\noindent Moreover one can go further and introduce a root datum by setting ${\rm X} = {\rm X}^*({\rm T})$ and by taking ${\rm X}\,\,\check{}$ to be the group of all 1-parameter multiplicative subgroups of ${\rm T}$.
The roots $\alpha$ have just been introduced before, but distinguishing the coroots among the cocharacters in ${\rm X}\,\,\check{}$\, is less immediate (over algebraically closed fields or more generally in the split case, they can be defined by means of computation in copies of ${\rm SL}_2$ attached to roots as in Example \ref{ex - root group datum GL(n)} below).
We won't need this but, as already mentioned, in the split case the resulting quadruple $\Psi = ({\rm X}, R, {\rm X}\,\,\check{}, R\,\,\check{}\,)$ characterizes, up to isomorphism, the reductive group we started with (see \cite{SGA3} or \cite[Chap. 9 and 10]{Springer}).

\smallskip

One of the main results of Borel-Tits theory \cite{BoTi} about reductive groups over arbitrary fields is the existence of a very precise combinatorics on the groups of rational points.
The definition of this combinatorial structure -- called a {\it root group datum} -- is given in a purely group-theoretic context.
It is so to speak a collection of subgroups and classes modulo an abstract subgroup ${\rm T}$, all indexed by an abstract root system and subject to relations which generalize and formalize the presentation of ${\rm SL}_n$ (or of any split simply connected simple group) over a field by means of elementary unipotent matrices \cite{Steinberg}.
This combinatorics for the rational points ${\rm G}(k)$ of an isotropic reductive group ${\rm G}$ is indexed by the root system $R({\rm T}, {\rm G})$ with respect to a maximal split torus which we have just introduced; in that case, the abstract group ${\rm T}$ of the root group datum can be chosen to be the group of rational points of the maximal split torus (previously denoted by the same letter!).
More precisely, the axioms of a root group datum are given in the following definition, taken from \cite[6.1]{BT1a}\footnote{Though the notion is taken from \cite{BT1a}, the terminology we use here is not the exact translation of the French ``donn\'ee radicielle" as used in [loc. cit.]: this is because we have already used the terminology ``root datum" in the combinatorial sense of \cite{SGA3}.
Accordingly, we use the notation of \cite{SGA3} instead of that of  \cite{BT1a}, e.g. a root system is denoted by the letter $R$ instead of $\Phi$.}.

\begin{Def}
\label{def - RD}
Let $R$ be a root system and let ${\rm G}$ be a group.
Assume that we are given a system $\bigl( {\rm T}, ({\rm U}_\alpha, {\rm M}_\alpha)_{\alpha \in R}\bigr)$ where ${\rm T}$ and each ${\rm U}_\alpha$ is a subgroup in ${\rm G}$, and each ${\rm M}_\alpha$ is a right congruence class modulo ${\rm T}$.
We say that this system is a {\rm root group datum} of type $R$ for ${\rm G}$ if it satisfies the following axioms:
\begin{itemize}
\item[RGD 1] For each $\alpha \in R$, we have ${\rm U}_\alpha \neq \{1 \}$.
\item[RGD 2] For any $\alpha,\beta \in R$, the commutator group $[{\rm U}_\alpha, {\rm U}_\beta]$ is contained in the group generated by the groups ${\rm U}_\gamma$ indexed by roots $\gamma$ in $R \cap ({\bf Z}_{>0}\alpha + {\bf Z}_{>0}\beta)$.
\item[RGD 3] If both $\alpha$ and $2 \alpha$ belong to $R$, we have ${\rm U}_{2\alpha}\subsetneq {\rm U}_\alpha$.
\item[RGD 4] For each $\alpha \in R$, the class ${\rm M}_\alpha$ satisfies ${\rm U}_{-\alpha} {\bf -} \{1 \}\subset {\rm U}_\alpha {\rm M}_\alpha {\rm U}_\alpha$.
\item[RGD 5] For any $\alpha,\beta \in R$ and each $n \in {\rm M}_\alpha$, we have $n{\rm U}_\beta n^{-1} = {\rm U}_{s_\alpha(\beta)}$.
\item[RGD 6] We have ${\rm T}{\rm U}^+ \cap {\rm U}^- = \{1 \}$, where ${\rm U}^\pm$ is the subgroup generated by the groups ${\rm U}_\alpha$ indexed by the roots $\alpha$ of sign $\pm$.
\end{itemize}
The groups ${\rm U}_\alpha$ are called the {\rm root groups} of the root group datum.
\end{Def}

This list of axioms is probably a bit hard to swallow in one stroke, but the example of ${\rm GL}_n$ can help a lot to have clearer ideas.
We use the notation of Example \ref{ex - RS of type A} (root system of type ${\rm A}_n$).

\begin{Example}
\label{ex - root group datum GL(n)}
Let ${\rm G}= {\rm GL}_{n+1}$ and let ${\rm T}$ be the group of invertible diagonal matrices.
To each root $\alpha_{i,j}$ of the root system $R$ of type ${\rm A}_n$, we attach the subgroup of elementary unipotent matrices ${\rm U}_{i,j} = {\rm U}_{\alpha_{i,j}}= \{ {\rm I}_n + \lambda {\rm E}_{i,j} : \lambda \in k\}$.
We can see easily that ${\rm N}_{\rm G}({\rm T}) = \{$monomial matrices$\}$, that ${\rm Z}_{\rm G}({\rm T}) = {\rm T}$ and finally that ${\rm N}_{\rm G}({\rm T})/{\rm Z}_{\rm G}({\rm T}) \simeq \mathcal{S}_{n+1}$.
Acting by conjugation, the group ${\rm N}_{\rm G}({\rm T})$ permutes the subgroups ${\rm U}_{\alpha_{i,j}}$ and the corresponding action on the indexing roots is nothing else than the action of the Weyl group $\mathcal{S}_{n+1}$ on $R$.
The axioms of a root group datum follow from matrix computation, in particular checking axiom {\rm (RGD4)} can be reduced to the following equality in ${\rm SL}_2$:

\smallskip

\centerline{$
\begin{pmatrix}1&0\\1&1\end{pmatrix} = \begin{pmatrix}1&1\\0&1\end{pmatrix} \left(\begin{array}{lr}0&-1\\1&0\end{array}\right)\begin{pmatrix}1&1\\0&1\end{pmatrix}
$.}
\end{Example}

We can now conclude this subsection by quoting a general result due to A.~Borel and J.~Tits (see  \cite[6.1.3 c)]{BT1a} and \cite{BoTi}).

\begin{Thm}
\label{th - RD in isotropic red gps}
Let ${\rm G}$ be a connected reductive group over a field $k$, which we assume to be $k$-isotropic.
Let ${\rm T}$ be a maximal $k$-split torus in ${\rm G}$, which provides a root system $R = R({\rm T},{\rm G})$.
\begin{itemize}
\item[(i)] For every root $\alpha \in R$ the connected subgroup ${\rm U}_\alpha$ with Lie algebra $\mathcal{L}({\rm G})_\alpha$ is unipotent; moreover it is abelian or two-step nilpotent.
\item[(ii)] The subgroups ${\rm T}(k)$ and ${\rm U}_\alpha(k)$, for $\alpha \in R$, are part of a root group datum of type $R$ in the group of rational points ${\rm G}(k)$.
\end{itemize}
\end{Thm}

Recall that we say that a reductive group is {\it isotropic over $k$}~if it contains a non-central $k$-split torus of positive dimension (the terminology is inspired by the case of orthogonal groups and is compatible with the notion of isotropy for quadratic forms \cite[23.4]{Borel}).
Note finally that the structure of a root group datum implies that (coarser) of a Tits system (also called BN-pair) \cite[IV.2]{Lie456}, which was used by J.~Tits to prove, in a uniform way, the simplicity (modulo center) of the groups of rational points of isotropic simple groups (over sufficiently large fields) \cite{TitsSimple}.

\subsubsection{Valuations on root group data}
\label{sss - DRV}
Bruhat-Tits theory deals with isotropic reductive groups over valued fields.
As for Borel-Tits theory (arbitrary ground field), a substantial part of this theory can also be summed up in combinatorial terms.
This can be done by using the notion of a {\it valuation}~of a root group datum, which formalizes among other things the fact that the valuation of the ground field induces a filtration on each root group.
The definition is taken from  \cite[6.2]{BT1a}.

\begin{Def}
\label{def - V}
Let ${\rm G}$ be an abstract group and let $\bigl( {\rm T}, ({\rm U}_\alpha, {\rm M}_\alpha)_{\alpha \in R}\bigr)$ be a root group datum of type $R$ for it.
A {\rm valuation} of this root group datum is a collection $\varphi = (\varphi_\alpha)_{\alpha \in R}$ of maps $\varphi_\alpha : {\rm U}_\alpha \to {\bf R}\cup \{\infty \}$ satisfying the following axioms.
\begin{itemize}
\item[V0] For each $\alpha \in R$, the image of $\varphi_\alpha$ contains at least three elements.
\item[V1] For each $\alpha \in R$ and each $\ell \in {\bf R}\cup \{\infty \}$, the preimage $\varphi_\alpha^{-1}([\ell;\infty])$ is a subgroup of ${\rm U}_\alpha$, which we denote by ${\rm U}_{\alpha, \ell}$; moreover we require ${\rm U}_{\alpha,\infty} = \{1 \}$.
\item[V2] For each $\alpha \in R$ and each $n \in {\rm M}_\alpha$, the map $u \mapsto \varphi_{-\alpha}(u) - \varphi_\alpha(nun^{-1})$ is constant on the set ${\rm U}_{-\alpha}^* = {\rm U}_{-\alpha}-\{1\}$.
\item[V3] For any $\alpha,\beta \in R$ and $\ell,\ell' \in {\bf R}$ such that $\beta \not\in -{\bf R}_+\alpha$, the commutator group $[{\rm U}_{\alpha,\ell}, {\rm U}_{\beta,\ell'}]$ lies in the group generated by the groups ${\rm U}_{p\alpha + q\beta, p\ell+q\ell'}$ where $p,q \in {\bf Z}_{>0}$ and $p\alpha + q\beta \in R$.
\item[V4] If both $\alpha$ and $2 \alpha$ belong to $R$, the restriction of $2\varphi_\alpha$ to ${\rm U}_{2\alpha}$ is equal to $\varphi_{2\alpha}$.
\item[V5] For $\alpha \in R$, $u \in {\rm U}_\alpha$ and $u', u'' \in {\rm U}_{-\alpha}$ such that $u'uu'' \in {\rm M}_\alpha$, we have $\varphi_{-\alpha}(u') = -\varphi_\alpha(u)$.
\end{itemize}
\end{Def}

The geometric counterpart to this list of technical axioms is the existence, for a group endowed with a valued root group datum, of a Euclidean building (called the {\it Bruhat-Tits building} of the group) on which it acts by isometries with remarkable transitivity properties \cite[\S 7]{BT1a}.
For instance, if the ground field is discretely valued, the corresponding building is simplicial and a fundamental domain for the group action is given by a maximal (poly)simplex, also called an {\it alcove} (in fact, if the ground field is discretely valued, the existence of a valuation on a root group datum can be conveniently replaced by the existence of an affine Tits system \cite[\S 2]{BT1a}).
As already mentioned, the action turns out to be strongly transitive, meaning that the group acts transitively on the inclusions of an alcove in an apartment (Remark \ref{rk - BrT} in \ref{sss - simplicial}).


\subsection{Bruhat-Tits buildings}
\label{ss - Bruhat-Tits buildings}

The purpose of this subsection is to roughly explain how Bruhat-Tits theory attaches a Euclidean building to a suitable reductive group defined over a valued field.
This Bruhat-Tits building comes equipped with a strongly transitive action by the group of rational points, which in turn implies many interesting decompositions of the group.
The latter decompositions are useful for instance to doing harmonic analysis or studying various classes of linear representations of the group.
We roughly explain the descent method used to perform the construction of the Euclidean buildings, and finally mention how some integral models are both one of the main tools and an important outcome of the theory.

\subsubsection{Foldings and gluing}
\label{sss - gluing}
We keep the (connected) semisimple group ${\rm G}$, defined over the (now, complete valued non-Archimedean) field $k$ but from now on, {\it we assume for simplicity that $k$ is a local field (i.e., is locally compact) and we denote by $\omega$ its discrete valuation}, normalized so that $\omega(k^\times)=\mathbf{Z}$. Hence $\omega(\cdot) = - {\rm log}_q |\cdot|$, where $q >1$ is a generator of the discrete group $|k^\times|$.

We also assume that ${\rm G}$ contains a $k$-split torus of positive dimension: this is an isotropy assumption over $k$ already introduced at the end of \ref{sss - RS and RD} (in this situation, this algebraic condition is equivalent to the topological condition that the group of rational points ${\rm G}(k)$ is non-compact \cite{PrasadSMF}).
In order to associate to ${\rm G}$ a Euclidean building on which ${\rm G}(k)$ acts strongly transitively, according to \cite{TitsCorvallis} we need two things:

\begin{enumerate}
\item a model, say $\Sigma$, for the apartments;
\item a way to glue many copies of $\Sigma$ altogether in such a way that they will satisfy the incidence axioms of a building (\ref{sss - simplicial}).
\end{enumerate}

{\it Model for the apartment}.---~
References for what follows are \cite[\S 1]{TitsCorvallis} or \cite[Chapter I]{La}.
Let ${\rm T}$ be a maximal $k$-split torus in ${\rm G}$ and let ${\rm X}_*({\rm T})$ denote its group of 1-parameter subgroups (or {\it cocharacters}).
As a first step, we set $\Sigma_{\rm vect}= {\rm X}_{*}({\rm T}) \otimes_{\bf Z} {\bf R}$.

\begin{Prop}
\label{prop-affine_apartment}
There exists an affine space $\Sigma$ with underlying vector space $\Sigma_{\rm vect}$, equipped with an action by affine transformations  $\nu : {\rm N}(k) =  {\rm N}_{\rm G}({\rm T})(k) \to {\rm Aff}(\Sigma)$ and having the following properties.
\begin{itemize}
\item[(i)] There is a scalar product on $\Sigma$ such that $\nu\bigl( {\rm N}(k) \bigr)$ is an affine reflection group.
\item[(ii)] The vectorial part of this group is the Weyl group of the root system $R = R({\rm T}, {\rm G})$.
\item[(iii)] The translation (normal) subgroup acts cocompactly on $\Sigma$, it is equal to $\nu\bigl({\rm Z}(k)\bigr)$ and the vector $\nu(z)$ attached to an element $z \in {\rm Z}(k)$ is defined by $\chi\bigl(\nu(z)\bigr) = -\omega \bigl( \chi(z) \bigr)$ for any $\chi \in {\rm X}^*({\rm T})$.
\end{itemize}
\end{Prop}

If we go back to the example of ${\rm GL}({\rm V})$ acting by precomposition on the space of classes of norms $ \mathcal{X}({\rm V},k)$ as described in \ref{ss - SL(n) Bruhat-Tits}, we can see the previous statement as a generalization of the fact, mentioned in \ref{sss - actions on GI}, that for any basis $\mathbf{e}$ of ${\rm V}$, the group ${\rm N}_{\rm e}$ of monomial matrices with respect to ${\bf e}$ acts on the apartment $\mathbb{A}_{\bf e}$ as $\mathcal{S}_d \ltimes {\bf Z}^d$ where $d = {\rm dim}({\rm V})$.

\smallskip

{\it Filtrations and gluing}.---~
Still for this special case, we saw (Proposition~ref{prop - folding}) that any elementary unipotent matrix $u_{ij}(\lambda) = {\rm I}_d + \lambda {\rm E}_{ij}$ fixes pointwise a closed half-apartment in $\mathbb{A}_{\bf e}$ bounded by a hyperplane of the form $\{c_i - c_j =$ constant$\}$ (the constant depends on the valuation $\omega(\lambda)$ of the additive parameter $\lambda$), the rest of the apartment $\mathbb{A}_{\bf e}$ associated to ${\bf e}$ being ``folded" away from $\mathbb{A}_{\bf e}$.

In order to construct the Bruhat-Tits building in the general case, the gluing equivalence will impose this folding action for unipotent elements in root groups; this will be done by taking into account the ``valuation" of the unipotent element under consideration.
What formalizes this is the previous notion of a valuation for a root group datum (Definition \ref{def - V}), which provides a filtration on each root group.
For further details, we refer to the motivations given in \cite[1.1-1.4]{TitsCorvallis}.
It is not straightforward to perform this in general, but it can be done quite explicitly when the group ${\rm G}$ is {\it split} over $k$ (i.e., when it contains a maximal torus which is $k$-split).
For the general case, one has to go to a (finite, separable) extension of the ground field splitting ${\rm G}$ and then to use subtle descent arguments.
The main difficulty for the descent step is to handle at the same time Galois actions on the split group and on its ``split" building in order to descend the ground field both for the valuation of the root group datum and at the geometric level (see \ref{sss - descent and functoriality} for slightly more details).

Let us provisionally assume that ${\rm G}$ is split over $k$.
Then each root group ${\rm U}_\alpha (k)$ is isomorphic to the additive group of $k$ and for any such group ${\rm U}_\alpha(k)$ we can use the valuation of $k$ to define a decreasing filtration
$\{{\rm U}_{\alpha}(k)_\ell \}_{\ell \in {\bf Z}}$ satisfying:

\smallskip

\begin{center}
$\bigcup_{\ell \in {\bf Z}} {\rm U}_{\alpha}(k)_\ell = {\rm
  U}_\alpha(k) \hskip3mm {\rm and} \hskip3mm \bigcap_{\ell \in {\bf Z}} {\rm U}_{\alpha}(k)_\ell  = \{1 \},$
\end{center}

\smallskip

\noindent and further compatibilities, namely the axioms of a valuation (Definition~\ref{def - V}) for the root group datum structure on ${\rm G}(k)$ given by Borel-Tits theory (Theorem~\ref{th - RD in isotropic red gps}) -- the latter root group datum structure in the split case is easy to obtain by means of Chevalley bases \cite{Steinberg} (see remark below).
For instance, in the case of the general linear group, this can be merely done by using the parameterizations

\smallskip

\centerline{$(k,+) \simeq {\rm U}_{\alpha_{i,j}}(k) = \{{\rm id}+ \lambda {\rm E}_{i,j} : \lambda \in k\}$.}

\smallskip
\begin{Remark}
\label{rk-filtrations_split}
Let us be slightly more precise here. For a split group ${\rm G}$, each root group ${\rm U}_\alpha$ is $k$-isomorphic to the additive group $\mathbf{G}_a$, and the choice of a Chevalley basis of ${\rm Lie}({\rm G})$ determines a set of isomorphisms $\{p_\alpha : {\rm U}_\alpha \rightarrow \mathbf{G}_a\}_{\alpha \in R}$. It is easily checked that the collection of maps $$\xymatrix{\varphi_\alpha : {\rm U}_\alpha(k) \ar@{->}[r]^{p_\alpha} & \mathbf{G}_a(k) \ar@{->}[r]^\omega & \mathbf{R}}$$ defines a valuation on the root group datum $({\rm T}(k),({\rm U}_\alpha(k),{\rm M}_\alpha))$.

For each $\ell \in \mathbf{R}$, the condition $|p_{\alpha}| \leqslant q^{-s}$ defines an \emph{affinoid} subgroup ${\rm U}_{\alpha,s}$ in ${\rm U}_{\alpha}^{\rm an}$ such that ${\rm U}_{\alpha}(k)_\ell = {\rm U}_{\alpha,s}(k)$ for any $s \in (\ell-1, \ell]$. The latter identity holds after replacement of $k$ by any finite extension $k'$, \emph{as long as we normalize the valuation of $k'$ in such a way that is extends the valuation on $k$}. This shows that Bruhat-Tits filtrations on root groups, in the split case at this stage, comes from a decreasing, exhaustive and separated filtration of ${\rm U}_{\alpha}^{\rm an}$ by affinoid subgroups $\{{\rm U}_{\alpha,s}\}_{s \in \mathbf{R}}$.
\end{Remark}
\smallskip

Let us consider again the apartment $\Sigma$ with underlying vector space $\Sigma_{\rm vect}= {\rm X}_*({\rm T}) \otimes_{\bf Z} {\bf R}$.
We are interested in the affine linear forms $\alpha + \ell$ ($\alpha \in R$, $\ell \in {\bf Z}$).
We fix an origin, say $o$, such that $(\alpha + 0)(o) = 0$ for any root $\alpha \in R$.
We have ``level sets" $\{\alpha + \ell = 0 \}$ and ``positive half-spaces" $\{\alpha + \ell \geqslant 0 \}$ bounded by them.

For each $x \in \Sigma$, we set ${\rm N}_x = {\rm Stab}_{{\rm G}(k)}(x)$ (using the action $\nu$ of Proposition~\ref{prop-affine_apartment})
and for each root $\alpha$ we denote by ${\rm U}_{\alpha}(k)_x$ the biggest subgroup ${\rm U}_{\alpha}(k)_\ell$ such that $x \in \{\alpha + \ell \geqslant 0\}$ (i.e. $\ell$ is minimal for the latter property).
At last, we define ${\rm P}_x$ to be the subgroup of ${\rm G}(k)$ generated by ${\rm N}_x$ and by $\{{\rm U}_{\alpha}(k)_x \}_{\alpha \in R}$.
We are now in position to define a binary relation, say $\sim$, on ${\rm G}(k) \times \Sigma$ by:

\smallskip

\centerline{$(g,x) \sim (h,y)$ \quad $\Longleftrightarrow$ \quad there exists $n \in {\rm N}_{\rm G}(T)(k)$ such that $y=\nu(n).x$ and $g^{-1}hn \in {\rm P}_x$.}

\smallskip

{\it Construction of the Bruhat-Tits buildings}.---~
This relation is exactly what is needed in order to glue together copies of $\Sigma$ and to finally obtain the desired Euclidean building.

\begin{Thm}
\label{th - gluing}
The relation $\sim$ is an equivalence relation on the set ${\rm G}(k) \times \Sigma$ and the quotient space $\displaystyle {\mathcal{B}} = {\mathcal{B}}({\rm G},k) = \left({\rm G}(k) \times \Sigma\right) / \sim$ is a Euclidean building whose apartments are isomorphic to $\Sigma$ and whose Weyl group is the affine reflection group $W = \nu\bigl( {\rm N}(k) \bigr)$.
Moreover the ${\rm G}(k)$-action by left multiplication on the first factor of ${\rm G}(k) \times \Sigma$ induces an action of ${\rm G}(k)$ by isometries on ${\mathcal{B}}({\rm G},k)$.
\end{Thm}

{\it Notation}.---~ According to Definition \ref{defi - non-simplicial building}, copies of $\Sigma$ in $\mathcal{B}({\rm G},k)$  are called \emph{apartments}; they are the maximal flat (i.e., euclidean) subspaces. Thanks to ${\rm G}(k)$-conjugacy of maximal split tori \ref{th - conjugacy}, apartments of  $\mathcal{B}({\rm G},k)$ are in bijection with maximal split tori of ${\rm G}$. Therefore, we will speak of the \emph{apartment of a maximal split torus} ${\rm S}$ of ${\rm G}$ and write ${\rm A}({\rm S},k)$. By construction, this is an affine space under the $\mathbf{R}$-vector space ${\rm Hom}_{\bf Ab}({\rm X}^*({\rm S}),\mathbf{R})$.

\begin{proof}[Reference for the proof]
As already explained, the difficulty is to check the axioms of a valuation (Def \ref{def - V}) for a suitable choice of filtrations on the root groups of a Borel-Tits root group datum (Th. \ref{th - RD in isotropic red gps}). Indeed, the definition of the equivalence relation $\sim$, hence the construction of a suitable Euclidean building, for a valued root group datum can be done in this purely abstract context \cite[\S 7]{BT1a}.
The existence of a valued root group datum for reductive groups over suitable valued (not necessarily complete) fields was announced in \cite[6.2.3 c)]{BT1a} and was finally settled in the second IH\'ES paper (1984) by F.~Bruhat and J.~Tits \cite[Introduction and Th. 5.1.20]{BT1b}.
\end{proof}

One way to understand the gluing equivalence relation $\sim$ is to see that it prescribes stabilizers.
Actually, it can eventually be proved that {\it a posteriori} we have:

\smallskip

\centerline{$\Sigma^{{\rm U}_{\alpha,\ell}(k)} = \{\alpha + \ell \geqslant 0 \}$ \quad and \quad ${\rm Stab}_{{\rm G}(k)}(x) = {\rm P}_x$ for any $x \in {\mathcal{B}}$.}

\smallskip

A more formal way to state the result is to say that to each valued root group datum on a group is associated a Euclidean building, which can be obtained by a gluing equivalence relation defined as above \cite[\S 7]{BT1a}.

\begin{Example}
\label{SL-building}
In the case when ${\rm G}={\rm SL}({\rm V})$, it can be checked that the building obtained by the above method is equivariantly isomorphic to the Goldman-Iwahori space $\mathcal{X}(\V,k)$ \cite[10.2]{BT1a}.
\end{Example}

\subsubsection{Descent and functoriality}
\label{sss - descent and functoriality}
Suitable filtrations on root groups so that an equivalence relation $\sim$ as above can be defined do not always exist.
Moreover, even when things go well, the way to construct the Bruhat-Tits building is not by first exhibiting a valuation on the root group datum given by Borel-Tits theory
and then by using the gluing relation $\sim$.
As usual in algebraic group theory, one has first to deal with the split case, and then to apply various and difficult arguments of descent of the ground field. Bruhat and Tits used a two-step descent, allowing a fine description of smooth integral models of the
group associated with facets. A one-step descent was introduced by Rousseau in
his thesis {\cite{RousseauOrsay}}, whose validity in full generality now
follows from recent work connected to Tits' Center Conjecture
({\cite{Struyve}}).

\smallskip

{\it Galois actions}.---~
More precisely, one has to find a suitable (finite) Galois extension $k'/k$ such that ${\rm G}$ splits over $k'$  (or, at least, {\it quasi-splits} over $k'$, i.e. admits a Borel subgroup defined over $k'$) and, which is much more delicate, which enables one:

\begin{enumerate}
\item to define a ${\rm Gal}(k'/k)$-action by isometries on the ``(quasi)-split" building ${\mathcal{B}}({\rm G},k')$;
\item to check that a building for ${\rm G}(k)$ lies in the Galois fixed point set ${\mathcal{B}}({\rm G},k')^{{\rm Gal}(k'/k)}$.
\end{enumerate}

\noindent Similarly, the group ${\rm G}(k')$ of course admits a ${\rm Gal}(k'/k)$-action.

\begin{Remark}
\label{rk - enough fixed points}
Recall that, by completeness and non-positive curvature, once step 1 is settled we know that we have sufficiently many Galois-fixed points in ${\mathcal{B}}({\rm G},k')$ (see the discussion of the Bruhat-Tits fixed point theorem in \ref{sss - geometry of buildings}).
\end{Remark}

F.~Bruhat and J.~Tits found a uniform procedure to deal with various situations of that kind.
The procedure described in \cite[9.2]{BT1a} formalizes, in abstract terms of buildings and group combinatorics, how to exhibit a valued root group datum structure (resp. a Euclidean building structure) on a subgroup of a bigger group with a valued root group datum (resp. on a subspace of the associated Bruhat-Tits building).
The main result \cite[Th. 9.2.10]{BT1a} says that under some sufficient conditions, the restriction of the valuation to a given sub-root group datum ``descends" to a valuation and its associated Bruhat-Tits building is the given subspace.
These sufficient conditions are designed to apply to subgroups and convex subspaces obtained as fixed-points of ``twists" by Galois actions (and they can also be applied to non-Galois twists ``\`a la Ree-Suzuki").

\smallskip

{\it Two descent steps}.---~
As already mentioned, this needn't work over an arbitrary valued field $k$ (even when $k$ is complete).
Moreover F.~Bruhat and J.~Tits do not perform the descent in one stroke, they have to argue by a two step descent.

The first step is the so-called {\it quasi-split} descent \cite[\S 4]{BT1b}.
It consists in dealing with field extensions splitting an initially quasi-split reductive group.
The Galois twists here (of the ambient group and building) are shown, by quite concrete arguments, to fit in the context of \cite[9.2]{BT1a} mentioned above.
This is possible thanks to a deep understanding of quasi-split groups: they can even be handled via a presentation (see \cite{Steinberg} and \cite[Appendice]{BT1b}).
In fact, the largest part of the chapter about the quasi-split descent \cite[\S 4]{BT1b} is dedicated to another topic which will be presented below (\ref{sss - models}), namely the construction of suitable integral models (i.e. group schemes over $k^\circ$  with generic fiber ${\rm G}$) defined by geometric conditions involving bounded subsets in the building.
The method chosen by F.~Bruhat and J.~Tits to obtain these integral models is by using a linear representation of ${\rm G}$ whose underlying vector space contains a suitable $k^\circ$-lattice, but they mention themselves that this could be done by Weil's techniques of group chunks.
Since then, thanks to the developments of N\'eron model techniques \cite{BLR}, this alternative method has been written down \cite{La}.

The second step is the so-called  {\it \'etale} descent \cite[\S 5]{BT1b}.
By definition, an \'etale extension, in the discretely valued case (to which we stick here), is unramified with separable residual extension; let us denote by $k^{\rm sh}$ the maximal \'etale extension of $k$.
This descent step consists in considering situations where the semisimple $k$-group ${\rm G}$ is such that ${\rm G} \otimes_k k^{\rm sh}$ is quasi-split (so that, by the first step, we already have a valued root group datum and a Bruhat-Tits building for ${\rm G}(k^{\rm sh})$, together with integral structures).
Checking that this fits in the geometric and combinatorial formalism of \cite[9.2]{BT1a} is more difficult in that case.
In fact, this is the place where the integral models over the valuation ring $k^\circ$ are used, in order to find a suitable torus in ${\rm G}$ which become maximal split in ${\rm G} \otimes_k k'$ for  some \'etale extension $k'$ of $k$ \cite[Cor. 5.1.12]{BT1b}.

\begin{Remark}
\label{rk-filtrations}
In the split case, we have noticed that the Bruhat-Tits filtrations on
rational points of root groups come from filtrations by affinoid
subgroups (\ref{rk-filtrations_split}). This fact holds in general and can be checked as follows: let $k'/k$ be a finite Galois
extension splitting ${\rm G}$ and consider a maximal torus ${\rm T}$
of ${\rm G}$ which splits over $k'$ and contains a maximal split torus
${\rm S}$. The canonical projection ${\rm X}^*({\rm T} \otimes_k k')
\rightarrow {\rm X}^*({\rm S} \otimes_k k') \tilde{=} {\rm X}^*({\rm S})$
induces a surjective map $$p: R({\rm T}
\otimes_k k',{\rm G} \otimes_k k') \longrightarrow R({\rm S},{\rm G}) \cup \{0\}$$ and there is a natural $k'$-isomorphism
$$\prod_{\beta \in p^{-1}(\alpha)} {\rm U}_\beta \times \prod_{\beta \in p^{-1}(2\alpha)} {\rm U}_\beta \simeq {\rm U}_\alpha \otimes_k k'$$ for any ordering of the factor.

\smallskip
\emph{A posteriori}, Bruhat-Tits two-step descent proves that any
maximal split torus ${\rm S}$ of ${\rm G}$ is contained in a maximal
torus ${\rm T}$ which splits over a finite Galois extension $k'/k$
such that ${\rm Gal}(k'/k)$ fixes a point in the apartment of
${\rm T} \otimes_k k'$ in $\mathcal{B}({\rm G},k')$. \emph{If the
valuation on $k'$ is normalized in such a way that it extends the
valuation on $k$}, then, for any
$\ell \in \mathbf{R}$, the affinoid
subgroup $$\prod_{\beta \in p^{-1}(\alpha)} {\rm U}_{\beta,\ell}
\times \prod_{\beta \in p^{-1}(2\alpha)} {\rm U}_{\beta,2\ell}$$ of
the left hand side corresponds to an affinoid subgroup of the
right hand side which does not depend on the ordering of the factors
and is preserved by the natural action of ${\rm Gal}(k'|k)$; this
can be checked by using calculations in \cite[6.1]{BT1a} at the level
of $k''$ points, for any finite extension $k''/k'$. By Galois
descent, we obtain an affinoid subgroup ${\rm U}_{\alpha,\ell}$ of
${\rm U}_{\alpha}^{\rm an}$ such that $${\rm U}_{\alpha,\ell}(k) = {\rm
  U}_\alpha(k) \cap \left(\prod_{\beta \in p^{-1}(\alpha)} {\rm U}_{\beta,\ell}(k')
\times \prod_{\beta \in p^{-1}(2\alpha)}{\rm U}_{\beta,2\ell}(k') \right).$$
By \cite[5.1.16 and 5.1.20]{BT1b}, the filtrations $\{{\rm U}_{\alpha,\ell}(k)\}_{\ell
  \in \mathbf{R}}$ are induced by a valuation on the root group datum $\left({\rm
  S}(k), \{{\rm U}_\alpha (k)\} \right)$.
\end{Remark}

Let us finish by mentioning why this two-step strategy is well-adapted to the case we are interested in, namely that of a semisimple group ${\rm G}$ defined over a complete, discretely valued field $k$ with perfect residue field $\widetilde{k}$: thanks to a result of R.~Steinberg's \cite[III, 2.3]{SerreGalois}, such a group is known to quasi-split over $k^{\rm sh}$.
Compactifications of Bruhat-Tits buildings fit in this more specific context for ${\rm G}$ and $k$.
Indeed, the Bruhat-Tits building ${\mathcal{B}}({\rm G},k)$ is locally compact if and only if so is $k$, see the discussion of the local structure of buildings below (\ref{sss - models}).
Note finally that the terminology ``henselian" used in \cite{BT1b} is a well-known algebraic generalization of ``complete" (the latter ``analytic" condition is the only one we consider seriously here, since we will use Berkovich geometry).

\smallskip

{\it Existence of Bruhat-Tits buildings}.---~
Here is at last a general statement on existence of Bruhat-Tits buildings which will be enough for our purposes; this result was announced in \cite[6.2.3 c)]{BT1a} and is implied by  \cite[Th. 5.1.20]{BT1b}.

\begin{Thm}
\label{th - existence of BT buildings}
Assume that $k$ is complete, discretely valued, with perfect residue field. The root group datum of ${\rm G}(k)$ associated with a split maximal torus admits a valuation satisfying the conditions of Definition \ref{def - V}.  
\end{Thm}

Let us also give now an example illustrating both the statement of the theorem and the general geometric approach characterizing Bruhat-Tits theory.

\begin{Example}
\label{ex - hermitian trees}
Let $h$ be a Hermitian form of index $1$ in three variables, say on the vector space ${\rm V} \simeq k^3$.
We assume that $h$ splits over a quadratic extension, say $E/k$, so that ${\rm SU}({\rm V},h)$ is isomorphic to ${\rm SL}_3$ over $E$, and we denote ${\rm Gal}(E/k) = \{1;\sigma \}$.
Then the building of ${\rm SU}({\rm V},h)$ can be seen as the set of fixed points for a suitable action of the Galois involution $\sigma$ on the $2$-dimensional Bruhat-Tits building of type $\widetilde {\rm A}_2$ associated to ${\rm V} \otimes_k E$ as in \ref{ss - SL(n) Bruhat-Tits}.
If $k$ is local and if $q$ denotes the cardinality of the residue field, then the Euclidean building ${\mathcal{B}}({\rm SU}({\rm V},h),k)$ is a locally finite tree: indeed, it is a Euclidean building of dimension $1$ because the $k$-rank of ${\rm SU}({\rm V},h)$, i.e. the dimension of maximal $k$-split tori, is $1$.
The tree is homogeneous of valency $1+q$ when $E/k$ is ramified, in which case the type of the group is {\rm C}-{\rm BC}${}_1$ in Tits' classification \cite[p. 60, last line]{TitsCorvallis}.
The tree is semi-homogeneous of valencies $1+q$ and $1+q^3$ when $E/k$ is unramified, and then the type is ${}^2\! {\rm A}_2'$ \cite[p. 63, line 2]{TitsCorvallis}.
For the computation of the valencies, we refer to \ref{sss - models} below.
\end{Example}


\smallskip

{\it Functoriality}.---~
For our purpose (i.e. embedding of Bruhat-Tits buildings in analytic spaces and subsequent compactifications), the existence statement is not sufficient.
We need a stronger result than the mere existence; in fact, we also need a good behavior of the building with respect to field extensions.

\begin{Thm}
\label{th - functoriality of BT buildings}
Whenever $k$ is complete, discretely valued, with perfect residue field, the Bruhat-Tits building ${\mathcal{B}}({\rm G},K)$ depends functorially on the non-Archimedean extension $K$ of $k$.
\end{Thm}

More precisely, let us denote by ${\rm G}-\mathbf{Sets}$ the category whose objets are pairs $(K/k,{\rm X})$, where $K/k$ is a non-Archimedean extension and ${\rm X}$ is a topological space endowed with a continuous action of ${\rm G}(K)$, and arrows $(K/k,{\rm X}) \rightarrow (K'/k,{\rm X}')$ are pairs $(\iota, f)$, where $\iota$ is an isometric embedding of $K$ into $K'$ and $f$ is a ${\rm G}(K)$-equivariant and continous map from ${\rm X}$ to ${\rm X}'$. We see the building of ${\rm G}$ as a \emph{section} $\mathcal{B}({\rm G}, -)$ of the forgetful functor $${\rm G}-\mathbf{Sets} \longrightarrow \left( \begin{array}{c} {\rm non-Archimedean} \\ {\rm extensions }\  K/k \end{array} \right).$$

\begin{remark}[Reference]
It is explained in \cite[1.3.4]{RTW1} how to deduce this from the general theory.
\end{remark}

One word of caution is in order here. If $k'/k$ is a Galois extension, then there is a natural action of ${\rm Gal}(k'/k)$ on $\mathcal{B}({\rm G},k')$ by functoriality and the smaller building $\mathcal{B}({\rm G},k)$ is contained in the Galois-fixed point set in ${\mathcal{B}}({\rm G},k')$. In general, this inclusion is strict, even when the group is split
\cite[III]{RousseauOrsay} (see also {\bf 5.2}). However, one can show that there is equality if the extension $k'/k$ is \emph{tamely ramified} [{\bf loc. cit.}] and \cite{Prasad}.

\vskip2mm We will need to have more precise information about the behavior of apartments. As above, we assume that $k$ is complete, discretely valued and with perfect residue field.

\vskip2mm
\begin{Def} Let ${\rm T}$ be a maximal torus of ${\rm G}$ and let $k_{\rm T}$ be the minimal Galois extension of $k$ (in some fixed algebraic closure) which splits ${\rm T}$. We denote by $k_{\rm T}^{\rm ur}$ the maximal unramified extension of $k$ in $k_{\rm T}$. 

The torus ${\rm T}$ is \emph{well-adjusted} if the maximal split subtori of ${\rm T}$ and ${\rm T} \otimes_k k_{\rm T}^{\rm ur}$ are maximal split tori of ${\rm G}$ and ${\rm G} \otimes_k k_{\rm T}^{\rm ur}$.
\end{Def} 

\vskip2mm
\begin{Lemma} 
\label{functor-apartment} 1. Every maximal split torus ${\rm S}$ of ${\rm G}$ is contained in a well-adjusted maximal torus ${\rm T}$. \\
\indent 2. Assume that ${\rm S}$ and ${\rm T}$ are as above, and let $K/k$ be any non-Archimedean field extension which splits ${\rm T}$. The embedding $\mathcal{B}({\rm G},k) \hookrightarrow \mathcal{B}({\rm G},{\rm K})$ maps ${\rm A}({\rm S},k)$ into ${\rm A}({\rm T},{\rm K})$. 
\end{Lemma}

\begin{proof}
1. For each unramified finite Galois extension $k'/k$, we can find a torus ${\rm S}' \subset {\rm G}$ which contains ${\rm S}$ and such that ${\rm S}' \otimes_k k'$ is a maximal split torus of ${\rm G} \otimes_k k'$ \cite[Corollaire 5.1.12]{BT1b}. We choose a pair $(k',{\rm S}')$ such that the rank of ${\rm S}'$ is maximal, equal to the relative rank of ${\rm G} \otimes_k k^{\rm ur}$; this means that ${\rm S}' \otimes_k k''$ is a maximal split torus of ${\rm G} \otimes_k k''$ for any unramified extension $k''/k$ containing $k'$.

The centralizer of ${\rm S}' \otimes_k k'$ in ${\rm G} \otimes_k k'$ is a maximal torus of ${\rm G} \otimes_k k'$, hence ${\rm T} = {\rm Z}({\rm S}')$ is a maximal torus of ${\rm G}$. By construction, ${\rm S}'$ splits over $k_{\rm T}^{\rm ur}$ and ${\rm S'} \otimes_k k_{\rm T}^{\rm ur}$ is a maximal split torus of ${\rm G} \otimes_k k_{\rm T}^{\rm ur}$. Since ${\rm S} \subset {\rm S}'$, this proves that ${\rm T}$ is well-adjusted.

2. We keep the same notation as above. The extension $K/k$ contains $k_{\rm T}$, hence it is enough by functoriality to check that the embedding $\mathcal{B}({\rm G},k) \hookrightarrow \mathcal{B}({\rm G},k_{\rm T})$ maps ${\rm A}({\rm S},k)$ into ${\rm A}({\rm T},k_{\rm T})$.

Let us consider the embeddings $$\mathcal{B}({\rm G},k) \hookrightarrow \mathcal{B}({\rm G},k_{\rm T}^{\rm ur}) \hookrightarrow \mathcal{B}({\rm G},k_{\rm T}).$$ The first one maps ${\rm A}({\rm S},k)$ into ${\rm A}({\rm S}',k_{\rm T}^{\rm ur})$  by \ By \cite[Proposition 5.1.14]{BT1b} and the second one maps ${\rm A}({\rm S}',k_{\rm T}^{\rm ur})$ into ${\rm A}({\rm T},k_{\rm T})$ by \cite[Th\'eor\`eme 2.5.6]{RousseauOrsay}, hence their composite has the required property. \end{proof}

\subsubsection{Compact open subgroups and integral structures}
\label{sss - models}
In what follows, we maintain the previous assumptions, in particular the group ${\rm G}$ is semisimple and $k$-isotropic.
The building ${\mathcal{B}}({\rm G}, k)$ admits a strongly transitive ${\rm G}(k)$-action by isometries.
Moreover it is a {\it labelled} simplicial complex in the sense that, if $d$ denotes the number of codimension 1 facets (called {\it panels}) in the closure of a given alcove, we can choose $d$ colors and assign one of them to each panel in  ${\mathcal{B}}({\rm G}, k)$ so that each color appears exactly once in the closure of each alcove.
For some questions, it is convenient to restrict oneself to the finite index subgroup ${\rm G}(k)^\bullet$ consisting of the color-preserving (or {\it type-preserving}) isometries in ${\rm G}(k)$.

\smallskip

{\it Compact open subgroups}.---~
For any facet $F \subset {\mathcal{B}}({\rm G},k)$ we denote by ${\rm P}_F$ the stabilizer ${\rm Stab}_{{\rm G}(k)}(F)$: it is a bounded subgroup of ${\rm G}(k)$ and when $k$ is local, it is a compact, open subgroup.
It follows from the Bruhat-Tits fixed point theorem (\ref{sss - geometry of buildings}) that the conjugacy classes of maximal compact subgroups in ${\rm G}(k)^\bullet$ are in one-to-one correspondence with the vertices in the closure of a given alcove.
The fact that there are usually several conjugacy classes of maximal compact subgroups in ${\rm G}(k)$ makes harmonic analysis more delicate than in the classical case of real Lie groups.
Still, for instance thanks to the notion of a special vertex, many achievements can also be obtained in the non-Archimedean case \cite{Macdonald}.
Recall that a point $x \in {\mathcal{B}}({\rm G},k)$ is called {\it special} if for any apartment $\mathbb{A}$ containing $x$, the stabilizer of $x$ in the affine Weyl group is the full vectorial part of this affine reflection group, i.e. is isomorphic to the (spherical) Weyl group of the root system $R$ of ${\rm G}$ over $k$.

\smallskip

{\it Integral models for some stabilizers}.---~
In what follows, we are more interested in algebraic properties of compact open subgroups obtained as facet stabilizers.
The following statement is explained in \cite[5.1.9]{BT1b}.

\begin{Thm}
\label{th - integral structures}
For any facet $F \subset {\mathcal{B}}({\rm G},k)$ there exists a smooth $k^\circ$-group scheme $\mathcal{G}_F$ with generic fiber ${\rm G}$ such that $\mathcal{G}_F(k^\circ) = {\rm P}_F$.
\end{Thm}

As already mentioned, the point of view of group schemes over $k^\circ$ in Bruhat-Tits theory is not only an important tool to perform the descent, but it is also an important outcome of the theory.
Here is an example.
The ``best" structure {\it a priori} available for a facet stabilizer is only of topological nature (and even for this, we have to assume that $k$ is locally compact).
The above models over $k^\circ$ provide an algebraic point of view on these groups, which allows one to define a filtration on them leading to the computation of some cohomology groups of great interest for the congruence subgroup problem, see for instance \cite{PraRag1} and \cite{PraRag2}.
Filtrations are also of great importance in the representation theory of non-Archimedean Lie groups, see for instance \cite{MoyPrasad1} and \cite{MoyPrasad2}.

\smallskip

{\it Closed fibres and local combinatorial description of the building}.---~
We finish this brief summary of Bruhat-Tits theory by mentioning quickly two further applications of integral models for facet stabilizers.

First let us pick a facet $F \subset {\mathcal{B}}({\rm G}, k)$ as above and consider the associated $k^\circ$-group scheme $\mathcal{G}_F$.
As a scheme over $k^\circ$, it has a closed fibre (so to speak obtained by reduction modulo $k^{\circ\circ}$) which we denote by $\overline{\mathcal{G}_F}$.
This is a group scheme over the residue field $\widetilde{k}$.
It turns out that the rational points $\overline{\mathcal{G}_F}(\widetilde{k})$ have a nice combinatorial structure (even though the $\widetilde{k}$-group $\overline{\mathcal{G}_F}$ needn't be reductive in general); more precisely, $\overline{\mathcal{G}_F}(\widetilde{k})$ has a Tits system structure (see  the end of \ref{sss - RS and RD}) with finite Weyl group.
One consequence of this is that $\overline{\mathcal{G}_F}(\widetilde{k})$ admits an action on a spherical building (a {\it spherical building} is merely a simplicial complex satisfying the axioms of Definition~\ref{defi - simplicial building} with the Euclidean tiling $\Sigma$ replaced by a spherical one).
The nice point is that this spherical building naturally appears in the (Euclidean) Bruhat-Tits building ${\mathcal{B}}({\rm G}, k)$.
Namely, the set of closed facets containing $F$ is a geometric realization of the spherical building of $\overline{\mathcal{G}_F}(\widetilde{k})$ \cite[Prop. 5.1.32]{BT1b}.
In particular, for a complete valued field $k$, the building  ${\mathcal{B}}({\rm G}, k)$ is locally finite if and only if the spherical building of $\overline{\mathcal{G}_F}(\widetilde{k})$ is actually finite for each facet $F$, which amounts to requiring that the residue field $\widetilde{k}$ be finite.
Note that a metric space admits a compactification if, and only if, it is locally compact.
Therefore from this combinatorial description of neighborhoods of facets, we see that {\it the Bruhat-Tits building ${\mathcal{B}}({\rm G},k)$ admits a compactification if and only if $k$ is a local field}.

\begin{Remark}
\label{rk - parahorics}
Let us assume here that $k$ is discretely valued.
This is the context where the more classical combinatorial structure of an (affine) Tits system is relevant \cite[IV.2]{Lie456}.
Let us exhibit such a structure.
First, a parahoric subgroup in ${\rm G}(k)$ can be defined to be the image of $(\mathcal{G}_F)^\circ (k^\circ)$ for some facet $F$ in ${\mathcal{B}}({\rm G},k)$, where $(\mathcal{G}_F)^\circ$ denotes the identity component of $\mathcal{G}_F$ \cite[5.2.8]{BT1b}.
We also say for short that a parahoric subgroup is the connected stabilizer of a facet in the Bruhat-Tits building ${\mathcal{B}}({\rm G}, k)$.
If ${\rm G}$ is simply connected (in the sense of algebraic groups), then the family of parahoric subgroups is the family of abstract parabolic subgroups of a Tits system with affine Weyl group \cite[Prop. 5.2.10]{BT1b}.
An Iwahori subgroup corresponds to the case when $F$ is a maximal facet.
At last, if moreover $k$ is local with residual characteristic $p$, then an Iwahori subgroup can be characterized as the normalizer of a maximal pro-$p$ subgroup and an arbitrary parahoric subgroup as a subgroup containing an Iwahori subgroup.
\end{Remark}

Finally, the above integral models provide an important tool in the realization of Bruhat-Tits buildings in analytic spaces (and subsequent compactifications). Indeed, the fundamental step (see Theorem~\ref{thm - subgroup}) for the whole thing consists in attaching injectively to {\it any} point $x \in {\mathcal{B}}({\rm G},K)$ an affinoid subgroup ${\rm G}_x$ of the analytic space ${\rm G}^{\rm an}$ attached to ${\rm G}$, and the definition of ${\rm G}_x$ makes use of the integral models attached to vertices. But one word of caution is in order here since the connexion with integral models avoids all their subtleties! For our construction, only smooth $k^\circ$-group schemes $\mathcal{G}_F$ which are \emph{reductive} are of interest; this is not the case in general, but one can easily prove the following statement: {\it given a vertex $x \in \mathcal{B}({\rm G},k)$, there exists a finite extension $k'/k$  such that the ${k'}^{\circ}$-group scheme $\mathcal{G}'_x$, attached to the point $x$ seen as a vertex of $\mathcal{B}({\rm G},k')$, is a Chevalley-Demazure group scheme over ${k'}^\circ$}. In this situation, one can define $({\rm G} \otimes_k k')_x$ as the \emph{generic fibre} of the formal completion of $\mathcal{G}'_x$ along its special fibre; this is a $k'$-affinoid subgroup of $({\rm G} \otimes_k k')^{\rm an}$ and one invokes descent theory to produce a $k$-affinoid subgroup of ${\rm G}^{\rm an}$.

\subsubsection{A characterization of apartments}
\label{sss - characterization_apartments}
For later use, we end this section on Bruhat-Tits theory by a useful characterization of apartments inside buildings.

\vskip1mm Given a torus ${\rm S}$ over $k$, we denote by ${\rm S}^1(k)$ the maximal bounded subgroup of ${\rm S}(k)$. It is the subgroup of ${\rm S}(k)$ defined by the equations $|\chi|=1$, where $\chi$ runs over the character group of ${\rm S}$.

\begin{Prop}
\label{characterization_apartments}
Let ${\rm S}$ be a maximal split torus and let $x$ be a point of $\mathcal{B}({\rm G},k)$. If the residue field of $k$ contains at least four elements, then the following conditions are equivalent:
\begin{itemize}
\item[(i)] $x$ belongs to the apartment ${\rm A}({\rm S},k)$;
\item[(ii)] $x$ is fixed under the action of ${\rm S}^1(k)$.
\end{itemize}
\end{Prop}

\begin{proof}
Condition (i) always implies condition (ii). With our hypothesis on the cardinality of the residue field, the converse implication holds by \cite[Proposition 5.1.37]{BT1b}. \end{proof}


\section{Buildings and Berkovich spaces}
\label{s - general compactifications}

As above, we consider a semisimple group ${\rm G}$ over some non-Archimedean field $k$. In this section, we explain how to realize the Bruhat-Tits building $\mathcal{B}({\rm G},k)$ of ${\rm G}(k)$ in non-Archimedean analytic spaces deduced from ${\rm G}$, and we present two procedures that can be used to compactify Bruhat-Tits buildings in full generality; as we pointed out before, the term ``compactification'' is abusive if $k$ is not a local field (see the discussion before Remark \ref{rk - parahorics}).

\smallskip
Assuming that $k$ is locally compact, let us describe very briefly those two ways of compactifying a building. The first is due to V. Berkovich when ${\rm G}$ is split \cite[Chap. V]{Ber1} and it consists in two steps:

1. to define a closed embedding of the building into the analytification of the group (\ref{ss - closed embedding});

2. to compose this closed embedding with an analytic map from the group to a (compact) flag variety (\ref{ss-maps_to_flags}).

By taking the closure of the image of the composed map, we obtain an equivariant compactification which admits a Lie-theoretic description (as expected). For instance, there is a convenient description of this ${\rm G}(k)$-topological space (convergence of sequences, boundary strata etc.) by means of invariant fans in $\left({\rm X}_*({\rm S}) \otimes_{\mathbf{Z}} \mathbf{R}, {\rm W} \right)$, where ${\rm X}_*({\rm S})$ denotes the cocharacter group of a maximal split torus ${\rm S}$ endowed with the natural action of the Weyl group ${\rm W}$ (\ref{ss - fans}).
The finite family of compactifications obtained in this way is indexed by ${\rm G}(k)$-conjugacy classes of parabolic subgroups.

\smallskip These spaces can be recovered from a different point of view, using representation theory and the concrete compactification $\overline{\mathcal{X}}({\rm V},k)$ of the building $\mathcal{X}({\rm V},k)$ of ${\rm SL}({\rm V},k)$ which was described in Section 2. It mimics the original strategy of I. Satake in the case of symmetric spaces \cite{Satake1}: we pick a faithful linear representation of ${\rm G}$ and, relying on analytic geometry, we embed $\mathcal{B}({\rm G},k)$ in $\mathcal{X}({\rm V},k)$; by taking the closure in $\overline{\mathcal{X}}({\rm V},k)$, we obtain our compactification.

\smallskip \noindent \textbf{Caution} --- 1. We need some functoriality assumption on the building with respect to the field: in a sense which was made precise after the statement of Theorem \ref{th - functoriality of BT buildings}, this means that $\mathcal{B}({\rm G}, -)$ is functor on the category of non-Archimedean extensions of $k$.

\smallskip
As explained in \cite[1.3.4]{RTW1}, these assumptions are fulfilled if
$k$ quasi-splits over a tamely ramified extension of $k$. This is in
particular the case is $k$ is discretely valued with perfect residue field, or if ${\rm G}$ is split.

2. There is no other restriction on the non-Archimedean field $k$ considered in 4.1. From 4.2 on, we assume that $k$ is \emph{local}. In any case, the reader should keep in mind that non-local extensions of $k$ do always appear in the study of Bruhat-Tits buildings from Berkovich's point of view (see Proposition 4.2).

\smallskip
The references for the results quoted in this section are \cite{RTW1} and \cite{RTW2}.

\subsection{Realizing buildings inside Berkovich spaces}
\label{ss - closed embedding}
Let $k$ be a field which is complete with respect to a non-trivial non-Archimedean absolute value.
We fix a semisimple group  $\G$ over $k$.
Our first goal is to define a continuous injection of the Bruhat-Tits building $\mathcal{B}(\G,k)$ in the Berkovich space $\G^{\rm an}$ associated to the algebraic group $\G$. Since $\G$ is affine with affine coordinate ring $\mathcal{O}(\G)$, its analytification consists of all multiplicative seminorms on $\mathcal{O}(\G)$ extending the absolute value on $k$ \cite{Temkin}.

\subsubsection{Non-Archimedean extensions and universal points}
\label{sss - n-A extensions}

\medskip
We will have to consider infinite non-Archimedean extensions of $k$ as in the following example.

\begin{Example}
\label{ex - field}
Let $\mathbf{r} = (r_1, \ldots, r_n)$ be a tuple of positive real numbers such that
$r_1^{i_1} \ldots r_n^{i_n} \notin |k^\times|$ for all choices of $(i_1, \ldots, i_n) \in \mathbf{Z}^n - \{0\}$. Then the $k$-algebra
\[k_{\mathbf{r}} = \left\{\sum_{I = (i_1 \ldots, i_n)} a_I x_1^{i_1}\ldots x_n^{i_n} \in k[[x_1^{\pm 1},\ldots, x_n^{\pm 1}]] \  ; \ |a_I| r_1^{i_1} \ldots r_n^{i_n} \rightarrow 0 \text{ when }|i_1| + \ldots + |i_n| \rightarrow \infty \right\} \]
is a  non-Archimedean field extension of $k$ with absolute value $|f| = \max_I \{|a_I| r_1^{i_1} \ldots r_n^{i_n} \}$.
\end{Example}

\vskip2mm We also need to recall the notion of a \emph{universal} point \footnote{This notion was introduced by Berkovich, who used the adjective \emph{peaked} \cite[5.2]{Ber1}. Its study was carried on by Poineau, who prefered the adjective \emph{universal} \cite{Poineau}.}. Let $z$ be a point in ${\rm G}^{\rm an}$, seen as a multiplicative $k$-seminorm on $\mathcal{O}({\rm G})$. For a given non-Archimedean field extension $K/k$, there is a natural $K$-seminorm $||.|| = z \otimes 1$ on $\mathcal{O}({\rm G}) \otimes_k K$, defined by $$||a|| = \inf \max_i |a_i(z)|\cdot |\lambda_i|$$
where the infimum is taken over the set of all expressions $\sum_i a_i \otimes \lambda_i$ representing $a$, with $a_i \in \mathcal{O}({\rm G})$ and $\lambda_i \in K$. The point $z$ is said to be \emph{universal} if, for any non-Archimedean field extension $K/k$, the above $K$-seminorm on $\mathcal{O}({\rm G}) \otimes_k K$ is multiplicative. One writes $z_K$ for the corresponding point in ${\rm G}^{\rm an} \widehat{\otimes}_k K$. We observe that this condition depends only on the completed residue field $\mathcal{H}(z)$ of ${\rm G}^{\rm an}$ at $z$.

\vskip2mm \begin{Remark} \label{Rk-universal} 1. Obviously, points of ${\rm G}^{\rm an}$ coming from $k$-rational points of ${\rm G}$ are universal.
\vskip1pt
2. Let $x \in {\rm G}^{\rm an}$ be universal. For any finite Galois extension $k'/k$, the canonical extension $x_{k'}$ of $x$ to ${\rm G}^{\rm an} \otimes_k k'$ is invariant under the action of ${\rm Gal}(k'/k)$: indeed, the $k'$-norm $x \otimes 1$ on $\mathcal{O}({\rm G}) \otimes_k k'$ is Galois invariant.
\vskip1pt
3. If $k$ is algebraically closed, Poineau proved that every point of ${\rm G}^{\rm an}$ is universal \cite[Corollaire 4.10]{Poineau}.
\end{Remark}

\subsubsection{Improving transitivity}
\label{sss-transitivity}

Now let $\G^{\rm an}$ be the Berkovich analytic space associated to the algebraic group $\G$. Our goal is the first step mentioned in the introduction, namely the definition of a continuous injection
\[\vartheta: \mathcal{B}(\G,k) \longrightarrow \G^{\rm an}.\]
We proceed as follows.
For every point $x$ in the building $\mathcal{B}(\G,k)$ we construct an affinoid subgroup $\G_x$ of $\G^{\rm an}$ such that, for any non-Archimedean extension $K/k$, the subgroup ${\rm G}_x(K)$ of ${\rm G}(K)$ is precisely the stabilizer of $x$ in the building over $K$. Then we define $\vartheta(x)$ as the (multiplicative) seminorm on $\mathcal{O}(\G)$ defined by taking the maximum over the compact subset ${\rm G}_x$ of ${\rm G}^{\rm an}$.

\medskip If the Bruhat-Tits building $\mathcal{B}({\rm G},k)$ can be seen as non-Archimedean analogue of a Riemannian symmetric space, it is not homogeneous under ${\rm G}(k)$; for example, if $k$ is discretely valued, the building carries a polysimplicial structure which is preserved by the action of ${\rm G}(k)$. There is a very simple way to remedy at this situation using field extensions, and this is where our functoriality assumption comes in.

\vskip2mm Let us first of all recall that the notion of a special point was defined in Section 1, just before Definition \ref{defi - non-simplicial building}. Its importance comes from the fact that, when ${\rm G}$ is split, the stabilizer of a special point is particularly nice (see the discussion after Theorem~\ref{thm - subgroup}). As simple consequences of the definition, one should notice the following two properties: if a point $x \in \mathcal{B}({\rm G},k)$ is special, then
\begin{itemize}
\item[-] every point in the ${\rm G}(k)$-orbit of $x$ is again special;
\item[-] if moreover ${\rm G}$ is \emph{split}, then $x$ remains special in $\mathcal{B}({\rm G},K)$ for any non-Archimedean field extension $K/k$ (indeed: the local Weyl group at $x$ over ${\rm K}$ contains the local Weyl group at $x$ over $k$, and the full Weyl group of ${\rm G}$ is the same over $k$ and over $K$).
\end{itemize}

\vskip1mm We can now explain how field extensions allow to improve transitivity of the group action on the building.

\begin{Prop}
\label{prop - special point}
1. Given any two points $x, y \in \mathcal{B}({\rm G},k)$, there exists a non-Archimedean field extension $K/k$ such that $x$ and $y$, identified with points of $\mathcal{B}({\rm G},K)$ via the canonical injection $\mathcal{B}({\rm G},k) \hookrightarrow \mathcal{B}({\rm G},K)$, belong to the same orbit under ${\rm G}(K)$.

\vskip1mm 2. For every point $x \in \mathcal{B}(\G,k)$, there exists a non-Archimedean field extension $K/k$ such that the following conditions hold:

{\rm (i)}~The group $\G \otimes_{k} K$ is split;
{\rm (ii)}~The canonical injection $\mathcal{B}(\G,k) \rightarrow \mathcal{B}(\G, K)$ maps $x$ to a special point.
\end{Prop}

\medskip
We give a proof of this Proposition since it is the a key result for the investigation of Bruhat-Tits buildings from Berkovich's point of view. The second assertion follows easily from the first: just pick a finite separable field extension $k'/k$ splitting ${\rm G}$ and a special point $y$ in $\mathcal{B}({\rm G},k')$, then consider a non-Archimedean field extension $K/k'$ such that $x$ and $y$ belong to the same ${\rm G}(K)$-orbit. In order to prove the first assertion, we may and do assume that ${\rm G}$ is split. Let ${\rm S}$ denote a maximal split torus of ${\rm G}$ whose apartment ${\rm A}({\rm S},k)$ contains both $x$ and $y$. As recalled in Proposition 3.17, this apartment is an affine space under ${\rm X}_*({\rm S}) \otimes_{\mathbf{Z}} \mathbf{R}$, where ${\rm X}_*({\rm S})$ denotes the cocharacter space of ${\rm S}$, and ${\rm S}(k)$ acts on ${\rm A}({\rm S},k)$ by translation via a map $\nu : {\rm S}(k) \rightarrow {\rm X}_*({\rm S}) \otimes_{\mathbf{Z}} \mathbf{R}$. Using a basis of characters to identify ${\rm X}_*({\rm S})$ (resp. ${\rm S}$) with $\mathbf{Z}^n$ (resp. $\mathbf{G}_{\rm m}^n$), it turns out that $\nu$ is simply the map $$k^{\times} \longrightarrow \mathbf{R}^n, \ \ \ (t_1,\ldots, t_n) \mapsto (-\log |t_1|, \ldots, -\log |t_n|).$$
By combining finite field extensions and transcendental extensions as described in Example 4.1, we can construct a non-Archimedean field extension $K/k$ such that the vector $x-y \in \mathbf{R}^n$ belongs to the subgroup $\log |(K^\times)^n|$. This implies that $x$ and $y$, seen as points of ${\rm A}({\rm S},K)$, belong to the same orbit under ${\rm S}(K)$, hence under ${\rm G}(K)$.

\begin{Remark}
If $|K^\times| = \mathbf{R}_{>0}$, then ${\rm G}(K)$ acts transitively on $\mathcal{B}({\rm G},K)$. However, it is more natural to work functorially than to fix arbitrarily an algebraically closed non-Archimedean extension $\Omega/k$ such that $|\Omega^\times| = \mathbf{R}_{>0}$.
\end{Remark}

\subsubsection{Affinoid subgroups}
\label{sss - affinoid}
Let us now describe the key fact explaining the relationship between Bruhat-Tits theory and non-Archimedean analytic geometry. This result is crucial for all subsequent constructions.

\begin{Thm}
\label{thm - subgroup}
For every point $x \in \mathcal{B}(\G,k)$ there exists a unique $k$-affinoid subgroup $\G_x$ of $\G^{\rm an}$ satisfying the following condition: for every non-Archimedean field extension $K/k$, the group $\G_x(K)$ is the stabilizer in $\G(K)$ of the image of $x$ under the injection
$\mathcal{B}(\G,k) \rightarrow \mathcal{B}(\G,K)$.
\end{Thm}

The idea of the proof is the following (see \cite[Th. 2.1]{RTW1} for details). If $\G$ is split and $x$ is a special point in the building, then the integral model $\mathcal{G}_x$ of ${\rm G}$ described in (3.2.3) is a Chevalley group scheme, and we define ${\rm G}_x$ as the generic fibre $\widehat{\mathcal{G}_x}_\eta$ of the formal completion of $\mathcal{G}_x$ along its special fibre. This is a $k$-affinoid subgroup of ${\rm G}^{\rm an}$, and it is easy to check that it satisfies the universal property in our claim. Thanks to Proposition \ref{prop - special point}, we can achieve this situation after a suitable non-Archimedean extension $K/k$, and we apply faithfully flat descent to obtain the $k$-affinoid subgroup ${\rm G}_x$ \cite[App. A]{RTW1}. Let us remark that, in order to perform this descent step, it is necessary to work with an extension which is not too big (technically, the field $K$ should be a $k$-affinoid algebra); since one can obtain $K$ by combining finite extensions with the transcendental ones described in Example \ref{ex - field}, this is fine.

\subsubsection{Closed embedding in the analytic group}
\label{sss - closed embedding}
The $k$-affinoid subgroup $\G_x$ is the Berkovich spectrum of a $k$-affinoid algebra ${\rm A}_x$, i.e., $\G_x$ is the Gelfand spectrum $\mathcal{M}({\rm A}_x)$ of bounded multiplicative seminorms on ${\rm A}_x$. This is a compact and Hausdorff topological space over which elements of ${\rm A}_x$ define non-negative real valued functions. For any non-zero $k$-affinoid algebra ${\rm A}$, one can show that its Gelfand spectrum $\mathcal{M}({\rm A})$ contains a smallest non-empty subset, called its \emph{Shilov boundary} and denoted $\Gamma({\rm A})$, such that each element $f$ of ${\rm A}$ reaches its maximum at some point in $\Gamma({\rm A})$.

\medskip \begin{Remark} (i) If ${\rm A} = k\{{\rm T}\}$ is the Tate algebra of restricted power series in one variable, then $\mathcal{M}({\rm A})$ is Berkovich's closed unit disc and its Shilov boundary is reduced to the point $o$ defined by the Gauss norm: for $f = \sum_{n \in \mathbf{N}} a_n{\rm T}^n$, one has $|f(o)| = \max_n |a_n|$.

\vskip1pt
(ii) Let $a \in k$ with $0 < |a| < 1$. If ${\rm A} = k\{{\rm T}, {\rm S}\}/({\rm ST}-a)$, then $\mathcal{M}({\rm A})$ is an annulus of modulus $|a|$ and $\Gamma({\rm A})$ contains two points $o, o'$: for $f = \sum_{n \in \mathbf{Z}} a_n {\rm T}^n$, where ${\rm T}^{-1} = a^{-1}{\rm S}$, one has $|f(o)| = \max_n |a_n|$ and $|f(o')| = \max_n |a_n|.|a|^n$.

\vskip1pt
(iii) For any non-zero $k$-affinoid algebra ${\rm A}$, its Shilov boundary $\Gamma({\rm A})$ is reduced to a point if and only if the seminorm $${\rm A} \rightarrow \mathbf{R}_{\geqslant 0}, \ \ \ f \mapsto \sup_{x \in \mathcal{M}({\rm A})} |f(x)|$$ is multiplicative.
\end{Remark}

For every point $x$ of $\mathcal{B}({\rm G}, k)$, it turns out that the Shilov boundary of ${\rm G}_x = \mathcal{M}({\rm A}_x)$ is reduced to a unique point, denoted $\vartheta(x)$. This is easily seen by combining the nice behavior of Shilov boundaries under non-Archimedean extensions, together with a natural bijection between the Shilov boundary of $\mathcal{V}_\eta$ and the set of irreducible components of $\mathcal{V} \otimes_{k^\circ} \widetilde{k}$ if $\mathcal{V}$ is a normal $k^\circ$-formal scheme; indeed, the smooth $k^\circ$-group scheme $\mathcal{G}_x$ has a connected special fibre when it is a Chevalley group scheme. Let us also note that the affinoid subgroup ${\rm G}_x$ is completely determined by the single point $\vartheta(x)$ via

$${\rm G}_x = \{z \in {\rm G}^{\rm an} \ ; \ \forall f \in \mathcal{O}({\rm G}), \ |f(z)| \leqslant |f(\vartheta(x))|\}.$$

\medskip
In this way we define the desired map
\[\vartheta: \mathcal{B}(\G, k) \rightarrow \G^{\rm an},\]

\noindent and we show \cite[Prop. 2.7]{RTW1} that it is injective, continuous and $\G(k)$-equivariant (where $\G(k)$ acts on $\G^{\rm an}$ by conjugation). If $k$ is a local field,  $\vartheta$ induces a homeomorphism from $\mathcal{B}(\G,k)$ to a closed subspace of $\G^{\rm an}$ \cite[Prop. 2.11]{RTW1}.

\medskip
Finally, the map $\vartheta $ is also compatible with non-Archimedean extensions $K/k$, i.e., the following diagram

$$\xymatrix{\mathcal{B}({\rm G},K) \ar@{->}[r]^{\vartheta_K} & ({\rm G} \otimes_k {\rm K})^{\rm an} \ar@{->}[d]^{p_{K/k}} \\ \mathcal{B}({\rm G},k) \ar@{->}[r]_\vartheta \ar@{->}[u]^{\iota_{K/k}} & {\rm G}^{\rm an}}$$
where $\iota_{K/k}$ (resp. $p_{K/k}$) is the canonical embedding (resp. projection) is commutative. In particular, we see that this defines a \emph{section} of $p_{K/k}$ over the image of $\vartheta$. In fact, any point $z$ belonging to this subset of ${\rm G}^{\rm an}$ is universal (\ref{sss - n-A extensions}) and $\vartheta_K(\iota_{K/k}(x))$ coincides with the canonical lift $\vartheta(x)_K$ of $\vartheta(x)$ to $({\rm G} \otimes_k K)^{\rm an}$ for any $x \in \mathcal{B}({\rm G},k)$.

Moreover, if $K/k$ is a Galois extension, then the upper arrow in the diagram is ${\rm Gal}(K/k)$-equivariant by \cite[Prop. 2.7]{RTW1}.


\subsection{Compactifying buildings with analytic flag varieties}
\label{ss-maps_to_flags}

Once the building has been realized in the analytic space ${\rm G}^{\rm an}$, it is easy to obtain compactifications. In order not to misuse the latter word, we assume from now one that $k$ is \emph{locally compact}.

\subsubsection{Maps to flag varieties}
\label{sss - maps to flags}
The embedding $\vartheta: \mathcal{B}(\G,k) \rightarrow \G^{\rm an}$ defined in \ref{sss - closed embedding} can be used to compactify the Bruhat-Tits building $\mathcal{B}(\G,k)$.
We choose a parabolic subgroup $\rP$ of $\G$.
Then the  flag variety $\G/\rP$ is complete, and therefore the associated Berkovich space $(\G/\rP)^{\rm an}$ is compact.
Hence we can map the building to a compact space by the composition
\[\vartheta_P: \mathcal{B}(\G,k) \stackrel{\vartheta}{\longrightarrow} \G^{\rm an} \longrightarrow (\G / \rP)^{\rm an}.\]
The map $\vartheta_P$ is by construction $\G(k)$-equivariant and it depends only on the ${\rm G}(k)$-conjugacy class of ${\rm P}$: we have $\vartheta_{g{\rm P}g^{-1}} = g\vartheta_{\rm P}g^{-1}$ for any $g \in {\rm G}(k)$.

\smallskip However, $\vartheta_P$ may not be injective. By the structure theory of semisimple groups, there exists a finite family of normal reductive subgroups $\G_i$ of $\G$ (each of them quasi-simple), such that the product morphism
\[\prod_i \G_i \longrightarrow \G\]
is a central isogeny. Then the building $\mathcal{B}(\G,k)$ can be identified with the product of all $\mathcal{B}(\G_i,k)$. If one of the factors $\G_i$ is contained in $\rP$, then the factor $\mathcal{B}(\G_i,k)$ is squashed down to a point in the analytic flag variety $(\G / \rP)^{\rm an}$.

If we remove from $\mathcal{B}({\rm G},k)$ all factors $\mathcal{B}(\G_i, k)$ such that $\G_i$ is contained in ${\rm P}$, then we obtain a building $\mathcal{B}_t(\G,k)$, where $t$ stands for the type of the parabolic subgroup $\rP$, i.e., for its ${\rm G}(k)$-conjugacy class. The factor $\mathcal{B}_t(\G,k)$ is mapped injectively into $(\G/\rP)^{\rm an}$ via $\vartheta_\rP$.

\begin{Remark}
\label{rk - injectivity}
If $\G$ is almost simple, then $\vartheta_\rP$ is injective whenever $\rP$ is a proper parabolic subgroup in $\G$; hence in this case the map $\vartheta_\rP$ provides an embedding of $\mathcal{B}(\G,k)$ into $(\G/ \rP)^{\rm an}$.
\end{Remark}

\subsubsection{Berkovich compactifications}
\label{sss - Berkovich compactifications}
Allowing compactifications of the building in which some factors are squashed down to a point, we introduce the following definition.

\begin{Def}
\label{defi - Berkovich compactification}
Let $t$ be a ${\rm G}(k)$-conjugacy class of parabolic subgroups of ${\rm G}$. We define $\overline{\mathcal{B}}_t(\G,k)$ to be the closure of the image of $\mathcal{B}(\G,k)$ in $(\G/\rP)^{\rm an}$ under $\vartheta_P$, where ${\rm P}$ belongs to $t$, and we endow this space with the induced topology.
The compact space $\overline{\mathcal{B}}_t(\G,k)$ is called the {\rm Berkovich compactification of type $t$}~of the building $\mathcal{B}(\G,k)$.
\end{Def}

Note that we obtain one compactification for each ${\rm G}(k)$-conjugacy class of parabolic subgroups.

\begin{Remark}
\label{rk - non locally compact}
If we drop the assumption that $k$ is  locally compact, the map $\vartheta_\rP$ is continuous but the image of $\mathcal{B}_t({\rm G},k)$ is not locally closed. In this case, the right way to proceed is to compactify each apartment ${\rm A}_t({\rm S},k)$ of $\mathcal{B}_t({\rm G},k)$ by closing it in ${\rm G}^{\rm an}/{\rm P}^{\rm an}$ and to define $\overline{\mathcal{B}}_t({\rm G},k)$ as the union of all compactified apartments. This set is a quotient of ${\rm G}(k) \times \overline{A}_t({\rm S},k)$ and we endow it with the quotient topology \cite[3.4]{RTW1}.
\end{Remark}

\subsubsection{The boundary}
\label{sss - boundary}
Now we want to describe the boundary of the Berkovich compactifications. We fix a type $t$ (i.e., a ${\rm G}(k)$-conjugacy class) of parabolic subgroups.

\begin{Def}
\label{defi - osculatory}
Two parabolic subgroups $\rP$ and ${\rm Q}$ of $\G$ are called {\rm osculatory}~if their intersection $\rP \cap {\rm Q}$ is also a parabolic subgroup.
\end{Def}

Hence $\rP$ and ${\rm Q}$ are osculatory if and only if they contain a common Borel group after a suitable field extension. We can generalize this definition to semisimple groups over arbitrary base schemes.
Then for every parabolic subgroup ${\rm Q}$  there is a variety $\mathrm{Osc}_t({\rm Q})$ over $k$ representing the functor which associates to any base scheme ${\rm S}$ the set of all parabolics of type $t$ over ${\rm S}$ which are osculatory to ${\rm Q}$ \cite[Prop. 3.2]{RTW1}.

\begin{Def}
\label{defi - t-relevant}
Let ${\rm Q}$ be a parabolic subgroup.
We say that ${\rm Q}$ is {\rm $t$-relevant}~if there is no parabolic subgroup ${\rm Q}'$ strictly containing ${\rm Q}$ such that $\mathrm{Osc}_t({\rm Q}) = \mathrm{Osc}_t({\rm Q}')$.
\end{Def}
\vspace{0,3cm}

Let us illustrate this definition with the following example.

\begin{Example}
\label{example - relevant}
Let $\G$ be the group ${\rm SL}({\rm V})$, where ${\rm V}$ is a  $k$-vector space of dimension $d +1$.
The non-trivial parabolic subgroups of $\G$ are the stabilizers of flags
\[(0 \subsetneq {\rm V}_1 \subsetneq \ldots \subsetneq {\rm V}_r \subsetneq {\rm V}).\]
Let ${\rm H}$ be a hyperplane in ${\rm V}$, and let $\rP$ be the parabolic subgroup of ${\rm SL}({\rm V})$ stabilizing the flag $(0 \subset {\rm H} \subset {\rm V})$. We denote its type by $\delta$.
Let ${\rm Q}$ be an arbitrary parabolic subgroup, stabilizing a flag  $(0 \subsetneq {\rm V}_1 \subsetneq \ldots \subsetneq {\rm V}_r \subsetneq {\rm V})$.
Then ${\rm Q}$ and ${\rm P}$ are  osculatory if and only if ${\rm H}$ contains the linear subspace ${\rm V}_r$.
This shows that all parabolic subgroups ${\rm Q}$ stabilizing flags contained in the subspace ${\rm V}_r$ give rise to the same variety $\mathrm{Osc}_\delta({\rm Q})$.
Therefore, a non-trivial parabolic is $\delta$-relevant if and only if the corresponding flag has the form $0 \subsetneq {\rm W} \subsetneq {\rm V}$.
\end{Example}

Having understood how to parametrize boundary strata, we can now give the general description of the Berkovich compactification $\overline{\mathcal{B}}_t(\G,k)$. The following result is \cite[Theorem 4.1]{RTW1}.

\begin{Thm}
\label{thm - stratification}
For every $t$-relevant parabolic subgroup ${\rm Q}$, let ${\rm Q}_{\rm ss}$ be its semisimplification (i.e., ${\rm Q}_{\rm ss}$ is the quotient ${\rm Q}/\mathcal{R}({\rm Q})$ where $\mathcal{R}({\rm Q})$ denotes the radical of ${\rm Q}$).
Then $\overline{\mathcal{B}}_t(\G,k)$ is the disjoint union of all the buildings $\mathcal{B}_t({\rm Q}_{ss},k)$, where ${\rm Q}$ runs over the $t$-relevant parabolic subgroups of $\G$.
\end{Thm}

The fact that the Berkovich compactifications of a given group are contained in the flag varieties of this group enables one to have natural maps between compactifications: they are the restrictions to the compactifications of (the analytic maps associated to) the natural fibrations between the flag varieties.
The above combinatorics of $t$-relevancy is a useful tool to formulate which boundary components are shrunk when passing from a compactification to a smaller one \cite[Section 4.2]{RTW1}.

\begin{Example}
\label{example - stratification}
Let us continue the discussion in Example \ref{example - relevant} by describing the stratification of $\overline{\mathcal{B}}_\delta({\rm SL}({\rm V}),k)$.
Any $\delta$-relevant subgroup ${\rm Q}$ of $\G = {\rm SL}({\rm V})$ is either equal to ${\rm SL}({\rm V})$ or equal to the stabilizer of a linear subspace $0 \subsetneq {\rm W} \subsetneq {\rm V}$.
In the latter case ${\rm Q}_{\rm ss}$ is isogeneous to ${\rm SL}({\rm W}) \times {\rm SL}({\rm V}/{\rm W})$.
Now ${\rm SL}({\rm W})$ is contained in a parabolic of type $\delta$, hence $\mathcal{B}_\delta({\rm Q}_{\rm ss},k)$ coincides with $\mathcal{B}({\rm SL}({\rm V}/{\rm W}),k)$. Therefore
\[\overline{\mathcal{B}}_\delta({\rm SL}({\rm V}),k) = \bigcup_{{\rm W} \subsetneq {\rm V}} \mathcal{B}\bigl({\rm SL}({\rm V}/{\rm W},k)\bigr),\]
where ${\rm W}$ runs over all linear subspaces ${\rm W} \subsetneq {\rm V}$.
\end{Example}

Recall from \ref{SL-building} that the Euclidean building $\mathcal{B}({\rm SL}({\rm V}),k)$ can be identified with the Goldman-Iwahori space $\mathcal{X}({\rm V},k)$ defined in \ref{defi - GI}.
Hence $\overline{\mathcal{B}}_\delta({\rm SL}({\rm V}),k)$ is the disjoint union of all $\mathcal{X}({\rm V}/{\rm W},k)$.
Therefore we can identify the seminorm compactification $\overline{\mathcal{X}}(V,k)$ from \ref{ss - seminorm compactification} with the Berkovich compactification of type $\delta$.

\subsection{Invariant fans and other compactifications}
\label{ss - fans}

Our next goal is to compare our approach to compactifying building with another
one, developed in \cite{Wer07} without making use of Berkovich
geometry. In this work, compactified buildings are defined by a gluing procedure, similar to the one defining the Bruhat-Tits building in Theorem \ref{th - gluing}.  In a first step, compactifications of apartments are obtained by a cone decomposition. Then these compactified apartments are glued together with the help of subgroups which turn out to be the stabilizers of points in the compactified building.

Let $\G$ be a (connected) semisimple group over $k$ and $\mathcal{B}(\G,k)$ the associated Bruhat-Tits building. We fix a maximal split torus $\T$ in $\G$, giving rise to the cocharacter space $\Sigma_{\rm vect} = {\rm X}_\ast(\T) \otimes \mathbf{R}$. The starting point is a faithful, geometrically irreducible representation $\rho: \G \rightarrow {\rm GL}(\V)$ on some finite-dimensional $k$-vector space $\V$.

Let $R=R(\T,\G) \subset {\rm X}^\ast(\T)$ be the associated root system. We fix a basis $\Delta$ of $R$ and denote by $\lambda_0(\Delta)$ the highest weight of the representation $\rho$ with respect to $\Delta$. Then every other ($k$-rational) weight of $\rho$ is of the form $\lambda_0(\Delta) - \sum_{\alpha \in \Delta} n_\alpha \alpha$ with coefficients $n_\alpha \geqslant 0$. We write $[\lambda_0(\Delta) - \lambda] = \{ \alpha \in \Delta:  n_\alpha >0\}$. We call every such subset $Y$ of $\Delta$ of the form $Y = [\lambda_0(\Delta) - \lambda]$ for some weight $\lambda$ \emph{admissible}.

\begin{Def}
Let ${\rm Y} \subset \Delta$ be an admissible subset. We denote by ${\rm C}_{\rm Y}^\Delta$ the following cone in $\Sigma_{\rm vect}$:
\[{\rm C}_{\rm Y}^\Delta = \left\{ x \in \Sigma_{\rm vect} \ ; \  \begin{array}{ll} \alpha(x) = 0  & \text{ for all }\alpha \in {\rm Y}, \text{ and } \\
 (\lambda_0(\Delta) - \lambda)(x) \geqslant 0  & \text{ for all weights }\lambda \text{ such that } [\lambda_0(\Delta) - \lambda] \not\subset {\rm Y}\end{array} \right\} \]
\end{Def}

The collection of all cones ${\rm C}_{\rm Y}^\Delta$, where $\Delta$ runs over all bases of the root system and ${\rm Y}$ over all admissible subsets of $\Delta$, is a complete fan $\mathcal{F}_\rho$ in $\Sigma_{\rm vect}$. There is a natural compactification of $\Sigma_{\rm vect}$ associated to $\mathcal{F}_\rho$, which is defined as $\overline{\Sigma}_{\rm vect} = \bigcup_{{\rm C} \in \mathcal{F}_\rho} \Sigma_{\rm vect} / \langle {\rm C} \rangle$ endowed with a topology given by tubular neighborhoods around boundary points. For details see \cite[Section 2]{Wer07} or \cite[Appendix B]{RTW1}.

We will describe this compactification in two examples.

\begin{Example}
\label{regularweight} If the highest weight of $\rho$ is regular, then every subset ${\rm Y}$ of $\Delta$ is admissible. In this case, the fan $\mathcal{F}_\rho$ is the full Weyl fan. In the case of a root system of type ${\rm A}_2$, the resulting compactification is shown on Figure 1. The shaded area is a compactified Weyl chamber, whose interior contains the corresponding highest weight of $\rho$.
\end{Example}

\begin{Example}
\label{identityrep}
Let $\G = {\rm SL}(\V)$ be the special linear group of a $(d+1)$-dimensional $k$-vector space ${\rm V}$, and let $\rho$ be the identical representation. We look at the torus $\T$ of diagonal matrices in ${\rm SL}(\V)$, which gives rise to the root system $R = \{\alpha_{i,j}\}$ of type ${\rm A}_{d}$ described in Example \ref{ex - RS of type A}. Then $\Delta = \{\alpha_{0,1}, \alpha_{1,2}, \ldots, \alpha_{d-1, d}\}$ is a basis of $R$ and
$\lambda_0(\Delta) = \varepsilon_0$ in the notation of Example \ref{ex - RS of type A}. The other weights of the identical representation are $\varepsilon_1, \ldots, \varepsilon_d$. Hence the admissible subsets of $\Delta$ are precisely the sets ${\rm Y}_r  = \{\alpha_{0,1}, \ldots, \alpha_{r-1 ,r}\}$ for $r = 1, \ldots, d$, and
${\rm Y}_{0} = \varnothing$. Let $\eta_0, \ldots, \eta_d$ be the dual basis of $\varepsilon_0, \ldots, \varepsilon_d$. Then $\Sigma_{\rm vect}$ can be identified with $\bigoplus_{i = 0}^d \mathbf{R} \eta_i / \mathbf{R} (\sum_{i} \eta_i)$, and we find
\[{\rm C}_{{\rm Y}_r}^\Delta = \{ x = \sum_i x_i \eta_i \in \Sigma_{\rm vect}: x_0 = \ldots = x_r \text{ and }x_0 \geqslant x_{r+1},   x_0 \geqslant x_{r+2},  \ldots, x_0 \geqslant  x_d\} / \mathbf{R} (\sum_i \eta_i) \]
The associated compactification is shown in Figure 2. The shaded area is a compactified Weyl chamber and its codimension one face marked by an arrow contains the highest weight of $\rho$ (with respect to this Weyl chamber).
\end{Example}

\begin{figure}[h]
\centering
\input{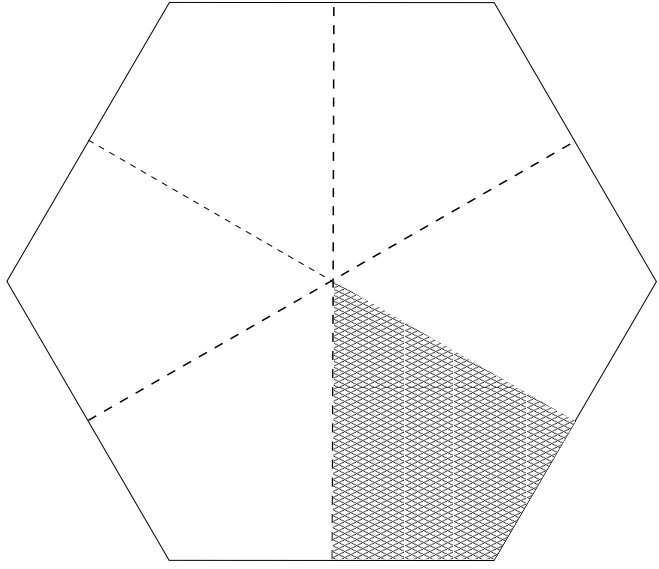}
\caption{Compactification of an apartment: regular highest weight}
\end{figure}

\begin{figure}[h]
\centering
\input{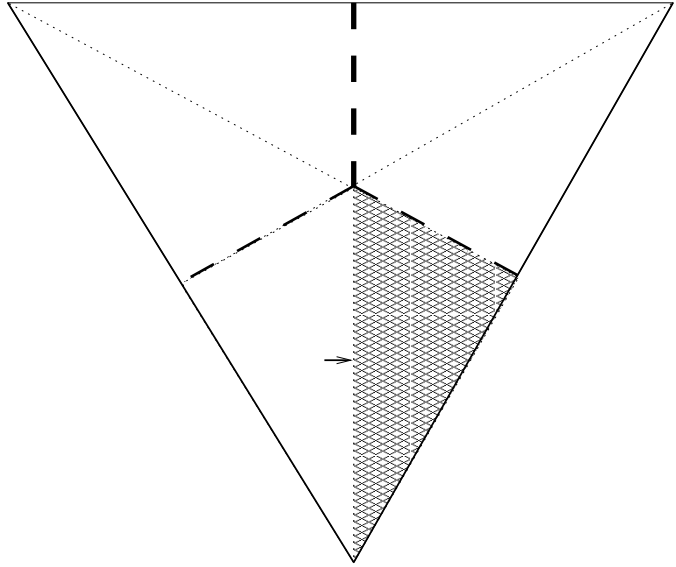}
\caption{Compactification of an apartment: singular highest weight}
\end{figure}

The compactification $\overline{\Sigma}_{\rm vect}$ induces a compactification $\overline{\Sigma}$ of the apartment $\Sigma = {\rm A}({\rm T},k)$, which is an affine space under $\Sigma_{\rm vect}$. Note that the fan $\mathcal{F}_\rho$ and hence the compactification $\overline{\Sigma}$ only depend on the Weyl chamber face containing the highest weight of $\rho$, see \cite[Theorem 4.5]{Wer07}.

Using a generalization of Bruhat-Tits theory one can define a subgroup ${\rm P}_x$ for all $x \in \overline{\Sigma}$ such that for $x \in \Sigma$ we retrieve the groups ${\rm P}_x$ defined in section 3.2, see \cite[section 3]{Wer07}. Note that by continuity the action of ${\rm N}_{\rm G}({\rm T},k)$ on $\Sigma$ extends to an action on $\overline{\Sigma}$.

\begin{Def}
\label{def - compact}
The compactification $\overline{\mathcal{B}}({\rm G},k)_\rho$
associated to the representation $\rho$ is defined as the quotient of the topological space ${\rm G}(k) \times \overline{\Sigma}$ by a similar equivalence relation as in Theorem \ref{th - gluing}:
\smallskip

\centerline{$(g,x) \sim (h,y)$ \quad $\Longleftrightarrow$ \quad there exists $n \in {\rm N}_{\rm G}({\rm T},k)$ such that $y=\nu(n).x$ and $g^{-1}hn \in {\rm P}_x$.}

\end{Def}

The compactification
of $\mathcal{B}({\rm G},k)$ with respect to a representation with regular highest weight coincides with the polyhedral compactification defined by Erasmus Landvogt in \cite{La}.

The connection to the compactifications defined with Berkovich spaces in section \ref{ss-maps_to_flags} is given by the following result, which is proved in \cite[Theorem 2.1]{RTW2}.

\begin{Thm} \label{comparison}
Let $\rho$ be a faithful, absolutely irreducible  representation of
${\rm G}$ with highest weight $\lambda_0(\Delta)$. Define
\[{\rm Z} = \{\alpha \in \Delta : \langle \alpha, \lambda_0(\Delta) \rangle = 0\},\]
where $\langle\, , \, \rangle$ is a scalar product associated to the root system as in Definition \ref{def - root system}.
We denote by $\tau$ the type of the standard parabolic subgroup of
${\rm G}$ associated to ${\rm Z}$.
Then there is a ${\rm G}(k)$-equivariant homeomorphism
\[\overline{\mathcal{B}}({\rm G},k)_\rho \rightarrow
\overline{\mathcal{B}}_\tau({\rm G},k)\]
restricting to the identity map on the building.
\end{Thm}

\begin{Example} In the situation of Example \ref{identityrep} we have
  $\lambda_0(\Delta) = \varepsilon_0$ and ${\rm Z} = \{\alpha_{1,2},
  \ldots, \alpha_{d-1, d}\}$. The associated standard parabolic is the
  stabilizer of a line. We denote its type by $\pi$. Hence the
  compactification of the building associated to ${\rm SL}(\V)$ given
  by the identity representation is the one associated to type $\pi$
  by  Theorem \ref{comparison}. This compactification was studied in
  \cite{Wer01}. It is  isomorphic to the seminorm compactification $\overline{\mathcal{X}}({\rm V}^\vee,k)$
  of the building $\mathcal{X}({\rm V}^\vee,k)$.
\end{Example}

\subsection{Satake's viewpoint}
\label{ss - Satake}

If ${\rm G}$ is a non-compact real Lie group with maximal compact subgroup ${\mr
  K}$, Satake constructed in \cite{Satake2} a compactification of the Riemannian symmetric space $\Sr = \G/{\rm K}$ in the following way:
\begin{itemize}
\item{(i)} First consider the symmetric space $\Hr$ associated to the group ${\rm PSL}(n,\mathbf{C})$ which can be identified with the space of all positive definite hermitian $n \times n$-matrices with determinant $1$. Then $\Hr$ has a natural compactification $\overline{\Hr}$ defined as the set of the homothety classes of all hermitian $n \times n$-matrices.
\item{(ii)} For an arbitrary symmetric space $\Sr = \G/{\rm K}$ use a faithful representation of $\G$ to embed $\Sr$ into $\Hr$ and consider the closure of $\Sr$ in $\overline{\Hr}$.
\end{itemize}

In the setting of Bruhat-Tits buildings we can imitate this strategy in two different ways.

{\it Functoriality of buildings} ---- The first strategy is a generalization of functoriality results for buildings developed by Landvogt \cite{LandvogtCrelle}. Let $\rho: \G \rightarrow {\rm SL}(\V)$ be a representation of the semisimple group $\G$.
Let $\Sr$ be a maximal split torus in $\G$ with normalizer $\N$, and
let ${\rm A}(\Sr,k)$ denote the corresponding apartment in
$\mathcal{B}(\G,k)$.  Choose a special vertex $o$ in ${\rm
  A}(\Sr,k)$. By \cite{LandvogtCrelle}, there exists a maximal split
torus $\T$ in ${\rm SL}({\V})$ containing $\rho(\Sr)$, and there
exists a point $o'$ in the apartment ${\rm A}(\T,k)$ of $\T$ in $\mathcal{B}({\rm SL}(V),k)$  such that the following properties hold:

\begin{enumerate}
\item There is a unique affine map between apartments $i: {\rm
  A}(\Sr,k) \rightarrow {\rm A}(\T,k)$ such that $i(o) = o'$. Its
  linear part is the map on cocharacter spaces ${\rm X}_*(\Sr)
  \otimes_\mathbf{Z} \mathbf{R} \rightarrow {\rm X}_*(\T)
  \otimes_{\mathbf{Z}} \mathbf{Z}$ induced by $\rho: {\rm S} \rightarrow \T$.
\item The map $i$ is such that $\rho(\rP_x) \subset \rP'_{i(x)}$ for
  all $x \in {\rm A}(\Sr,k)$, where $\rP_x$ denotes the stabilizer of the point $x$ with respect to the $\G(k)$-action on $\mathcal{B}(\G,k)$, and $\rP'_{i(x)}$ denotes the stabilizer of the point $i(x)$ with respect to the ${\rm SL}(\V,k)$-action on $\mathcal{B}({\rm SL}(\V),k)$.
\item The map $\rho_\ast: {\rm A}(\Sr,k) \rightarrow {\rm A}(\T,k)
  \rightarrow \mathcal{B}({\rm SL}(\V),k)$ defined by composing $i$
  with the natural embedding of the apartment ${\rm A}(\T,k)$ in the
  building $\mathcal{B}({\rm SL}({\V}),k)$ is $\N(k)$-equivariant,
  i.e., for all $x \in {\rm A}(\Sr,k)$ and $n \in \N(k)$ we have $\rho_\ast (nx) = \rho(n) \rho_\ast(x)$.
\end{enumerate}

These properties imply that $\rho_\ast: {\rm A}(\Sr,k) \rightarrow \mathcal{B}({\rm SL}(\V), k)$ can be continued to a map
$\rho_\ast: \mathcal{B}(\G,k) \rightarrow \mathcal{B}({\rm SL}({\V}),k)$, which is continuous and $\G(k)$-equivariant. By \cite[2.2.9]{LandvogtCrelle}, $\rho_\ast$ is injective.

Let $\mathcal{F}$ be the fan in ${\rm X}_*({\rm T}) \otimes_{\mathbf
  Z} \mathbf{R}$ associated to
the identity representation, which is described in Example
\ref{identityrep}. It turns out that the preimage of $\mathcal{F}$
under the map $\Sigma_{\rm vect} (\Sr,k) \rightarrow \Sigma_{\rm
  vect}(\T,k)$ induced by $\rho: \Sr \rightarrow \T$ is the fan
$\mathcal{F}_\rho$, see \cite[Lemma 5.1]{RTW2}. This implies that the
map $i$ can be extended to a map of compactified apartments
$\overline{\rm A}({\rm S},k) \rightarrow \overline{\rm A} ({\rm
  T},k)$.  An analysis of the stabilizers of boundary points
shows moreover that $\rho({\rm P}_x) \subset {\rm P}'_{i(x)}$ for all
$x \in \overline{\rm A}({\rm S},k)$, where ${\rm P}_x$ denotes the stabilizer
of $x$ in $\G(k)$, and ${\rm P}'_{i(x)}$ denotes the stabilizer of $i(x)$ in ${\rm SL}(\V,k)$ \cite[Lemma 5.2]{RTW2}.
Then it follows from the definition of $\overline{\mathcal{B}}(\G,k)_\rho$ in \ref{def - compact} that the embedding of buildings $\rho_\ast$ may be extended to a map
\[\overline{\mathcal{B}}({\rm G},k)_\rho \longrightarrow \overline{\mathcal{B}}({\rm SL}(\V),k)_{\rm id}.\]
It is shown in \cite[Theorem 5.3]{RTW2} that this map is a $\G(k)$-equivariant homeomorphism of $\overline{\mathcal{B}}(\G,k)_\rho$ onto the closure of the image of $\mathcal{B}(\G,k)$ in the right hand side.

{\it Complete flag variety} ---- Satake's strategy of embedding the
building in a fixed compactification of the building associated to
${\rm SL}(\V,k)$ can also be applied in the setting of Berkovich
spaces. Recall from \ref{SL-building} that the building $\mathcal{B}({\rm SL}(\V),k)$ can be
identified with the space $\mathcal{X}(\V,k)$ of (homothety classes
of) non-Archimedean norms on $\V$. In section \ref{ss - seminorm compactification}, we constructed a compactification
$\overline{\mathcal{X}}(\V,k)$ as the space of (homothety classes of) non-zero non-Archimedean seminorms on $\V$ and a retraction map $\tau: \mathbf{P}(\V)^{{\rm an}} \longrightarrow \overline{\mathcal{X}}(\V,k)$.

Now let $\G$ be a (connected) semisimple $k$-group together with an absolutely irreducible
projective representation $\rho: \G \rightarrow {\rm PGL}(\V,k)$. Let
$\Bor(\G)$ be the variety of all Borel groups of $\G$. We assume for
simplicity that $\G$ is quasi-split, i.e., that there exists a Borel
group $\rm B$ defined over $k$; this amounts to saying that ${\rm
  Bor}({\rm G})(k)$ is non-empty. Then $\Bor(\G)$ is isomorphic to $\G/{\rm B}$.  There is a natural morphism
\[\Bor(\G) \longrightarrow \mathbf{P}(\V)\]
such that any Borel subgroup ${\rm B}$ in $\G \otimes K$ for some field
extension $K$ of $k$ is mapped to the unique $K$-point in
$\mathbf{P}(\V)$ invariant under ${\rm B} \otimes_k K$, see \cite[Proposition 4.1]{RTW2}. Recall that in section \ref{sss - maps to flags}
we defined a map
\[\vartheta_\varnothing: \mathcal{B}(\G,k) \rightarrow {\rm Bor(G)}^{{\rm an}}\]
($\varnothing$ denotes the type of Borel subgroups). Now we consider the composition
\[ \mathcal{B}(\G,k) \stackrel{\vartheta_\varnothing}{\longrightarrow} \Bor(\G)^{{\rm an}} \rightarrow \mathbf{P}(\V)^{{\rm an}} \stackrel{\tau}{\longrightarrow} \overline{\mathcal{X}}(\V,k).\]
We can compactify the building $\mathcal{B}(\G,k)$ by taking the closure of the image.
If $\rho^\vee$ denotes  the contragredient representation of $\rho$, then it is shown in \cite[4.8 and 5.3]{RTW2} that in this way we obtain the compactification $\overline{\mathcal{B}}(\G,k)_{\rho^\vee}$.


\section{Erratum to \cite{RTW1} and \cite{RTW2}}
\label{s - Erratum}

\vskip2mm \noindent {\bf 5.1.} Tobias Schmidt pointed out that Lemma A.10 in Appendix A to \cite{RTW1} needed to be corrected. The problem comes from the fact that, for a finite Galois extension $\ell/k$ of non-Archimedean fields, the canonical map $$\lambda : \ell \otimes_k \ell \longrightarrow \prod_{{\rm Gal}(\ell|k)} \ell, \ \ \ \ a \otimes b \longmapsto \left(g(a)b\right)_{g \in {\rm Gal}(\ell|k)}$$ is not always an isometry when the left-hand side is equiped with the tensor product norm; this is the case if and only the extension is tamely ramified. 

\vskip1mm A first observation is that the algebraic isomorphism $\lambda$ is an isometry with respect to spectral norms on both sides since we are working with finite dimensional $k$-algebras. Therefore, the question amounts to understanding when the tensor product norm $|.|_\otimes$ on ${\rm A} = \ell \otimes_k \ell$ coincides with the spectral norm, which is the case if and only if $|.|_\otimes$ is power-multiplicative. Let us consider M.~Temkin's \emph{graded reduction} of $({\rm A}, |.|_\otimes)$ \cite{Temkin0}, which is to say the graded ring 
$$\tilde{\rm A}_\bullet = \bigoplus_{r \in {\bf R}_{>0}} {\rm A}_{\leqslant r}/{\rm A}_{<r}$$
where ${\rm A}_{\leqslant r} = \{a \in {\rm A} \ ; \ |a|_\otimes \leqslant r\}$ and ${\rm A}_{<r} = \{a \in {\rm A} \ ; \ |a|_\otimes < r\}$. The norm $|.|_\otimes$ is power-multiplicative, hence coincides with the spectral norm, if and only if $\tilde{\rm A}_\bullet$ is reduced. This graded ring is isomorphic to $\tilde{\ell}_\bullet \otimes_{\tilde{k}_\bullet} \tilde{\ell}_\bullet$ \cite[proof of Lemma 2.12]{Schmidt} and therefore is reduced if and only if the extension of graded fields $\tilde{\ell}_\bullet/\tilde{k}_\bullet$ is separable (since $\ell/k$ is Galois, separability of $\tilde{\ell}/\tilde{k}$ can be checked over $\tilde{\ell}$). This is the case if and only if the field extension $\ell/k$ is tamely ramified \cite[Proposition 2.10]{Ducros}.

\vskip2mm \noindent  {\bf 5.2.} By the arguments in 5.1, both the statement and the proof of Lemma 1.10 are correct if we restrict to a tamely ramified Galois extension.

\vskip2mm \noindent {\bf 5.3.} Lemma A.10 was not used in \cite{RTW1}. In the second paper \cite{RTW2}, we used it in Lemma 4.6 of \cite{RTW2}, a technical step in the proof of Proposition 4.5; therefore, both statements are proved only if the group ${\rm G}$ splits over a tamely ramified extension. Finally, the same restriction applies to Theorem 4.8 since the proof relies on Proposition 4.5.


\backmatter

\end{document}